\documentclass[pdflatex,sn-mathphys-num]{sn-jnl}

\usepackage[T1]{fontenc}
\usepackage[utf8]{inputenc}
\usepackage{amsmath,amssymb,amsthm,mathtools}
\usepackage{enumitem}
\usepackage{booktabs}
\usepackage{xcolor}
\hypersetup{
  hidelinks,
  bookmarksdepth=2,
  pdftitle={Markovian Lifts and Multifactor Approximations for Weighted Sub-fractional Brownian Motion},
  pdfauthor={Jose H. Ramirez-Gonzalez and Edison Cuba},
  pdfkeywords={Weighted sub-fractional Brownian motion, Markovian lifts, weighted Sobolev spaces, Young differential equations, multifactor approximation, strong convergence}
}
\allowdisplaybreaks[1]
\theoremstyle{thmstyleone}
\newtheorem{theorem}{Theorem}[section]
\newtheorem{proposition}[theorem]{Proposition}
\newtheorem{lemma}[theorem]{Lemma}
\newtheorem{corollary}[theorem]{Corollary}
\theoremstyle{thmstylethree}
\newtheorem{assumption}[theorem]{Assumption}
\theoremstyle{thmstyletwo}
\newtheorem{remark}[theorem]{Remark}
\newtheorem*{remark*}{Remark}
\theoremstyle{thmstyleone}

\newcommand{\R}{\mathbb{R}}
\newcommand{\Rplus}{\mathbb{R}_{+}}
\newcommand{\E}{\mathbb{E}}
\newcommand{\cE}{\mathcal{E}}
\newcommand{\cU}{\mathcal{U}}
\newcommand{\cL}{\mathcal{L}}
\newcommand{\cD}{\mathcal{D}}

\newcommand{\HS}{\mathrm{HS}}
\newcommand{\Tr}{\operatorname{Tr}}
\newcommand{\TV}{\operatorname{TV}}
\newcommand{\dd}{\,\mathrm{d}}
\newcommand{\drift}{\mathfrak{a}}
\newcommand{\norm}[1]{\left\lVert #1\right\rVert}
\newcommand{\abs}[1]{\left\lvert #1\right\rvert}
\newcommand{\ip}[2]{\left\langle #1,#2\right\rangle}
\newcommand{\Holder}[1]{C^{#1}}

\begin{document}
\raggedbottom

\title[Markovian Lifts and Multifactor Approximations]{Markovian Lifts and Multifactor Approximations for Weighted Sub-fractional Brownian Motion}

\author*[1]{\fnm{Jose H.} \sur{Ramirez-Gonzalez}}\email{hermenegildo.ramirez@usp.br}
\author[2]{\fnm{Edison} \sur{Cuba}}\email{edison.cuba@unsaac.edu.pe}

\affil[1]{\orgdiv{Institute of Mathematics and Statistics},
  \orgname{University of S\~ao Paulo (USP)},
  \orgaddress{\city{S\~ao Paulo}, \country{Brazil}}}
\affil[2]{\orgdiv{Departamento Acad\'emico de Matem\'aticas y Estad\'istica},
  \orgname{Universidad Nacional de San Antonio Abad del Cusco},
  \orgaddress{\city{Cusco}, \country{Per\'u}}}

\abstract{We construct a Hilbert-space Markovian lift of weighted sub-fractional
Brownian motion for $0<b<2$, including the logarithmic case $b=1$.
A Laplace factorization represents the process as the boundary value
of a forward curve. In the weighted Sobolev space with derivative
weight $r^\gamma/(1+r)^2$, $0<\gamma<1$, right translation is strongly
continuous and the noise operator is Hilbert--Schmidt if and only if
$b+\gamma>1$. The state norm bounds the total variation of initial
curves on compact intervals. For locally bounded deterministic temporal
weights, this bound allows the driven scalar equation to be restarted
from every lifted state. Under spatial Lipschitz assumptions, with
additional time and derivative regularity in the Young case, the
solution and forward curve form a time-inhomogeneous Markov process
with $C_b$-Feller transitions; both coefficients may have linear growth.
Discretization of the Laplace variable gives independent
Ornstein--Uhlenbeck factors coupled to the original process through
the same Gaussian noise. The nonlinear approximations converge
uniformly in time in every finite $L^p$. For uniform spectral
partitions with balanced cutoffs, the solution rate is
$M^{-(b-1)/3}$ for $1<b<2$, while the uniform driver rate is
$M^{-b/4}\sqrt{\log(2+M)}$ for $0<b\le1$.
Under stronger coefficient assumptions, we derive error bounds for
Young solutions and for the fully discrete scheme, retaining the
dependence on the spectral cutoff and specifying joint mesh conditions.}


\maketitle

\section{Introduction and main results}

Weighted sub-fractional Brownian motion is a Gaussian process whose
covariance combines a temporal weight with a memory kernel. Its current
value generally does not determine the conditional law of future
increments. We record the contribution of past Brownian noise to future
values in a forward curve, which provides a Markovian state for the
driver. Discretizing the same covariance representation gives finitely
many Ornstein--Uhlenbeck factors. We study the resulting approximations
of both the driver and the scalar equation it drives.

Fix $0<b<2$ and a non-negative Borel function
$f\in L^1_{\mathrm{loc}}(\Rplus)$.  We write $\zeta_{f,b}$ for the centered
Gaussian process with covariance
\begin{equation}\label{eq:covariance-intro}
 R_{f,b}(s,t)
 =\int_0^{s\wedge t}f(u)Q_b(s-u,t-u)\dd u,
 \qquad s,t\ge0,
\end{equation}
where
\begin{equation}\label{eq:Q-def-intro}
 Q_b(x,y)=
 \begin{cases}
 \dfrac{x^b+y^b-(x+y)^b}{1-b},&b\neq1,\\[2mm]
 (x+y)\log(x+y)-x\log x-y\log y,&b=1,
 \end{cases}
 \qquad x,y\ge0,
\end{equation}
with $0\log0:=0$. For $b\neq1$, this is the weighted sub-fractional
Brownian covariance in \cite[Section~2.1, Eq.~(4)]{RamirezSun2024}.
The power-weight family originates in occupation-time fluctuation
limits of age-dependent branching systems
\cite[Theorem~2.1(i), Eq.~(2.1)]{LopezMimbelaMurilloRamirez2024}.
The logarithmic case $b=1$ appears for power weights in
\cite[Theorem~2.1(ii), p.~596]{LopezMimbelaMurilloRamirez2024}.
Besides its origin in occupation-time fluctuation limits, this Gaussian
family has been used in animal movement modeling. Ram\'irez-Gonz\'alez,
Murillo-Salas, and Sun~\cite{RamirezMurilloSun2026} fitted its logarithmic
member to bat trajectories using longitude and latitude observations.
The same family is included among the Gaussian models considered in
the multidimensional telemetry analysis of
\cite{RamirezChaconMoragaSun2026}. Inference for stochastic differential equations driven by weighted
sub-fractional Brownian motion is studied in
\cite{RamirezNeuralInference2026}. The Euler approximation is used to
construct a trajectory likelihood that retains the dependence between
the Gaussian driving increments. Neural-network and radial-basis
representations are used to estimate the drift and diffusion coefficients
together with the normalized temporal weight and covariance exponent.

The regularity of the boundary process changes at the critical value
$b=1$.  If $f\in L^\infty_{\mathrm{loc}}(\Rplus)$ and $0<b<1$, then
$\zeta_{f,b}$ admits a modification whose paths are H\"older continuous
of every order
\[
 0<\eta<\frac{1+b}{2}.
\]
At $b=1$ the increment variance acquires a logarithmic correction, but
there is still a modification with H\"older continuous paths of every
order below one. For later
statements it is convenient to write
\[
 H_b^\circ:=
 \begin{cases}
 (1+b)/2,&0<b<1,\\
 1,&b=1,
 \end{cases}
\]
where the superscript indicates that the endpoint is not asserted.  Hence
the nonlinear equation is a Young equation for $0<b\le1$.  When $1<b<2$,
the boundary admits an absolutely continuous realization and the same
equation is interpreted in the Lebesgue--Stieltjes sense.

\paragraph{Relation to the literature.}
Laplace representations of fractional Gaussian processes lead to
infinite-dimensional Markovian states and finite-factor approximations;
see \cite{CarmonaCoutin1998,CarmonaCoutinMontseny2000}.
Affine representations of fractional processes are developed in
\cite{HarmsStefanovits2019}, and strong approximation rates in
\cite{Harms2021}. For stochastic Volterra equations, boundary
representations through first-order SPDEs are studied in
\cite{BenthDeteringKruhner2022}; weighted Sobolev realizations and
polynomially weighted shift spaces are treated in
\cite{Huber2024,BianchiBonaccorsiCanadasFriesen2026}.
Related results on affine Volterra processes and their Markovian lifts
appear in
\cite{AbiJaberLarssonPulido2019,CuchieroTeichmann2020,Hamaguchi2024}.

The nonlinear equations in
\cite{Huber2024,BianchiBonaccorsiCanadasFriesen2026} are formulated as
stochastic Volterra equations with Brownian It\^o integrals. In our
construction, the Gaussian forward curve is also an It\^o stochastic
evolution, while its boundary drives the scalar equation in the Young
or Lebesgue--Stieltjes sense. The main issue for the nonlinear equation
is to restart it from an arbitrary curve in the state space: such a
curve need not have the H\"older regularity of the Gaussian process.
The total-variation bound below provides the required control.
For the finite-factor approximations, we use the same Gaussian noise
and estimate the driver error in the path norms used by the solution
map. This gives convergence under the well-posedness assumptions;
quantitative bounds require the additional hypotheses specified below.

\paragraph{The state space and the Hilbert--Schmidt condition.}

For
\begin{equation}\label{eq:mu-intro}
 \mu_b(\dd\lambda)
 :=\frac{1}{\Gamma(2-b)}\lambda^{1-b}\dd\lambda,
 \qquad \lambda>0,
\end{equation}
and
\begin{equation}\label{eq:phi-intro}
 \phi_r(\lambda)
 :=\sqrt b\,\frac{1-e^{-\lambda r}}{\lambda},
 \qquad r\ge0,
\end{equation}
let $\cU_b:=L^2((0,\infty),\mu_b)$. Here $\lambda$ is the positive
relaxation rate in the Laplace representation; the term spectral below
refers to this variable. Proposition~\ref{prop:factorization}
proves
\[
 Q_b(x,y)=\langle\phi_x,\phi_y\rangle_{\cU_b}.
\]
For $0<\gamma<1$, we use the polynomial weight
\begin{equation}\label{eq:state-weight-intro}
 w_\gamma(r):=\frac{r^\gamma}{(1+r)^2},\qquad r>0,
\end{equation}
and the Hilbert space
\begin{equation}\label{eq:E-intro}
 \cE_\gamma
 :=\left\{\theta\in AC_{\mathrm{loc}}(\Rplus):
 \|\theta\|_\gamma^2
 :=|\theta(0)|^2+
 \int_0^\infty\frac{r^\gamma}{(1+r)^2}|\theta'(r)|^2\dd r<\infty\right\}.
\end{equation}
The value assigned to $w_\gamma$ at the single point $r=0$ is immaterial.
$r^\gamma$ controls the singularity of the feature derivative at the origin,
while $(1+r)^{-2}$ provides sufficient decay at infinity and is stable under
right translation.  Consider the feature rule
\begin{equation}\label{eq:G-intro}
 (G_bh)(r):=\langle\phi_r,h\rangle_{\cU_b},
 \qquad h\in\cU_b,\ r\ge0.
\end{equation}
This rule defines a Hilbert--Schmidt operator from $\cU_b$ to
$\cE_\gamma$ precisely when
\begin{equation}\label{eq:key-condition-intro}
 b+\gamma>1.
\end{equation}
\[
 \|G_b\|_{\HS(\cU_b,\cE_\gamma)}^2
 =b2^{b-2}\Gamma(b+\gamma-1)\Gamma(3-b-\gamma).
\]
The parameters of the lifted equation satisfy
\begin{equation}\label{eq:standing-gamma-range}
 0<\gamma<1,\qquad b+\gamma>1.
\end{equation}
Thus $1-b<\gamma<1$ for $0<b<1$ and $0<\gamma<1$ for $1\le b<2$.
Lemma~\ref{lem:G-HS} also covers $\gamma=0$ when $1<b<2$;
we use positive $\gamma$ to state the construction in one parameter range.

Write $(S_t\theta)(r)=\theta(t+r)$ for right translation and let $A$
be its generator. Lemma~\ref{lem:state-shift} proves that $(S_t)_{t\ge0}$
is strongly continuous on $\cE_\gamma$ and that boundary evaluation
$\delta_0\theta=\theta(0)$ is continuous. Theorem~\ref{thm:unified-lift}
then constructs the mild solution of
\[
 \dd\Theta_t=A\Theta_t\dd t+\sqrt{f(t)}G_b\dd W_t,
 \qquad \Theta_0=0,
\]
where $W$ is cylindrical Brownian motion on $\cU_b$.
Its boundary is $\Theta_t(0)=\zeta_{t,f,b}$, while
\[
 \Theta_t(r)=\E[\zeta_{t+r,f,b}\mid\mathcal F_t].
\]
Here $\mathcal F_t$ is the filtration carrying the driving Brownian motions.
The curve records the contribution of past noise to future values;
the remaining contribution uses independent future Brownian increments.
This decomposition gives the Markov property.

\paragraph{The nonlinear equation.}
Fix $T>0$ and consider
\begin{equation}\label{eq:main-equation-intro}
 X_t=x_0+\int_0^t\drift(X_s,s)\dd s
       +\int_0^t\sigma(X_s,s)\dd\zeta_{s,f,b},
 \qquad 0\le t\le T,
\end{equation}
where $x_0\in\R$ is deterministic and
$\drift,\sigma:\R\times\Rplus\to\R$ are the drift and diffusion coefficients.

Assume additionally that $f\in L^\infty_{\mathrm{loc}}(\Rplus)$.
Under Assumption~\ref{ass:coefficients-mf},
Theorem~\ref{thm:nonlinear-lift} gives the Markov property for
$(X,\Theta)$ and $C_b$-Feller continuity of its time-dependent transitions. In the Young regime this assumption
allows linear growth of both coefficients. The drift and diffusion are
continuous and globally Lipschitz in space; the diffusion is Lipschitz
in time and has a bounded, uniformly $\alpha$-H\"older spatial derivative, with
$\alpha>(1-b)/(1+b)$ for $b<1$ and $\alpha>0$ for $b=1$.
For $b>1$, only continuity and uniform global spatial Lipschitz continuity
of both coefficients are required.
Table~\ref{tab:b-gamma-ranges} summarizes the parameter ranges for this
Hilbert--Schmidt realization and the corresponding pathwise integrals.

\begin{table}[!htbp]
\centering
\caption{State exponents and pathwise interpretations used in the lifted theory.}
\label{tab:b-gamma-ranges}
\begin{tabular}{@{}lll@{}}
\toprule
Range of $b$ & State exponent $\gamma$ & Pathwise integral \\
\midrule
$0<b<1$ & $1-b<\gamma<1$ & Young \\
$b=1$ & $0<\gamma<1$ & Young \\
$1<b<2$ & $0<\gamma<1$ & Lebesgue--Stieltjes \\
\bottomrule
\end{tabular}
\end{table}

\paragraph{Approximation and coefficient regularity.}
Approximating the Laplace features by step functions in the spectral
variable produces independent Ornstein--Uhlenbeck factors. The finite
driver is the time integral of their weighted sum; the Markovian state
for the nonlinear approximation consists of the scalar solution and
these factors. Write $M$ for the number of factors. The spectral rate
bounds also allow nonuniform partitions of $[0,\Lambda]$ with maximum
cell width $\delta_\lambda\leq C_{\mathrm{mesh}}\Lambda/M$ and arbitrary
nodes in each cell, where $C_{\mathrm{mesh}}$ is independent of $M$
and $\Lambda$.
For a uniform partition of $[0,\Lambda]$,
Theorem~\ref{thm:finite-variation-rate}
gives the strong solution rate $M^{-(b-1)/3}$ when $1<b<2$ and
$\Lambda=M^{2/3}$. In the Young regime,
Theorem~\ref{thm:young-driver-rate} gives the uniform driver rate
$M^{-b/4}\sqrt{\log(2+M)}$ for $\Lambda=M^{1/2}$.
The H\"older estimate
and quantitative stability yield the distinct solution bound of
Theorem~\ref{thm:young-X-rate}:
\[
 \left(\E\sup_{t\leq T}|X_t-X_t^{M,\Lambda}|^p\right)^{1/p}
 \leq C\left(M^{-b/4}\sqrt{\log(2+M)}\right)^{1-\eta/\eta'},
 \qquad \Lambda=M^{1/2},
\]
for $p\geq2$ and $1/2<\eta<\eta'<H_b^\circ$, under
Assumption~\ref{ass:young-quantitative} and the spectral mesh condition
of that theorem. Proposition~\ref{prop:ou-weak-coefficients}
establishes qualitative convergence in every finite $L^p$, uniformly
in time, under the weaker Markov--Feller assumptions.

For $0<b\le1$, the coefficient assumptions form the hierarchy
\[
 \text{Assumption~\ref{ass:coefficients}}
 \ \Longrightarrow\
 \text{Assumption~\ref{ass:young-quantitative}}
 \ \Longrightarrow\
 \text{Assumption~\ref{ass:coefficients-mf}},
\]
where the last implication holds with $\alpha=1$.
Assumption~\ref{ass:coefficients-mf} gives the Markov--Feller property
and qualitative convergence. Assumption~\ref{ass:young-quantitative}
also bounds the diffusion and imposes joint spatial and temporal
Lipschitz continuity on its derivative; these bounds are used for
quantitative stability and the spectral solution rate.
Assumption~\ref{ass:coefficients} adds a bounded drift, uniformly
Lipschitz in time, for the Euler estimate and the fully discrete bound.
For $1<b<2$, the spectral solution rate in
Theorem~\ref{thm:finite-variation-rate} already holds under
Assumption~\ref{ass:coefficients-mf}, whereas the fully discrete bound
in Theorem~\ref{thm:joint-fully-discrete} assumes bounded coefficients
jointly Lipschitz in space and time.

Theorem~\ref{thm:joint-fully-discrete} bounds separately the spectral
error, the factors' time-discretization error, and the scalar Euler error.
The dependence on the spectral cutoff remains explicit, since increasing
that cutoff requires finer time resolution. Corollary~\ref{cor:joint-meshes}
gives sufficient conditions on the spectral and time meshes to preserve
the spectral convergence rate.

\paragraph{Organization.}
Section~2 develops the Laplace factorization, the polynomially weighted
state space, and the forward-curve lift.  Section~3 proves the
path regularity, the boundary-slope criterion, and the Markov--Feller lift of
\eqref{eq:main-equation-intro}.  Section~4 establishes the finite-factor and
fully discrete convergence bounds.
\paragraph{Notation and conventions.}\label{par:notation}
All stochastic processes are defined on a complete filtered probability space
$(\Omega,\mathcal F,(\mathcal F_t)_{t\geq0},\mathbb P)$ satisfying the usual
conditions, and $\E$ denotes expectation with respect to $\mathbb P$.  The symbol $\Gamma$ denotes the Euler gamma function, and $\Tr$ denotes the trace of a trace-class operator.  We set
$\mathbb N:=\{1,2,\ldots\}$ and write $\Rplus=[0,\infty)$ and $AC_{\mathrm{loc}}(\Rplus)$ for the locally absolutely
continuous functions on $\Rplus$.  If $E$ is a metric
space, $C_b(E)$ denotes the bounded continuous real-valued functions on $E$.
For normed spaces $U$ and $E$, $\cL(U,E)$ denotes the bounded linear
operators from $U$ to $E$, with operator norm
$\|T\|_{\mathrm{op}}:=\sup_{\|u\|_U\leq1}\|Tu\|_E$.
For a function $g:\R\to\R$, write
$\operatorname{Lip}(g):=\sup_{x\ne y}|g(x)-g(y)|/|x-y|$.
If $U$ and $E$ are Hilbert spaces, an operator
$T\in\cL(U,E)$ is Hilbert--Schmidt when, for an orthonormal basis $(e_k)$
of $U$,
\[
 \sum_{k=1}^\infty\|Te_k\|_E^2<\infty.
\]
This condition and the value of the sum are independent of the chosen basis;
we denote the corresponding space by $\cL_2(U,E)$ and its norm by
$\|T\|_{\HS(U,E)}$.  We write $T^*$ for the
Hilbert-space adjoint of $T$.  The domain of an operator $T$ is denoted by
$\cD(T)$, and
$\langle\cdot,\cdot\rangle_U$ and $\|\cdot\|_U$ denote the inner product and
norm on $U$.  For a real random variable $Z$ and $p\geq1$,
$\|Z\|_{L^p(\Omega)}:=(\E|Z|^p)^{1/p}$.  For $x\in C([0,T])$ and $\eta\in(0,1)$, set
\[
 \|x\|_{\infty;[0,T]}:=\sup_{0\le t\le T}|x(t)|,
 \qquad
 [x]_{\eta;[0,T]}
 :=\sup_{0\le s<t\le T}\frac{|x(t)-x(s)|}{|t-s|^\eta}.
\]
We write
$\Holder{\eta}([0,T]):=\{x\in C([0,T]):[x]_{\eta;[0,T]}<\infty\}$,
equipped with
\[
 \|x\|_{\Holder{\eta}([0,T])}
 :=\|x\|_{\infty;[0,T]}+[x]_{\eta;[0,T]}.
\]
For $x\in C([0,T])$ and $p\geq1$, define
\[
 \|x\|_{p\text{-var};[0,T]}
 :=\left(\sup_{\Pi}\sum_{i=0}^{m-1}
 |x(t_{i+1})-x(t_i)|^p\right)^{1/p},
\]
where the supremum is over all finite partitions
$\Pi=\{0=t_0<\cdots<t_m=T\}$, and set
\[
 \TV_{[0,T]}(x)
 :=\sup_{\Pi}\sum_{i=0}^{m-1}|x(t_{i+1})-x(t_i)|.
\]
For a finite partition $\pi$ of an interval, its mesh width is
$|\pi|:=\max_{[u,v]\in\pi}(v-u)$.
For $f\in L^\infty(0,T)$, write
\[
 \|f\|_{\infty,T}:=\operatorname*{ess\,sup}_{0\le t\le T}|f(t)|.
\]
For processes $X$ and $Y$, $X\stackrel{\mathrm{f.d.d.}}{=}Y$ means that
$(X_{t_1},\ldots,X_{t_m})\stackrel{d}{=}(Y_{t_1},\ldots,Y_{t_m})$ for
every $m\in\mathbb N$ and $t_1,\ldots,t_m\geq0$.  The letter $C$ denotes a finite deterministic constant whose
value may change from line to line.  A subscript on $C$ records the parameters
on which the constant may depend.  A time-inhomogeneous transition
evolution $(P_{s,t})_{s\le t}$ on $E$ is called Feller when
$P_{s,t}C_b(E)\subseteq C_b(E)$ for every $s\le t$.
We use the term $C_b$-Feller in this sense; no strong continuity in the
supremum norm on $C_b$ is asserted.
Local integrability and local absolute continuity include the origin.
Unless explicitly stated otherwise, every weight $f$ appearing below is Borel
measurable and non-negative; the displayed $L^1_{\mathrm{loc}}$ or
$L^\infty_{\mathrm{loc}}$ assumptions are additional integrability
requirements.

\section{Covariance factorization and the Markovian lift}

We first factor the covariance into Laplace features.  We then choose a
Sobolev space on which right translation is strongly continuous and the
feature operator is Hilbert--Schmidt.  The resulting stochastic convolution
has the Gaussian driver as its boundary value and the Markov property as a
curve-valued process.

\subsection{The covariance features}

For $0<b<2$, let $\mu_b$ and $\phi_r$ be given by
\eqref{eq:mu-intro}--\eqref{eq:phi-intro}.  Set
\begin{equation*}
 \kappa_b:=
 \begin{cases}
 \dfrac{2-2^b}{1-b},&b\neq1,\\[2mm]
 2\log2,&b=1.
 \end{cases}
\end{equation*}

For $b\ne1$, the Laplace identity below is given in
\cite[proof of Theorem~1, Eq.~(6)]{RamirezSun2024}. We formulate it in the
spectral Hilbert space used below, with the same normalization at $b=1$.
These features define both the forward curve and its finite-factor
approximations.  Their linear independence will also rule out exact
representations by a fixed finite family of convolution kernels.

\begin{proposition}\label{prop:factorization}
For every $0<b<2$ and $x,y\geq0$,
\begin{equation}\label{eq:factorization}
 Q_b(x,y)=\ip{\phi_x}{\phi_y}_{\cU_b}
 =\int_0^\infty\phi_x(\lambda)\phi_y(\lambda)\mu_b(\dd\lambda),
\end{equation}
where $Q_b$ is defined by \eqref{eq:Q-def-intro}.  In particular,
\begin{equation}\label{eq:phi-norm}
 \norm{\phi_r}_{\cU_b}^2=Q_b(r,r)=\kappa_b r^b.
\end{equation}
The kernel $Q_b$ has infinite rank on every nontrivial interval.  More
precisely, for any distinct $x_1,\ldots,x_n>0$, the Gram matrix
$\bigl(Q_b(x_i,x_j)\bigr)_{i,j=1}^n$ is strictly positive definite.
\end{proposition}

\begin{proof}
By definition,
\[
 \ip{\phi_x}{\phi_y}_{\cU_b}
 =\frac{b}{\Gamma(2-b)}
 \int_0^\infty
 (1-e^{-\lambda x})(1-e^{-\lambda y})\lambda^{-1-b}\dd\lambda.
\]
For $x,y>0$, the integrand is of order $\lambda^{1-b}$ at zero and
$\lambda^{-1-b}$ at infinity, so it is integrable for $0<b<2$.
Using
\[
 (1-e^{-\lambda x})(1-e^{-\lambda y})
 =\lambda^2\int_0^x\int_0^y e^{-\lambda(v+w)}\dd w\dd v
\]
and Tonelli's theorem, we obtain
\begin{align*}
 \ip{\phi_x}{\phi_y}_{\cU_b}
 &=\int_0^x\int_0^y
 \frac{b}{\Gamma(2-b)}
 \int_0^\infty e^{-\lambda(v+w)}\lambda^{1-b}\dd\lambda\dd w\dd v\\
 &=\int_0^x\int_0^y b(v+w)^{b-2}\dd w\dd v.
\end{align*}
The singularity at $(v,w)=(0,0)$ is integrable: the length of a slice
$v+w=s$ in the rectangle is at most $s$, so the integral near zero is
bounded by a constant times $\int_0^{x+y}s^{b-1}\dd s$.
For $b\ne1$, direct integration gives
\[
 \int_0^x\int_0^y b(v+w)^{b-2}\dd w\dd v
 =\frac{(x+y)^b-x^b-y^b}{b-1}=Q_b(x,y).
\]
For $b=1$, the same integral is
\[
 \int_0^x\log\frac{v+y}{v}\dd v
 =(x+y)\log(x+y)-x\log x-y\log y=Q_1(x,y).
\]
If either coordinate is zero, both sides vanish.  This proves
\eqref{eq:factorization} in all cases, and setting $x=y=r$ gives
\eqref{eq:phi-norm}.

To prove infinite rank, fix distinct positive numbers
$x_1,\ldots,x_n$ and suppose that
$\sum_{j=1}^n c_j\phi_{x_j}=0$ in $\cU_b$.  Since the left-hand side is a
continuous function of $\lambda>0$ and $\mu_b$ has full support, it vanishes
for every $\lambda>0$.  After multiplying by $\lambda/\sqrt b$ we obtain
\[
 \sum_{j=1}^n c_j(1-e^{-\lambda x_j})=0,
 \qquad \lambda>0.
\]
Letting $\lambda\to\infty$ gives $\sum_jc_j=0$, and hence
$\sum_jc_je^{-\lambda x_j}=0$ for every $\lambda>0$.  The exponentials with
distinct exponents are linearly independent: ordering
$x_1<\cdots<x_n$, multiplying by $e^{\lambda x_1}$, and letting
$\lambda\to\infty$ gives $c_1=0$; iteration gives
$c_1=\cdots=c_n=0$.  Thus every finite collection of distinct features is
linearly independent.  For a nonzero vector $c=(c_1,\ldots,c_n)$,
\[
 \sum_{i,j=1}^n c_ic_jQ_b(x_i,x_j)
 =\norm{\sum_{j=1}^n c_j\phi_{x_j}}_{\cU_b}^2>0,
\]
so the corresponding Gram matrix is strictly positive definite.  Since every
nontrivial interval contains arbitrarily large finite collections of distinct
points, the Gram kernel has infinite rank there.
\end{proof}

\begin{corollary}
\label{cor:no-finite-factor}
Fix $0<b<2$.  For every $N\in\mathbb N$, there are no functions
$k_{b,1},\ldots,k_{b,N}:\Rplus\to\R$ such that
\begin{equation}\label{eq:no-finite-factor-kernel}
 Q_b(x,y)=\sum_{j=1}^N k_{b,j}(x)k_{b,j}(y),
 \qquad x,y\geq0.
\end{equation}
Consequently, there is no representation valid for every non-negative,
locally bounded weight $f$ of the form
\begin{equation}\label{eq:no-finite-factor-process}
 \zeta_{t,f,b}\stackrel{\mathrm{f.d.d.}}{=}
 \sum_{j=1}^N\int_0^t\sqrt{f(u)}k_{b,j}(t-u)\dd B_u^j,
\end{equation}
where $B^1,\ldots,B^N$ are independent real Brownian motions and the kernels
$k_{b,j}$ are Borel measurable, belong to
$L^2_{\mathrm{loc}}(\Rplus)$, and are independent of $f$.
\end{corollary}

\begin{proof}
If \eqref{eq:no-finite-factor-kernel} held, every Gram matrix of $Q_b$ would
have rank at most $N$.  Proposition~\ref{prop:factorization}, applied to
$N+1$ distinct positive points, gives instead a strictly positive definite
$(N+1)\times(N+1)$ Gram matrix, a contradiction.  If
\eqref{eq:no-finite-factor-process} held for every non-negative locally
bounded $f$, comparison of covariances would give, for every $s,t>0$,
\begin{equation}\label{eq:localized-factor-equality}
 \int_0^{s\wedge t} f(u)D_{s,t}(u)\dd u=0,
\end{equation}
where
\[
 D_{s,t}(u):=Q_b(s-u,t-u)
 -\sum_{j=1}^Nk_{b,j}(s-u)k_{b,j}(t-u).
\]
The local $L^2$ assumption and Cauchy--Schwarz show that
$D_{s,t}\in L^1(0,s\wedge t)$.  Taking in
\eqref{eq:localized-factor-equality} arbitrary non-negative bounded functions
with compact support in $(0,s\wedge t)$ yields
$D_{s,t}=0$ almost everywhere.

Fix $d\geq0$ and $L>0$, and take $s=L$ and $t=L+d$.  The change of variables
$x=L-u$ gives
\begin{equation}\label{eq:factorization-on-lines}
 Q_b(x,x+d)=\sum_{j=1}^Nk_{b,j}(x)k_{b,j}(x+d)
 \quad\text{for almost every }x\in(0,L).
\end{equation}
Since $L$ is arbitrary, this holds for almost every $x>0$ for each fixed
$d\geq0$.  Fubini's theorem in the variables $(x,d)$, followed by the
change of variables $(x,d)\mapsto(x,x+d)$ and symmetry, shows that
\begin{equation}\label{eq:factorization-ae-plane}
 Q_b(x,y)=\sum_{j=1}^Nk_{b,j}(x)k_{b,j}(y)
 \quad\text{for almost every }(x,y)\in(0,\infty)^2.
\end{equation}
Equation \eqref{eq:factorization-on-lines} with $d=0$ also gives the diagonal
identity for almost every $x>0$.

Let $I$ be a bounded nontrivial interval in $(0,\infty)$ and choose
$N+1$ independent points uniformly distributed on $I$.  With
probability one they are distinct, avoid the exceptional diagonal set, and
every pair avoids the null set in \eqref{eq:factorization-ae-plane}.
Consequently, there exist distinct deterministic points
$x_1,\ldots,x_{N+1}\in I$ for which
\[
 \bigl(Q_b(x_i,x_j)\bigr)_{i,j=1}^{N+1}
 =\left(\sum_{\ell=1}^N
 k_{b,\ell}(x_i)k_{b,\ell}(x_j)\right)_{i,j=1}^{N+1}.
\]
The matrix on the right has rank at most $N$, whereas the matrix on the left
is strictly positive definite by Proposition~\ref{prop:factorization}.  This
contradiction proves the second assertion.
\end{proof}

The corollary applies to a fixed finite family of convolution kernels that
works for every non-negative locally bounded weight.  It does not rule out other
finite-dimensional states or kernels chosen for a fixed weight.

\subsection{Polynomially weighted Sobolev spaces}

For $0\leq\gamma<1$, let $w_\gamma$ and $\cE_\gamma$ be given by
\eqref{eq:state-weight-intro}--\eqref{eq:E-intro}.  Put
\begin{equation}\label{eq:J-gamma}
 J_\gamma(T)
 :=\int_0^T\frac{1}{w_\gamma(r)}\dd r
 =\frac{T^{1-\gamma}}{1-\gamma}
  +\frac{2T^{2-\gamma}}{2-\gamma}
  +\frac{T^{3-\gamma}}{3-\gamma}<\infty.
\end{equation}

\begin{lemma}\label{lem:state-shift}
Let $0\leq\gamma<1$.  The space $\cE_\gamma$ is a separable Hilbert space.
Every $\theta\in\cE_\gamma$ has finite variation on compact intervals and,
for every $T>0$,
\begin{equation}\label{eq:state-TV}
 \TV_{[0,T]}(\theta)\leq J_\gamma(T)^{1/2}\|\theta\|_\gamma.
\end{equation}
Boundary evaluation $\delta_0\theta=\theta(0)$ is continuous.  Moreover,
$(S_t\theta)(r):=\theta(t+r)$ defines a strongly continuous semigroup on
$\cE_\gamma$.  Its generator is
\begin{equation}\label{eq:shift-generator}
 A\theta=\theta',\qquad
 \cD(A)=\{\theta\in\cE_\gamma:\theta'\in\cE_\gamma\}.
\end{equation}
Here $\theta'\in\cE_\gamma$ means that the almost-everywhere derivative
has a locally absolutely continuous representative belonging to that space.
For every $t\geq0$,
\begin{equation}\label{eq:shift-bound-explicit}
 \|S_t\|_{\cL(\cE_\gamma)}^2
 \leq\bigl(1+J_\gamma(t)^{1/2}\bigr)^2+(1+t)^2.
\end{equation}
In particular, for each $T>0$ there is $C_{\gamma,T}<\infty$ such that
$\sup_{0\le t\le T}\|S_t\|_{\cL(\cE_\gamma)}\le C_{\gamma,T}$.
\end{lemma}

\begin{proof}
Set
\[
 H_\gamma:=L^2\!\left(\Rplus,\frac{r^\gamma}{(1+r)^2}\dd r\right).
\]
With the inner product
\[
 \langle\theta,\vartheta\rangle_\gamma
 :=\theta(0)\vartheta(0)
   +\int_0^\infty w_\gamma(r)\theta'(r)\vartheta'(r)\dd r,
\]
the map $\mathcal J\theta:=(\theta(0),\theta')$ is an isometry from
$\cE_\gamma$ into $\R\times H_\gamma$.  Conversely, if
$(a,h)\in\R\times H_\gamma$, then
\[
 \theta(r):=a+\int_0^r h(q)\dd q
\]
is locally absolutely continuous.  Indeed, by Cauchy--Schwarz and
\eqref{eq:J-gamma},
\[
 \int_0^T|h(r)|\dd r
 \leq J_\gamma(T)^{1/2}
 \left(\int_0^T\frac{r^\gamma}{(1+r)^2}|h(r)|^2\dd r\right)^{1/2}<\infty.
\]
Thus $\mathcal J$ is onto.  Since $\R\times H_\gamma$ is a separable
Hilbert space, so is $\cE_\gamma$.  An absolutely continuous curve has
total variation $\int_0^T|\theta'(r)|\dd r$, so the same estimate gives
\eqref{eq:state-TV}.  Also $|\delta_0\theta|\leq\|\theta\|_\gamma$.

For $u>t\ge0$,
\[
 \frac{w_\gamma(u-t)}{w_\gamma(u)}
 =\left(\frac{u-t}{u}\right)^\gamma
  \left(\frac{1+u}{1+u-t}\right)^2
 \le(1+t)^2,
\]
where the first factor is at most one and
$(1+u)/(1+u-t)=1+t/(1+u-t)\le1+t$.  Therefore
\begin{align*}
 \int_0^\infty w_\gamma(r)|(S_t\theta)'(r)|^2\dd r
 &=\int_t^\infty w_\gamma(u-t)|\theta'(u)|^2\dd u\\
 &\le(1+t)^2\int_t^\infty w_\gamma(u)|\theta'(u)|^2\dd u.
\end{align*}
Moreover,
\[
 |\theta(t)|
 \le |\theta(0)|+\int_0^t|\theta'(r)|\dd r
 \le\bigl(1+J_\gamma(t)^{1/2}\bigr)\|\theta\|_\gamma.
\]
Combining the last two estimates gives \eqref{eq:shift-bound-explicit}.
In particular, every $S_t$ is bounded and the family is locally uniformly
bounded in operator norm.

To prove strong continuity, let $\tau_th(r):=h(t+r)$ on $H_\gamma$.  If
$h\in C_c(\Rplus)$, then $\tau_th\to h$ pointwise and uniformly on a fixed
compact interval containing the supports for all sufficiently small $t$;
the weight is locally integrable, so dominated convergence gives
$\|\tau_th-h\|_{H_\gamma}\to0$.  Since $C_c(\Rplus)$ is dense in
$H_\gamma$ and the preceding shift bound is locally uniform, the same
conclusion holds for every $h\in H_\gamma$.  Under $\mathcal J$,
\begin{equation}\label{eq:shift-under-isometry}
 \mathcal J S_t\mathcal J^{-1}(a,h)
 =\left(a+\int_0^t h(s)\dd s,\tau_th\right).
\end{equation}
The first component converges to $a$ by local integrability of $h$, and the
second converges to $h$ in $H_\gamma$; hence $(S_t)_{t\ge0}$ is strongly
continuous.

It remains to identify the generator.  Suppose first that
$\theta'\in\cE_\gamma$ and put $h=\theta'$.  Then $h$ is locally absolutely
continuous, $h'\in H_\gamma$, and $h(0)$ is well defined.  From
\eqref{eq:shift-under-isometry},
\[
 \mathcal J\frac{S_t\theta-\theta}{t}
 =\left(\frac1t\int_0^t h(s)\dd s,
         \frac{\tau_th-h}{t}\right).
\]
The first component converges to $h(0)$.  For the second, the fundamental
theorem of calculus gives
\[
 \frac{\tau_th-h}{t}-h'
 =\int_0^1(\tau_{vt}h'-h')\dd v.
\]
In particular,
\[
 \left\|\frac{\tau_th-h}{t}-h'\right\|_{H_\gamma}
 \leq\sup_{0\leq s\leq t}\|\tau_sh'-h'\|_{H_\gamma}\longrightarrow0
\]
by strong continuity at zero.  Thus
$(S_t\theta-\theta)/t\to\theta'$ in $\cE_\gamma$.

Conversely, suppose that $(S_t\theta-\theta)/t$ converges in $\cE_\gamma$
as $t\downarrow0$, and put $h=\theta'$.  Then
$(\tau_th-h)/t\to g$ in $H_\gamma$ for some $g$.  On every compact
subinterval of $(0,\infty)$ the weight is bounded above and below by positive
constants, so this convergence also holds in $L^2_{\mathrm{loc}}$.  For
$\varphi\in C_c^\infty((0,\infty))$, a change of variables yields
\[
 \int_0^\infty g(r)\varphi(r)\dd r
 =\lim_{t\downarrow0}\int_0^\infty
   \frac{h(r+t)-h(r)}{t}\varphi(r)\dd r
 =-\int_0^\infty h(r)\varphi'(r)\dd r.
\]
Hence $g$ is the weak derivative of $h$ on $(0,\infty)$.  Since both
$h$ and $g$ belong to $H_\gamma$, Cauchy--Schwarz with
$w_\gamma^{-1}$ and \eqref{eq:J-gamma} gives
\[
 \int_0^T(|h(r)|+|g(r)|)\dd r
 \leq J_\gamma(T)^{1/2}
       (\|h\|_{H_\gamma}+\|g\|_{H_\gamma})<\infty.
\]
To see the behavior at the origin explicitly, fix $r_0>0$ at which an
absolutely continuous interior representative of $h$ is defined.  On
$(0,r_0]$ this representative satisfies
$h(r)=h(r_0)-\int_r^{r_0}g(q)\dd q$.  The integrability of $g$ gives a
finite limit as $r\downarrow0$ and an absolutely continuous extension to
$[0,r_0]$.  Thus $h\in AC_{\mathrm{loc}}(\Rplus)$, and
$|h(0)|^2+\|g\|_{H_\gamma}^2<\infty$ gives $h\in\cE_\gamma$.
The boundary component
$t^{-1}\int_0^t h(s)\dd s$ converges to $h(0)$, so the generator limit is
precisely the element $h=\theta'$.  Hence both inclusions in \eqref{eq:shift-generator} hold.
\end{proof}

For $h\in\cU_b$, define
\begin{equation}\label{eq:G-def}
 (G_bh)(r):=\langle\phi_r,h\rangle_{\cU_b}.
\end{equation}
For $r>0$,
\[
 \phi_r'(\lambda)=\sqrt b\,e^{-\lambda r},
\]
and hence
\begin{equation}\label{eq:phi-derivative-norm}
 \|\phi_r'\|_{\cU_b}^2
 =\frac{b}{\Gamma(2-b)}
   \int_0^\infty e^{-2\lambda r}\lambda^{1-b}\dd\lambda
 =b(2r)^{b-2}.
\end{equation}
For every $T>0$ and $h\in\cU_b$, it follows that
\begin{align*}
 \int_0^T|\langle\phi_r',h\rangle_{\cU_b}|\dd r
 &\le\|h\|_{\cU_b}\int_0^T\|\phi_r'\|_{\cU_b}\dd r\\
 &=\|h\|_{\cU_b}\sqrt b\,2^{(b-2)/2}
   \int_0^T r^{(b-2)/2}\dd r<\infty,
\end{align*}
because $b>0$.  Also $\phi_r=\int_0^r\phi_q'\dd q$ as a Bochner integral
in $\cU_b$.  Consequently $G_bh\in AC_{\mathrm{loc}}(\Rplus)$ and
\[
 (G_bh)'(r)=\langle\phi_r',h\rangle_{\cU_b}
 \quad\text{for almost every }r>0.
\]

\begin{lemma}\label{lem:G-HS}
Let $0<b<2$ and $0\le\gamma<1$.  The rule \eqref{eq:G-def} defines an
operator $G_b\in\cL_2(\cU_b,\cE_\gamma)$ if and only if
\begin{equation}\label{eq:G-HS-condition}
 b+\gamma>1.
\end{equation}
In that case,
\begin{equation}\label{eq:G-HS-norm}
 \|G_b\|_{\HS(\cU_b,\cE_\gamma)}^2
 =b2^{b-2}\Gamma(b+\gamma-1)\Gamma(3-b-\gamma).
\end{equation}
\end{lemma}

\begin{proof}
If $b+\gamma>1$, Cauchy--Schwarz and \eqref{eq:phi-derivative-norm} give
\[
 \|G_bh\|_\gamma^2
 \leq b2^{b-2}\|h\|_{\cU_b}^2
 \int_0^\infty\frac{r^{b+\gamma-2}}{(1+r)^2}\dd r<\infty,
\]
so \eqref{eq:G-def} maps $\cU_b$ continuously into $\cE_\gamma$.
Let $(e_n)$ be an orthonormal basis of $\cU_b$.  Since $\phi_0=0$,
Parseval's identity, Tonelli's theorem, and
\eqref{eq:phi-derivative-norm} give
\begin{align*}
 \sum_{n=1}^\infty\|G_be_n\|_\gamma^2
 &=\int_0^\infty\frac{r^\gamma}{(1+r)^2}
   \sum_{n=1}^\infty|\langle\phi_r',e_n\rangle_{\cU_b}|^2\dd r\\
 &=b2^{b-2}\int_0^\infty
   \frac{r^{b+\gamma-2}}{(1+r)^2}\dd r.
\end{align*}
Near zero the integrand is equivalent to $r^{b+\gamma-2}$ and is integrable
exactly when $b+\gamma>1$.  At infinity it is equivalent to
$r^{b+\gamma-4}$, which is integrable automatically because
$b<2$ and $\gamma<1$.  If $b+\gamma>1$, set
$a=b+\gamma-1\in(0,2)$.  Then
\[
 \int_0^\infty\frac{r^{a-1}}{(1+r)^2}\dd r
 =\int_0^1v^{a-1}(1-v)^{1-a}\dd v
 =\Gamma(a)\Gamma(2-a),
\]
where $v=r/(1+r)$ and the last equality is Euler's beta integral.
This gives \eqref{eq:G-HS-norm}.  Conversely, if \eqref{eq:G-def}
defined a Hilbert--Schmidt map into $\cE_\gamma$, the same Parseval--Tonelli
calculation would have a finite left-hand side.  This is impossible when
$b+\gamma\le1$, because the integral diverges at zero.
\end{proof}

\begin{remark}
For $b=1$, $\mu_1(\dd\lambda)=\dd\lambda$ and
$\phi'_r(\lambda)=e^{-\lambda r}$, so
\[
 \|\phi'_r\|_{\cU_1}^2
 =\int_0^\infty e^{-2\lambda r}\dd\lambda=\frac1{2r}.
\]
The condition in Lemma~\ref{lem:G-HS} reduces to $\gamma>0$.
For $0<\gamma<1$, formula~\eqref{eq:G-HS-norm} gives
\[
 \|G_1\|_{\HS(\cU_1,\cE_\gamma)}^2
 =\frac12\Gamma(\gamma)\Gamma(2-\gamma)
 =\frac{\pi(1-\gamma)}{2\sin(\pi\gamma)}<\infty.
\]
At $\gamma=0$, the Parseval--Tonelli integral
$\int_0^\infty r^{-1}(1+r)^{-2}\dd r$ diverges logarithmically at zero,
since $r^{-1}(1+r)^{-2}\sim r^{-1}$ as $r\downarrow0$.
Consequently, $G_1\notin\cL_2(\cU_1,\cE_0)$.
Thus $\gamma=0$ is admissible for the Hilbert--Schmidt property only
when $1<b<2$, whereas positive $\gamma$ also covers the logarithmic case.
\end{remark}

\begin{remark}\label{rem:weight-choice}
More generally, the derivative weight
$w_{\gamma,m}(r)=r^\gamma/(1+r)^m$ yields the Hilbert--Schmidt integral
$\int_0^\infty r^{b+\gamma-2}(1+r)^{-m}\dd r$.
Its finiteness is equivalent to $1<b+\gamma<m+1$.
The choice $m=2$ makes the condition at infinity automatic for every
$0<b<2$ and $0\le\gamma<1$, leaving only the restriction at zero.
\end{remark}

\subsection{The transport lift}
\label{subsec:transport-lift}

The term transport refers to right translation in the horizon variable $r$,
whose generator is $A\theta=\theta'$.

Fix an orthonormal basis $(e_k)_{k\ge1}$ of $\cU_b$ and independent
standard $(\mathcal F_t)$-Brownian motions $(\beta_k)_{k\ge1}$.  The associated
cylindrical Brownian motion $W$ on $\cU_b$ is the centered Gaussian family
\[
 W_t(h):=\sum_{k=1}^\infty
 \langle h,e_k\rangle_{\cU_b}\,\beta_k(t),
 \qquad h\in\cU_b,
\]
where the series converges in $L^2(\Omega)$; in particular,
$\E[W_t(h)W_s(g)]=(t\wedge s)\langle h,g\rangle_{\cU_b}$.

The Hilbert--Schmidt condition makes the variance of the Brownian
series summable in the state norm. For a finite partial sum, the It\^o isometry
gives
\[
 \E\left\|\sum_{k=1}^nG_be_k\beta_k(t)\right\|_\gamma^2
 =t\sum_{k=1}^n\|G_be_k\|_\gamma^2.
\]
The series therefore converges in $L^2(\Omega;\cE_\gamma)$.
The semigroup bound gives the same integrability for the stochastic
convolution constructed below.

With this realization,
$\zeta_{f,b}$ is given by the scalar It\^o series
\begin{equation}\label{eq:volterra-intro}
 \zeta_{t,f,b}
 =\sum_{k=1}^\infty\int_0^t\sqrt{f(u)}\,
   \langle\phi_{t-u},e_k\rangle_{\cU_b}\dd\beta_k(u),
\end{equation}
where the series converges in $L^2(\Omega)$.
If $E$ is a Hilbert space and
$K\in L^2(0,t;\cL_2(\cU_b,E))$ is deterministic, its stochastic integral is
\[
 \int_0^t K_u\dd W_u
 :=\sum_{k=1}^\infty\int_0^t K_ue_k\dd\beta_k(u),
 \qquad\text{with convergence in }L^2(\Omega;E).
\]
Indeed, the scalar It\^o isometry and independence of the Brownian
coordinates give, for $n>m$,
\[
 \E\left\|\sum_{k=m+1}^n\int_0^t K_ue_k\dd\beta_k(u)\right\|_E^2
 =\int_0^t\sum_{k=m+1}^n\|K_ue_k\|_E^2\dd u\longrightarrow0.
\]
This also gives the Hilbert-space It\^o isometry
$\E\|\int_0^tK_u\dd W_u\|_E^2=\int_0^t\|K_u\|_{\HS}^2\dd u$.
For the transport equation below,
$K_u=S_{t-u}\sqrt{f(u)}G_b$; the state-space estimates above give the
required integrability.

\begin{theorem}\label{thm:unified-lift}
Let $0<b<2$, let $0<\gamma<1$ satisfy $b+\gamma>1$, and let
$f\in L^1_{\mathrm{loc}}(\Rplus)$ be non-negative.  The stochastic evolution equation
\begin{equation}\label{eq:lift-spde}
 \dd\Theta_t=A\Theta_t\dd t+\sqrt{f(t)}G_b\dd W_t,
 \qquad \Theta_0=0,
\end{equation}
where the cylindrical noise is understood through the Brownian coordinates
$(\beta_k)$ fixed above, has a predictable mild solution, unique up to modification, given by
\begin{equation}\label{eq:lift-mild}
 \Theta_t
 =\sum_{k=1}^\infty\int_0^t
 S_{t-u}\sqrt{f(u)}G_be_k\dd\beta_k(u)
 \quad\text{in }L^2(\Omega;\cE_\gamma).
\end{equation}
For every $r\ge0$,
\begin{equation}\label{eq:lift-explicit}
 \Theta_t(r)
 =\sum_{k=1}^\infty\int_0^t\sqrt{f(u)}\,
 \langle\phi_{t+r-u},e_k\rangle_{\cU_b}\dd\beta_k(u),
\end{equation}
and the boundary process satisfies
\begin{equation}\label{eq:lift-boundary}
 \delta_0\Theta_t=\Theta_t(0)=\zeta_{t,f,b}.
\end{equation}
Moreover,
\begin{equation}\label{eq:lift-flow}
 \Theta_{t+h}(r)=\Theta_t(r+h)
 +\sum_{k=1}^\infty\int_t^{t+h}\sqrt{f(u)}\,
 \langle\phi_{t+h+r-u},e_k\rangle_{\cU_b}\dd\beta_k(u).
\end{equation}
The flow identity defines a time-inhomogeneous Gaussian Markov process with a
Feller transition evolution on $\cE_\gamma$.  If, in addition,
$f\in L^\infty_{\mathrm{loc}}(\Rplus)$, the mild solution admits a continuous
$\cE_\gamma$-valued modification.  Adding the deterministic clock makes the
process time-homogeneous Markov.
\end{theorem}

\begin{proof}
Fix $T>0$ and write $E:=\cE_\gamma$.  By
Lemmas~\ref{lem:state-shift} and~\ref{lem:G-HS}, for $t\le T$,
\[
 \int_0^tf(u)\|S_{t-u}G_b\|_{\HS(\cU_b,E)}^2\dd u
 \le C_{\gamma,T}^2\|G_b\|_{\HS(\cU_b,E)}^2
       \|f\|_{L^1(0,T)}<\infty.
\]
Thus $u\mapsto S_{t-u}\sqrt{f(u)}G_b$ belongs to
$L^2(0,t;\cL_2(\cU_b,E))$, so \eqref{eq:lift-mild} is a well-defined
$E$-valued stochastic integral.  The linear mild equation determines this integral uniquely at each time.
It remains to choose a predictable version and establish the stated path
properties.

The map $t\mapsto\Theta_t$ is continuous from $[0,T]$ into
$L^2(\Omega;E)$.  Indeed, extend each stochastic-integral kernel by zero to
$[0,T]$.  The It\^o isometry reduces the claim to convergence in
$L^2(0,T;\cL_2(\cU_b,E))$.  On the common part of the integration intervals,
strong continuity of $S$ implies
$S_{t_n-u}G_b\to S_{t-u}G_b$ in Hilbert--Schmidt norm, because
\[
 \|(S_{t_n-u}-S_{t-u})G_b\|_{\HS}^2
 =\sum_j\|(S_{t_n-u}-S_{t-u})G_be_j\|_E^2\to0
\]
by dominated convergence over an orthonormal basis $(e_j)$, using the local
semigroup bound.  The same bound together with $f\in L^1(0,T)$ supplies an
integrable majorant.  The contribution from the symmetric difference of the
integration intervals tends to zero by absolute continuity of the Lebesgue
integral.

We construct a predictable modification directly.  The $L^2(\Omega;E)$
continuity is uniform on $[0,T]$.  Choose successively finer deterministic
partitions and let $Y_n(t)$ equal $\Theta_{t_i}$ on each interval
$(t_i,t_{i+1}]$, with $Y_n(0)=0$, so that
\[
 \sup_{t\leq T}\E\|Y_n(t)-\Theta_t\|_E^2\leq2^{-3n}.
\]
Each $Y_n$ is predictable because $\Theta_{t_i}$ is
$\mathcal F_{t_i}$-measurable.  For every fixed $t$ and $\varepsilon>0$,
Chebyshev's inequality gives
\[
 \sum_{n=1}^\infty
 \mathbb P\bigl(\|Y_n(t)-\Theta_t\|_E>\varepsilon\bigr)
 \leq\varepsilon^{-2}\sum_{n=1}^\infty2^{-3n}<\infty.
\]
The Borel--Cantelli lemma, applied with $\varepsilon=1/m$ for
$m\in\mathbb N$, therefore gives $Y_n(t)\to\Theta_t$ almost surely.
Define the limit to be zero where the sequence does not converge.
Since $E$ is complete and separable, this defines a predictable process
and a modification of $\Theta$.  Applying this construction on successive
integer time intervals gives a predictable version on $\Rplus$.
Uniqueness at each fixed time follows from the mild formula, as asserted.

For every $r\ge0$, evaluation $\delta_r\theta:=\theta(r)$ is continuous on
$E$, because $\delta_r=\delta_0S_r$.  Moreover,
\[
 \delta_rS_{t-u}G_bh
 =(G_bh)(t-u+r)=\langle\phi_{t-u+r},h\rangle_{\cU_b}.
\]
Applying the bounded functional $\delta_r$ to \eqref{eq:lift-mild} proves
\eqref{eq:lift-explicit}.  In particular, Proposition~\ref{prop:factorization}
and the It\^o isometry give
\[
 \E[\Theta_s(0)\Theta_t(0)]
 =\int_0^{s\wedge t}f(u)
   \langle\phi_{s-u},\phi_{t-u}\rangle_{\cU_b}\dd u
 =R_{f,b}(s,t).
\]
Thus the boundary in \eqref{eq:lift-boundary} is a realization of the
Gaussian process defined by \eqref{eq:covariance-intro}.

Assume now that $f\in L^\infty(0,T)$.  Choose $p>2$ and
$a\in(1/p,1/2)$.  Then
\begin{equation}\label{eq:factorization-continuity-check}
 \sup_{t\le T}\int_0^t(t-u)^{-2a}f(u)
 \|S_{t-u}G_b\|_{\HS}^2\dd u
 \le\frac{C_{\gamma,T}^2\|f\|_{\infty,T}\|G_b\|_{\HS}^2
 T^{1-2a}}{1-2a}<\infty.
\end{equation}
Define
\[
 Z_a(t):=\sum_{k=1}^\infty\int_0^t
 (t-u)^{-a}S_{t-u}\sqrt{f(u)}G_be_k\dd\beta_k(u).
\]
The deterministic kernels belong to
$L^2(\{(t,u):0<u<t<T\};\cL_2(\cU_b,E))$ by
\eqref{eq:factorization-continuity-check}.  Approximation by simple
kernels therefore gives a jointly measurable representative of $Z_a$.
For a centered Gaussian $E$-valued random variable $Z$, set
$q=\E\|Z\|_E^2$.  Diagonalizing its covariance gives
$\|Z\|_E^2\stackrel{d}{=}\sum_jq_j\xi_j^2$, where
$q_j\geq0$, $\sum_jq_j=q$, and the $\xi_j$ are standard normal.
For $q>0$, convexity yields
\[
 \E\|Z\|_E^p
 \leq q^{p/2}\sum_j\frac{q_j}{q}\E|\xi_j|^p
 =\E|\xi_1|^p\,q^{p/2};
\]
the case $q=0$ is immediate.  Applying this bound to $Z_a(t)$ gives
\[
 \E\|Z_a(t)\|_E^p
 \leq C_p\bigl(\E\|Z_a(t)\|_E^2\bigr)^{p/2}.
\]
Together with \eqref{eq:factorization-continuity-check} and Tonelli's
theorem, this shows that $Z_a\in L^p(\Omega;L^p(0,T;E))$.
To interchange the integrals below, the local semigroup bound gives
\begin{align*}
 &\int_0^t(t-v)^{a-1}
 \left(\int_0^v(v-u)^{-2a}f(u)
       \|S_{t-u}G_b\|_{\HS}^2\dd u\right)^{1/2}\dd v\\
 &\qquad\leq C\int_0^t(t-v)^{a-1}v^{1/2-a}\dd v<\infty.
\end{align*}
For fixed $t$, the integrand, viewed as a function of $v$ with values in
$L^2(0,t;\cL_2(\cU_b,E))$, is thus Bochner integrable.  The It\^o
integral is a bounded linear map from that space into $L^2(\Omega;E)$,
so it commutes with this Bochner integral.  The semigroup law then gives
\begin{align*}
 &\int_0^t(t-v)^{a-1}S_{t-v}Z_a(v)\dd v\\
 &\quad=\int_0^t
 \left(\int_u^t(t-v)^{a-1}(v-u)^{-a}\dd v\right)
 S_{t-u}\sqrt{f(u)}G_b\dd W_u\\
 &\quad=\frac{\pi}{\sin(\pi a)}\Theta_t.
\end{align*}
The inner integral equals
$\int_0^1(1-z)^{a-1}z^{-a}\dd z
=\Gamma(a)\Gamma(1-a)=\pi/\sin(\pi a)$, by the substitution
$z=(v-u)/(t-u)$.

For completeness, define
\[
 (R_az)(t):=\int_0^t(t-v)^{a-1}S_{t-v}z(v)\dd v.
\]
Let $q=p/(p-1)$.  Since $a>1/p$, we have $(a-1)q>-1$, and
H\"older's inequality gives
\[
 \sup_{t\leq T}\|(R_az)(t)\|_E
 \leq C_{\gamma,T}
 \left(\frac{T^{(a-1)q+1}}{(a-1)q+1}\right)^{1/q}
 \|z\|_{L^p(0,T;E)}.
\]
If $z$ is continuous, then $R_az$ is continuous: after writing the
integral as $\int_0^t r^{a-1}S_rz(t-r)\dd r$, dominated convergence
applies on $[0,T]$, since $r^{a-1}$ is integrable.
Continuous $E$-valued functions are dense in $L^p(0,T;E)$; the displayed
bound therefore extends this conclusion to every $z\in L^p(0,T;E)$.
It follows that $\sin(\pi a)R_aZ_a/\pi$ is a continuous modification of
$\Theta$ on $[0,T]$.  On overlapping integer intervals, these
modifications agree almost surely at rational times and hence everywhere
by continuity.  They therefore give a continuous modification on
$\Rplus$.

Finally, split \eqref{eq:lift-mild} at time $t$ and use the semigroup law:
\[
 \Theta_{t+h}=S_h\Theta_t+
 \sum_{k=1}^\infty\int_t^{t+h}
 S_{t+h-u}\sqrt{f(u)}G_be_k\dd\beta_k(u).
\]
Applying $\delta_r$ gives \eqref{eq:lift-flow}.  For $0\le s\le t$, define
\[
 \Xi_{s,t}:=\sum_{k=1}^\infty\int_s^t
 S_{t-u}\sqrt{f(u)}G_be_k\dd\beta_k(u).
\]
Then $\Xi_{s,t}$ is independent of $\mathcal F_s$.  For bounded Borel
$F:E\to\R$, set
\[
 (P_{s,t}F)(\theta):=\E F(S_{t-s}\theta+\Xi_{s,t}).
\]
The preceding decomposition gives
$\E[F(\Theta_t)\mid\mathcal F_s]=(P_{s,t}F)(\Theta_s)$ almost surely.
Furthermore, for $s\le r\le t$,
\[
 \Xi_{s,t}=S_{t-r}\Xi_{s,r}+\Xi_{r,t},
\]
with $\Xi_{r,t}$ independent of $\mathcal F_r$; hence
$P_{s,t}=P_{s,r}P_{r,t}$.  These identities establish the Markov
transition evolution.  If $F\in C_b(E)$ and
$\theta_n\to\theta$ in $E$, then
$F(S_{t-s}\theta_n+\Xi_{s,t})\to F(S_{t-s}\theta+\Xi_{s,t})$ almost surely,
and bounded convergence shows $P_{s,t}F\in C_b(E)$.  Thus the evolution is
Feller. For bounded Borel $F:\Rplus\times E\to\R$, the clock-augmented
transition operators are
\[
 (\mathcal P_hF)(s,\theta)
 :=\E F(s+h,S_h\theta+\Xi_{s,s+h}).
\]
The deterministic kernels are measurable in their time parameters, and
the convolution identity above gives
$\mathcal P_{h+k}=\mathcal P_h\mathcal P_k$.  These are the
time-homogeneous Markov transitions on $\Rplus\times E$.
\end{proof}

\begin{remark}
\label{rem:direct-HS-convention}
The forward curve has the conditional-mean interpretation
\begin{equation}\label{eq:forward-conditional-mean}
 \Theta_t(r)=\E[\zeta_{t+r,f,b}\mid\mathcal F_t],
 \qquad t,r\geq0.
\end{equation}
Indeed, split the stochastic integral for $\zeta_{t+r,f,b}$ at $t$.
Its first part is \eqref{eq:lift-explicit}; its second part is centered
and independent of $\mathcal F_t$. Both parts are square integrable,
so taking conditional expectations proves the identity.
In particular, for $h\geq0$,
\begin{equation}\label{eq:boundary-prediction-increment}
 \zeta_{t+h,f,b}-\zeta_{t,f,b}
 =\Theta_t(h)-\Theta_t(0)
   +\sum_{k=1}^{\infty}\int_t^{t+h}\sqrt{f(u)}
       \langle\phi_{t+h-u},e_k\rangle_{\cU_b}\dd\beta_k(u).
\end{equation}
Thus the present curve determines the contribution of past noise to
every future increment. The remaining Gaussian process uses only fresh
Brownian increments; its covariance at future lags $h,k\ge0$ is
\[
 C_t(h,k)=\int_t^{t+\min(h,k)}f(u)
 Q_b(t+h-u,t+k-u)\dd u.
\]
The full curve therefore determines the conditional law of future
increments. The boundary value alone need not do so.

Although $\delta_0G_b=0$, noise reaches the boundary through the shift:
$\delta_0S_qG_bh=\langle\phi_q,h\rangle_{\cU_b}$ for $q>0$.
In particular,
\[
 \E|\Theta_t(0)|^2
 =\kappa_b\int_0^t f(u)(t-u)^b\dd u,
\]
which is positive exactly when $f>0$ on a set of positive Lebesgue measure
in $(0,t)$. Applying
$\delta_0$ to a formal strong differential equation would additionally
require a meaning for $\delta_0A\Theta_t$, which the mild construction
does not assume. Proposition~\ref{prop:boundary-slope} gives the
precise mean-square boundary-derivative criterion.

All these conditional statements use the given Brownian filtration;
no identification with the natural filtration of the scalar boundary
is assumed. At $b=1$, the same construction uses
$\cU_1=L^2((0,\infty),\dd\lambda)$ and $0<\gamma<1$.
\end{remark}

\section{Young and finite-variation equations}

We first establish the path regularity needed to integrate against the
boundary process.  We then solve the equation from every lifted state
and prove that the pair $(X,\Theta)$ is Markov and $C_b$-Feller.
The weight $f$ is deterministic, non-negative and Borel measurable.
The pointwise boundary-slope criterion below only requires local
integrability; the path-regularity estimates and the nonlinear equation
use the boundedness assumptions stated in their respective results.

\subsection{The regularity dichotomy}

\begin{proposition}\label{prop:young-regularity}
Let $T>0$, $0<b<1$, and $f\in L^\infty(0,T)$.  For
$0\leq s<t\leq T$,
\begin{equation}\label{eq:increment-bound}
 \E\abs{\zeta_{t,f,b}-\zeta_{s,f,b}}^2
 \leq\frac{\norm{f}_{\infty,T}}{1-b^2}\abs{t-s}^{1+b}.
\end{equation}
Consequently, $\zeta_{f,b}$ has a modification whose paths are H\"older
continuous of every order
\begin{equation}\label{eq:young-holder-range}
 0<\eta<H_b^\circ.
\end{equation}
In particular, the paths can be used as Young integrators: for
$1/2<\eta<H_b^\circ$, integration is defined against every integrand
of finite $q$-variation with $\eta+1/q>1$.
\end{proposition}

\begin{proof}
Let $h=t-s$.  From \eqref{eq:volterra-intro}, the It\^o isometry gives
\[
 \E\abs{\zeta_{t,f,b}-\zeta_{s,f,b}}^2
 =\int_0^s f(u)\norm{\phi_{t-u}-\phi_{s-u}}_{\cU_b}^2\dd u
 +\int_s^t f(u)\norm{\phi_{t-u}}_{\cU_b}^2\dd u.
\]
Since
$\phi_{r+h}(\lambda)-\phi_r(\lambda)
 =\sqrt b\,e^{-\lambda r}(1-e^{-\lambda h})/\lambda$, integrating first in
$r=s-u$ over $(0,\infty)$ gives
\[
 \int_0^s\norm{\phi_{r+h}-\phi_r}_{\cU_b}^2\dd r
 \leq\frac{b}{2\Gamma(2-b)}\int_0^\infty
 (1-e^{-\lambda h})^2\lambda^{-2-b}\dd\lambda
 =\frac{2^b-1}{1-b^2}h^{1+b}.
\]
For the last equality, rescale $x=\lambda h$ and use
$\int_0^\infty(1-e^{-x})^2x^{-2-b}\dd x
=2\Gamma(1-b)(2^b-1)/(b(1+b))$, which follows by two integrations by parts.
By \eqref{eq:phi-norm}, the fresh-noise term is at most
$\norm{f}_{\infty,T}\kappa_bh^{1+b}/(1+b)
=\norm{f}_{\infty,T}(2-2^b)h^{1+b}/(1-b^2)$.  Adding the two bounds proves
\eqref{eq:increment-bound}.  Gaussian moment equivalence gives
$\E|\zeta_t-\zeta_s|^q\leq C_{q,b,f,T}|t-s|^{qH_b^\circ}$ for every $q\geq2$,
and Kolmogorov's theorem \cite[Section~3.2]{DaPratoZabczyk2014} gives
every exponent $\eta<H_b^\circ-1/q$.  Letting $q\to\infty$ along a sequence
and identifying the continuous modifications on rational times yields one
modification with every exponent in \eqref{eq:young-holder-range}.
\end{proof}

\begin{proposition}
\label{prop:critical-regularity}
Let $T>0$, $b=1$, and $f\in L^\infty(0,T)$.  If
$0\leq s<t\leq T$ and $h:=t-s$, then
\begin{equation}\label{eq:critical-increment-bound}
 \E\abs{\zeta_{t,f,1}-\zeta_{s,f,1}}^2
 \leq \norm{f}_{\infty,T}h^2
 \left(3\log2+\frac12\log\frac{T}{h}\right).
\end{equation}
Consequently, $\zeta_{f,1}$ has a modification whose paths are H\"older
continuous of every order below one.  Such paths have finite $p$-variation
for every $p>1$ and zero quadratic variation along every sequence of
partitions whose mesh tends to zero.
\end{proposition}

\begin{proof}
At $b=1$, $\cU_1=L^2((0,\infty),\dd\lambda)$ and
$\phi_v'(\lambda)=e^{-\lambda v}$.  For $r,h\geq0$, Tonelli's theorem gives
\begin{equation}\label{eq:critical-D1}
 D_1(r,h):=\norm{\phi_{r+h}-\phi_r}_{\cU_1}^2
 =\int_r^{r+h}\int_r^{r+h}\frac{\dd v\,\dd w}{v+w}.
\end{equation}
Hence $D_1(r,h)\leq D_1(0,h)=2h\log2$ and $D_1(r,h)\leq h^2/(2r)$ for
$r>0$.  By the It\^o isometry and \eqref{eq:phi-norm},
\[
 \E\abs{\zeta_{t,f,1}-\zeta_{s,f,1}}^2
 =\int_0^s f(u)D_1(s-u,h)\dd u
 +2\log2\int_s^t f(u)(t-u)\dd u.
\]
Write $\log_+x:=\max\{\log x,0\}$ for $x>0$, with $\log_+0:=0$.
Splitting the first integral at $s-u=h$ bounds it by
$\norm{f}_{\infty,T}\bigl(2h^2\log2+\tfrac{h^2}{2}\log_+(s/h)\bigr)$, and the
second is at most $\norm{f}_{\infty,T}h^2\log2$.  Since
$\log_+(s/h)\leq\log(T/h)$, \eqref{eq:critical-increment-bound} follows.

For every $H<1$, $h^2\log(eT/h)\leq C_{H,T}h^{2H}$, so Gaussian moment
equivalence and Kolmogorov's theorem give a $\rho$-H\"older modification for
every $\rho<1$; a countable sequence $\rho_m\uparrow1$ and identification on
rational times give a single modification.  If $p>1$, choosing $\rho>1/p$
gives finite $p$-variation.  Finally, for any partition $\pi$ and fixed
$\rho>1/2$,
$\sum_{[u,v]\in\pi}\abs{\zeta_v-\zeta_u}^2
 \leq[\zeta]_{\rho;[0,T]}^2T|\pi|^{2\rho-1}\to0$.
\end{proof}

\begin{proposition}
\label{prop:boundary-slope}
Let $0<b<2$, let $0<\gamma<1$ satisfy $b+\gamma>1$, and let
$f\in L^1_{\mathrm{loc}}(\Rplus)$ be non-negative.  Let $\Theta$ be the
transport lift of Theorem~\ref{thm:unified-lift}.  For a fixed $t>0$, the
right derivative $\partial_r\Theta_t(0)$, understood as the limit of
$h^{-1}(\Theta_t(h)-\Theta_t(0))$ in $L^2(\Omega)$ as $h\downarrow0$,
exists if and only if
\begin{equation}\label{eq:boundary-slope-criterion}
 I_{f,b}(t):=\int_0^t f(u)(t-u)^{b-2}\dd u<\infty.
\end{equation}
When it exists,
\[
 \partial_r\Theta_t(0)
 =\sum_{k=1}^\infty\int_0^t\sqrt{f(u)}\,
   \langle\phi'_{t-u},e_k\rangle_{\cU_b}\dd\beta_k(u),
 \qquad
 \E|\partial_r\Theta_t(0)|^2=b2^{b-2}I_{f,b}(t).
\]
If $0<b\leq1$ and $f\geq c>0$ almost everywhere on
$(t-\varepsilon,t)$ for some $0<\varepsilon\le t$, then
$I_{f,b}(t)=\infty$.
\end{proposition}

\begin{proof}
For $h,q>0$,
$h^{-1}(\phi_{q+h}-\phi_q)(\lambda)
 =\sqrt b\,e^{-\lambda q}(1-e^{-\lambda h})/(\lambda h)$
increases pointwise to $\phi'_q(\lambda)$ as $h\downarrow0$.  The It\^o
isometry, monotone convergence and \eqref{eq:phi-derivative-norm} give
\[
 \lim_{h\downarrow0}
 \E\left|\frac{\Theta_t(h)-\Theta_t(0)}{h}\right|^2
 =\int_0^tf(u)\|\phi'_{t-u}\|_{\cU_b}^2\dd u
 =b2^{b-2}I_{f,b}(t),
\]
with value $+\infty$ allowed.  Hence an $L^2$ limit requires
\eqref{eq:boundary-slope-criterion}; if that integral is finite, dominated
convergence gives the asserted limit.  Under the stated lower bound on $f$,
$I_{f,b}(t)\geq c\int_0^\varepsilon q^{b-2}\dd q=\infty$ for $0<b\leq1$.
\end{proof}

Proposition~\ref{prop:boundary-slope} concerns differentiation in
$L^2(\Omega)$ at a fixed time. The next proposition additionally assumes
$f\in L^\infty_{\mathrm{loc}}(\Rplus)$ and shows that, for $1<b<2$,
the derivative exists at every positive time and has a locally
square-integrable version whose time integral recovers the boundary.

\begin{proposition}\label{prop:velocity}
Let $1<b<2$, choose $0<\gamma<1$, let
$f\in L^\infty_{\mathrm{loc}}(\Rplus)$ be non-negative, and let $\Theta$ be
the transport lift of Theorem~\ref{thm:unified-lift}.  For every $t>0$,
the map $r\mapsto\Theta_t(r)$ is differentiable at $r=0$ in $L^2(\Omega)$
and
\begin{equation}\label{eq:velocity-def}
 V_t:=\partial_r\Theta_t(0)
 =\sum_{k=1}^\infty\int_0^t\sqrt{f(u)}\,
   \langle\phi'_{t-u},e_k\rangle_{\cU_b}\dd\beta_k(u).
\end{equation}
Set $V_0:=0$.  Its covariance is
\begin{equation}\label{eq:velocity-cov}
 \E[V_sV_t]
 =b\int_0^{s\wedge t}f(u)(s+t-2u)^{b-2}\dd u.
\end{equation}
Moreover,
\begin{equation}\label{eq:zeta-integrated-velocity}
 \zeta_{t,f,b}=\int_0^tV_s\dd s
 \qquad\text{for all }t\geq0\text{ on one event of probability one}.
\end{equation}
Hence $\zeta_{f,b}$ has absolutely continuous paths and finite variation on
compact intervals.
\end{proposition}

\begin{proof}
For $q,h>0$, Jensen's inequality and \eqref{eq:phi-derivative-norm} give
\[
 \norm{\frac{\phi_{q+h}-\phi_q}{h}}_{\cU_b}^2
 =\norm{\int_0^1\phi'_{q+vh}\dd v}_{\cU_b}^2
 \leq\int_0^1\norm{\phi'_{q+vh}}_{\cU_b}^2\dd v
 \leq b(2q)^{b-2},
\]
because $b-2<0$.  The right-hand side is integrable at $q=0$ exactly when
$b>1$, so dominated convergence in the It\^o isometry proves
\eqref{eq:velocity-def}.  The identity
$\ip{\phi'_r}{\phi'_q}_{\cU_b}=b(r+q)^{b-2}$ gives \eqref{eq:velocity-cov},
and
\[
 \E\int_0^T\abs{V_s}^2\dd s
 =b2^{b-2}\int_0^T\int_0^s f(u)(s-u)^{b-2}\dd u\dd s
 \leq \frac{2^{b-2}\norm{f}_{\infty,T}T^b}{b-1}<\infty.
\]
The map
$s\mapsto\mathbf1_{\{u<s\}}\sqrt{f(u)}\phi'_{s-u}$ is Borel with values
in the separable space $L^2(0,T;\cU_b)$.  Countably valued Borel
approximations of this map, followed by the stochastic integral, give
a jointly measurable version of $V$.  The preceding bound then implies
$V\in L^2(0,T)$ almost surely.
The condition of the stochastic Fubini theorem
\cite[Section~4.5]{DaPratoZabczyk2014} holds, since
\[
 \int_0^T\left(\int_0^s f(u)
       \|\phi'_{s-u}\|_{\cU_b}^2\dd u\right)^{1/2}\dd s
 \leq T^{1/2}\left(\E\int_0^T|V_s|^2\dd s\right)^{1/2}<\infty.
\]
Using $\phi_q=\int_0^q\phi'_r\dd r$, it gives, for every $t\leq T$,
\[
 \int_0^tV_s\dd s
 =\sum_{k=1}^\infty\int_0^t\sqrt{f(u)}
   \left(\int_u^t
   \langle\phi'_{s-u},e_k\rangle_{\cU_b}\dd s\right)\dd\beta_k(u)
 =\zeta_{t,f,b}.
\]
For each fixed $t$, this is an equality almost surely.  The primitive of
$V\in L^2(0,T)$ is continuous, as is the boundary of the chosen lift.
Equality at rational times therefore gives equality for all $t\in[0,T]$
on one event of probability one.  Taking a countable intersection over
integer $T$ proves \eqref{eq:zeta-integrated-velocity} locally in time.
\end{proof}

\subsection{Coefficient regularity and pathwise well-posedness}
\label{subsec:weaker-coefficients}

For $h:\R\to\R$, write $\|h\|_\infty:=\sup_{x\in\R}|h(x)|$.
For $0<\alpha\leq1$, its spatial H\"older seminorm is
\[
 [h]_\alpha:=\sup_{x\ne y}\frac{|h(x)-h(y)|}{|x-y|^\alpha}.
\]
At $\alpha=1$ this is the Lipschitz seminorm.

\begin{assumption}
\label{ass:coefficients-mf}
The functions $\drift,\sigma:\R\times\Rplus\to\R$ are continuous. For
every $T>0$, there are constants $0\leq L_T,K_T<\infty$ such that
\[
 \begin{aligned}
 |\drift(x,t)-\drift(y,t)|+|\sigma(x,t)-\sigma(y,t)|
     &\leq L_T|x-y|,\\
 |\drift(0,t)|+|\sigma(0,t)|&\leq K_T,
 \end{aligned}
 \qquad x,y\in\R,\quad t\in[0,T].
\]
Continuity gives the bound at the origin.  Together the two bounds imply
\[
 |\drift(x,t)|+|\sigma(x,t)|\leq K_T+L_T|x|.
\]
When $0<b\leq1$, assume in addition that there is an exponent
$\alpha\in(0,1]$ such that
\begin{equation}\label{eq:coefficient-alpha-threshold}
 \alpha>(H_b^\circ)^{-1}-1,
\end{equation}
and $\sigma$ is continuously differentiable in $x$, with
\begin{align}
 \sup_{0\leq t\leq T}
 \left(\|\partial_x\sigma(\cdot,t)\|_\infty
       +[\partial_x\sigma(\cdot,t)]_\alpha\right)
 &\leq K_T,\label{eq:coefficient-spatial-c1alpha}\\
 \sup_{x\in\R}|\sigma(x,t)-\sigma(x,s)|
 &\leq K_T|t-s|,\qquad s,t\in[0,T].
 \label{eq:coefficient-time-lipschitz}
\end{align}
The constants may be enlarged to satisfy all these inequalities.
For $1<b<2$, only the first two displayed bounds and continuity are
assumed.
\end{assumption}

The time regularity of $\sigma$ also controls that of $\partial_x\sigma$.
We need this implication to apply the nonautonomous Young theorem below.

\begin{lemma}
\label{lem:coefficient-time-interpolation}
Under \eqref{eq:coefficient-spatial-c1alpha}--\eqref{eq:coefficient-time-lipschitz},
put $\beta_*=\alpha/(1+\alpha)$. For every $T>0$,
\begin{equation}\label{eq:coefficient-time-derivative}
 \sup_{x\in\R}\bigl(
 |\sigma(x,t)-\sigma(x,s)|
 +|\partial_x\sigma(x,t)-\partial_x\sigma(x,s)|\bigr)
 \leq C_T|t-s|^{\beta_*},\qquad s,t\in[0,T].
\end{equation}
\end{lemma}

\begin{proof}
Fix $s\neq t$ in $[0,T]$, and put $d:=|t-s|$ and
$g:=\sigma(\cdot,t)-\sigma(\cdot,s)$.  Then $\|g\|_\infty\leq K_Td$ and
$[g']_\alpha\leq2K_T$.  For $x\in\R$ and $h>0$,
\[
 g'(x)=\frac{g(x+h)-g(x)}{h}
       -\frac1h\int_0^h\bigl(g'(x+r)-g'(x)\bigr)\dd r.
\]
Taking absolute values gives
\[
 \|g'\|_\infty\leq\frac{2K_Td}{h}+\frac{2K_T}{1+\alpha}h^\alpha.
\]
The choice $h=d^{1/(1+\alpha)}$ gives $\|g'\|_\infty\leq4K_Td^{\beta_*}$.
Since $d\leq T$, also $\|g\|_\infty\leq K_TT^{1-\beta_*}d^{\beta_*}$.  Thus
\eqref{eq:coefficient-time-derivative} holds with
$C_T=K_T(4+T^{1-\beta_*})$.
\end{proof}

\paragraph{Path norms.}
For $p\geq1$, let $C^{p\text{-var}}([0,T])$ denote the continuous paths of
finite $p$-variation (see the notation paragraph), with the norm
\[
 N_p(z):=|z(0)|+\|z\|_{p\text{-var};[0,T]}\geq\|z\|_{\infty;[0,T]}.
\]
For a continuous path $v$ of bounded variation, put
$B_T(v):=|v(0)|+\TV_{[0,T]}(v)$; then $N_p(v)\leq B_T(v)$.  If $z$ is
$\eta$-H\"older and $p\eta\geq1$, then
$\|z\|_{p\text{-var};[0,T]}\leq[z]_{\eta;[0,T]}T^\eta$, because
$\sum_i(t_{i+1}-t_i)^{p\eta}\leq T^{p\eta}$.  Young integration is the limit
of Riemann--Stieltjes sums when the integrand has finite $q$-variation and
$1/p+1/q>1$; see \cite[Section~2.1, (2.5)--(2.7)]{CongDucHong2018}.

\begin{proposition}
\label{prop:weak-coefficient-young}
Assume \eqref{eq:coefficient-spatial-c1alpha}--\eqref{eq:coefficient-time-lipschitz}
and the continuity and global spatial Lipschitz hypotheses on the
drift in Assumption~\ref{ass:coefficients-mf}.
Let $1<p<\min(2,1+\alpha)$, let $z\in C^{p\text{-var}}([0,T])$,
and let $v$ be a continuous path of bounded variation. Then
\begin{equation}\label{eq:weak-coefficient-young}
 x_t=x_0+\int_0^t\drift(x_r,r)\dd r
      +\int_0^t\sigma(x_r,r)\dd(z_r+v_r)
\end{equation}
has a unique solution in $C^{p\text{-var}}([0,T])$.
The last integral equals the sum of the Young integral against $z$ and
the Lebesgue--Stieltjes integral against $v$.
The solution depends continuously, in the uniform path topology, on
$(x_0,z,v)$ when $z$ carries the norm $N_p$ and $v$ the norm $B_T$.
\end{proposition}

\begin{proof}
Set $w=z+v$.  Since $N_p(v)\leq B_T(v)$, this driver has finite
$p$-variation.  We verify the hypotheses of
\cite[Theorems~3.6 and~3.9]{CongDucHong2018}; their diffusion coefficient
$g(t,x)$ is our $\sigma(x,t)$ and their drift $f(t,x)$ is our
$\drift(x,t)$.

In their notation, take $\delta=\alpha$ and
$\beta=\beta_*=\alpha/(1+\alpha)$.  The bound on $\partial_x\sigma$
implies global spatial Lipschitz continuity, and its $\alpha$-H\"older
bound gives $\mathbf H_1$(i)--(ii).  Lemma~\ref{lem:coefficient-time-interpolation}
gives $\mathbf H_1$(iii) with the additive control
$h(s,t)=C_T^{1/\beta_*}(t-s)$.  Choose $r\in(1/p,1)$ for the parameter
called $\alpha$ in their drift assumption; this choice is independent of
our spatial exponent $\alpha$.  Then $r>1/2$ because $p<2$, and
\[
 |\drift(x,t)|\leq L_T|x|+A_T,
 \qquad A_T:=\sup_{t\leq T}|\drift(0,t)|<\infty.
\]
The constant envelope $A_T$ belongs to $L^{1/(1-r)}(0,T)$.  This growth
bound and spatial Lipschitz continuity verify $\mathbf H_2$.
The three restrictions in $\mathbf H_3$ hold because
\[
 \alpha>p-1,\qquad
 \beta_* >1-\frac1p,\qquad
 \alpha r>\frac{\alpha}{p}>1-\frac1p.
\]
To check the variation exponents used in those theorems, choose $q_0$
and then $q$ so that
\[
 1-\frac1p<\frac1{q_0}
 <\min\left\{\beta_*,\alpha r,\frac\alpha p,\frac12\right\},
 \qquad
 \frac1{q_0\alpha}\leq\frac1q<\frac1p.
\]
Both choices are possible: every entry in the first minimum exceeds
$1-1/p$, and $1/(q_0\alpha)<1/p<r$.  These are precisely the choices
in \cite[(3.2)--(3.3)]{CongDucHong2018}, after replacing their drift
parameter by $r$.  In particular, $q>p$.  Their Theorem~3.6 gives
uniqueness in $C^{q\text{-var}}$ and a solution of finite
$\max\{p,1/r\}=p$ variation.  Since
$C^{p\text{-var}}\subseteq C^{q\text{-var}}$, it gives the claimed
existence and uniqueness in $C^{p\text{-var}}$.

Their Theorem~3.9 gives joint continuity of the evaluation map in time,
initial value and driver.  Suppose
$x_0^n\to x_0$, $N_p(z_n-z)\to0$, and $B_T(v_n-v)\to0$.
Then $N_p((z_n+v_n)-(z+v))\to0$.  This convergent data sequence and its
limit form a compact set.  The evaluation map is uniformly continuous
on its product with $[0,T]$, so the corresponding solutions converge
uniformly on $[0,T]$.

Finally, the coefficient bounds give
\[
 |\sigma(x_t,t)-\sigma(x_s,s)|
 \leq K_T|x_t-x_s|+K_T|t-s|.
\]
Thus the integrand has finite $p$-variation.  Since $p<2$, its Young
integrals against $w$ and $z$ exist.  Splitting each Riemann sum into
its $z$ and $v$ parts proves the asserted decomposition; the sums
against $v$ converge to the Lebesgue--Stieltjes integral because the
integrand is continuous and $v$ has bounded variation.
\end{proof}

\begin{remark}
For $0<b<1$, \eqref{eq:coefficient-alpha-threshold} is
$\alpha>(1-b)/(1+b)$; for $b=1$, it asks for any $\alpha>0$.  These are
sufficient conditions; we make no claim of necessity.  The
nonautonomous theorem allows time-dependent coefficients and a spatially
Lipschitz drift.  The
well-posedness and Markov--Feller result below uses
Assumption~\ref{ass:coefficients-mf}.  The quantitative estimates in
Section~\ref{sec:spectral-ou} impose the stronger, separately stated
Assumptions~\ref{ass:young-quantitative} and~\ref{ass:coefficients}.
\end{remark}

\subsection{The nonlinear lifted equation}

Throughout this subsection, fix $0<b<2$ and $0<\gamma<1$ with
$b+\gamma>1$, let $f\in L^\infty_{\mathrm{loc}}(\Rplus)$ be non-negative,
and impose Assumption~\ref{ass:coefficients-mf}.

\paragraph{Restarted states.}
Fix $s\ge0$ and $(x,\theta)\in\R\times\cE_\gamma$.  Let $W$ be a
cylindrical Brownian motion on $\cU_b$, adapted to its usual augmented
filtration $(\mathcal F_r)_{r\ge0}$, indexed by elapsed time $r$.
Define
\begin{equation}\label{eq:shifted-forward-explicit}
 \Theta_r^{s,\theta}(\rho)
 :=\theta(r+\rho)
 +\sum_{k=1}^\infty\int_0^r\sqrt{f(s+u)}\,
   \langle\phi_{r+\rho-u},e_k\rangle_{\cU_b}\dd\beta_k(u),
 \qquad r,\rho\ge0.
\end{equation}
Equivalently, $\Theta_r^{s,\theta}=S_r\theta+\Xi_r^s$ in $\cE_\gamma$, where
$\Xi^s:=\Theta^{s,0}$ is the stochastic convolution of
Theorem~\ref{thm:unified-lift} with the locally bounded weight
$f(s+\cdot)$.  By that theorem, $\Xi^s$ has a continuous
$\cE_\gamma$-valued modification; we fix it and define $\Theta^{s,\theta}$
by the preceding identity.  Since $\delta_0$ is bounded on $\cE_\gamma$,
the boundary
\[
 \begin{aligned}
 \zeta_r^{s,\theta}&:=\Theta_r^{s,\theta}(0)=\theta(r)+\zeta_r^{s,0},\\
 \zeta_r^{s,0}&:=\Xi_r^s(0)
 =\sum_{k=1}^\infty\int_0^r\sqrt{f(s+u)}\,
   \langle\phi_{r-u},e_k\rangle_{\cU_b}\dd\beta_k(u),
 \end{aligned}
\]
is continuous in $r$.  Any continuous modification of the fresh path
$\zeta^{s,0}$ supplied by Propositions~\ref{prop:young-regularity},
\ref{prop:critical-regularity} or~\ref{prop:velocity}, applied to
$f(s+\cdot)$, is indistinguishable from this one.

By Lemma~\ref{lem:state-shift}, for every $\vartheta\in\cE_\gamma$ and $T>0$,
\begin{equation}\label{eq:boundary-variation-control}
 B_T(\vartheta|_{[0,T]})
 \le\bigl(1+J_\gamma(T)^{1/2}\bigr)\|\vartheta\|_\gamma,
 \qquad
 \sup_{r\le T}\|S_r\vartheta\|_\gamma\le C_{\gamma,T}\|\vartheta\|_\gamma.
\end{equation}
In particular, restriction maps $\cE_\gamma$ continuously into the space
of continuous bounded-variation paths equipped with $B_T$.  We use this
space for the initial curve because a general state need not have the
Gaussian path's H\"older regularity.  For example,
$\theta(\rho)=\sqrt\rho$ satisfies
\[
 \|\theta\|_\gamma^2=\frac14\int_0^\infty
       \frac{\rho^{\gamma-1}}{(1+\rho)^2}\dd\rho<\infty,
 \qquad 0<\gamma<1,
\]
but is not $\eta$-H\"older at zero for $\eta>1/2$.  For $f\equiv0$,
the restarted boundary is $\sqrt r$.  The bounded-variation term in
Proposition~\ref{prop:weak-coefficient-young} allows this state as well
as every other state in $\cE_\gamma$.

\begin{remark}
Fix $r>0$.  Since $\theta$ is deterministic and $\Xi_r^s$ is centered, the
quotient $h^{-1}(\Theta_r^{s,\theta}(h)-\Theta_r^{s,\theta}(0))$ is the sum
of $h^{-1}(\theta(r+h)-\theta(r))$ and a centered term.  Taking expectations
isolates the deterministic quotient.  Consequently,
Proposition~\ref{prop:boundary-slope} shows that convergence in $L^2(\Omega)$ as
$h\downarrow0$ holds if and only if $\theta'_+(r)$ exists and
$\int_0^rf(s+u)(r-u)^{b-2}\dd u<\infty$; the limit then has second moment
$|\theta'_+(r)|^2+b2^{b-2}\int_0^rf(s+u)(r-u)^{b-2}\dd u$.
The equation below only uses bounded variation of the initial curve;
these pointwise derivative conditions are not required.
\end{remark}

\begin{theorem}
\label{thm:nonlinear-lift}
Let $0<b<2$, let $0<\gamma<1$ satisfy $b+\gamma>1$, let
$f\in L^\infty_{\mathrm{loc}}(\Rplus)$ be non-negative, and impose
Assumption~\ref{ass:coefficients-mf}. For every $s\ge0$ and
$(x,\theta)\in\R\times\cE_\gamma$, the equation
\begin{equation*}
 X_r^{s,x,\theta}
 =x+\int_0^r\drift(X_u^{s,x,\theta},s+u)\dd u
 +\int_0^r\sigma(X_u^{s,x,\theta},s+u)\dd\zeta_u^{s,\theta}
\end{equation*}
has the following pathwise solutions.  For $1<b<2$, it has a unique
continuous solution, and the last integral is Lebesgue--Stieltjes.
For $0<b\le1$ and every
$(H_b^\circ)^{-1}<p<\min(2,1+\alpha)$, it has a unique solution of
finite $p$-variation.  In this case the last integral is the sum of a
Lebesgue--Stieltjes integral against $\theta$ and a Young integral
against $\zeta^{s,0}$.  At $b=1$, the admissible range is
$p\in(1,1+\alpha)$.

The process
$\mathcal X_r^{s,x,\theta}:=\bigl(X_r^{s,x,\theta},\Theta_r^{s,\theta}\bigr)$
is time-inhomogeneous Markov on $\R\times\cE_\gamma$, and its transition
evolution is $C_b$-Feller:
\[
 P_{u,t}F\in C_b(\R\times\cE_\gamma)
 \quad\text{for every }F\in C_b(\R\times\cE_\gamma)
 \text{ and }0\le u\le t.
\]
For $s=0$ and $\theta=0$, the first component is the pathwise solution
of \eqref{eq:main-equation-intro} driven by the realization
\eqref{eq:lift-boundary}.
\end{theorem}

\begin{proof}
Fix $T>0$.  All coefficient bounds are taken on $\R\times[s,s+T]$, and
$L_T$ denotes a common spatial Lipschitz constant of $\drift$ and $\sigma$
there.

\paragraph{Finite-variation case.}
Suppose $1<b<2$.  Proposition~\ref{prop:velocity}, applied to
$f(s+\cdot)$, gives an absolutely continuous fresh path.  By
\eqref{eq:boundary-variation-control}, $\zeta^{s,\theta}$ is therefore
continuous and has finite variation.  Fix such a driving path and put
$\mu(\dd r):=\dd r+|\dd\zeta_r^{s,\theta}|$.  This is a finite atomless
measure on $[0,T]$.

Set $x^{(0)}_r=x$ and define the Picard iterates by
\[
 x^{(m+1)}_r=x+\int_0^r\drift(x^{(m)}_u,s+u)\dd u
                 +\int_0^r\sigma(x^{(m)}_u,s+u)\dd\zeta_u^{s,\theta}.
\]
They are continuous.  Let
$A_x:=\sup_{u\leq T}(|\drift(x,s+u)|+|\sigma(x,s+u)|)<\infty$.
The first difference is bounded by $A_x\mu([0,r])$.  Iterating the
Lipschitz estimate and using the fact that $\mu$ has no atoms gives
\[
 \sup_{u\leq r}|x^{(m+1)}_u-x^{(m)}_u|
 \leq A_xL_T^m\frac{\mu([0,r])^{m+1}}{(m+1)!},
 \qquad m\geq0.
\]
Indeed, each iteration integrates over an ordered simplex; its
$\mu^{\otimes(m+1)}$-mass is $\mu([0,r])^{m+1}/(m+1)!$ because its
permutations partition the product space up to the diagonals, which
have zero measure.  The series of these bounds converges, so the
iterates converge uniformly.  Spatial Lipschitz continuity lets us pass
to the limit in both integrals, giving a solution on $[0,T]$.

The same ordered-simplex argument proves the Gronwall inequality for
$\mu$: if a continuous $D\geq0$ satisfies
$D(r)\leq a+L\int_0^rD\dd\mu$, then
$D(r)\leq a\exp(L\mu([0,r]))$.  Apply this first to the difference of
two solutions, with $a=0$, to obtain uniqueness.  Apply it to
$1+|X_r^{s,x,\theta}|$, using linear growth of both coefficients, to
obtain, for a deterministic constant $C_{s,T}$,
\[
 1+\|X^{s,x,\theta}\|_{\infty;[0,T]}
 \leq(1+|x|)\exp\!\left(C_{s,T}
           [T+\TV_{[0,T]}(\zeta^{s,\theta})]\right).
\]
The construction therefore extends to every finite time interval.
It also gives measurability without choosing path-dependent time
partitions.  The set of continuous bounded-variation paths is Borel in
$C([0,T])$, since total variation is the supremum of the partition sums
over rational interior partition points and the two endpoints.
For continuous integrands and bounded-variation integrators, Stieltjes
integrals are limits of Riemann sums on deterministic refining
partitions.  Every Picard iterate, and hence its uniform limit, is
therefore Borel in $x$, $\theta$, and the fresh driving path.

\paragraph{Young case.}
For $0<b\le1$, fix $p$ as in the statement; this is possible because
\eqref{eq:coefficient-alpha-threshold} gives $1+\alpha>(H_b^\circ)^{-1}$
and $H_b^\circ>1/2$.  Choose $\eta\in(1/p,H_b^\circ)$.  By
Propositions~\ref{prop:young-regularity} and~\ref{prop:critical-regularity},
applied to $f(s+\cdot)$, the fresh path $z:=\zeta^{s,0}$ is
$\eta$-H\"older, hence in $C^{p\text{-var}}([0,T])$, and
$v:=\theta|_{[0,T]}$ has bounded variation.
Proposition~\ref{prop:weak-coefficient-young}, applied to
$(y,u)\mapsto(\drift(y,s+u),\sigma(y,s+u))$, gives existence and uniqueness
in $C^{p\text{-var}}([0,T])$, the stated splitting of the integral, and
continuity of the solution in $(x,z,v)$ for the norms $N_p$ on $z$ and
$B_T$ on $v$.
Solutions constructed with two admissible exponents agree by uniqueness
for the larger exponent.  A countable sequence of admissible exponents
decreasing to $(H_b^\circ)^{-1}$ therefore gives one version with finite
$p$-variation for every $p$ in the stated range.

For measurability, let $C^{0,p\text{-var}}([0,T])$ be the closure of smooth
paths in the norm $N_p$, a separable Banach space
\cite[Section~2.1, following (2.1)]{CongDucHong2018}.  For
$1/p<\eta_1<\eta$, let $z^{(n)}$ be the piecewise linear interpolant on
the dyadic grid of mesh $T2^{-n}$.  The interpolation estimate
\[
 \|z^{(n)}-z\|_{\Holder{\eta_1}([0,T])}
 \leq C_{\eta,\eta_1,T}2^{-n(\eta-\eta_1)}[z]_{\eta;[0,T]}
\]
follows from
$\|z^{(n)}-z\|_\infty\leq C(T2^{-n})^\eta[z]_\eta$ and
$[z^{(n)}-z]_\eta\leq C[z]_\eta$: use the second bound for increments
below the mesh size and the first for larger increments.
Hence $z^{(n)}\to z$ in $N_p$.  Each interpolant belongs
to $C^{0,p\text{-var}}$: smoothing its piecewise constant derivative in
$L^1$ gives convergence in total variation and thus in $N_p$.
It is also a continuous function, with values in that space, of
finitely many values of $z$.  The almost sure limit therefore makes
$z$ a Borel random element of $C^{0,p\text{-var}}([0,T])$.
Composing the
continuous restriction $\theta\mapsto\theta|_{[0,T]}$ into the norm $B_T$,
see \eqref{eq:boundary-variation-control}, with the continuous solution map
shows that the solution is a continuous, hence Borel, function of
$(x,\theta,z)\in\R\times\cE_\gamma\times C^{0,p\text{-var}}([0,T])$.

In the finite-variation case the Picard construction is nonanticipating.
In the Young case, the same solution construction on $[0,r]$ agrees with
the restriction of the solution on $[0,T]$ by uniqueness.  Together
with the measurability just proved, this shows that $X$ is adapted in
both cases.  Its continuous paths also give progressive measurability.

\paragraph{Restart and Markov property.}
For $u\le t$ and bounded Borel $F$ on
$\R\times\cE_\gamma$, define
\[
 P_{u,t}F(x,\theta)
 :=\E\!\left[F\!\left(X_{t-u}^{u,x,\theta},
                         \Theta_{t-u}^{u,\theta}\right)\right].
\]
By the measurability proved above and the continuity of
$\theta\mapsto S_r\theta$, these are
Borel transition kernels.  Fix $0\le r\le T$.  Splitting the stochastic
integral in \eqref{eq:shifted-forward-explicit} at elapsed time $r$, as in
the proof of \eqref{eq:lift-flow}, gives
\[
 \Theta_{r+h}^{s,\theta}
 =S_h\Theta_r^{s,\theta}+\Xi_h^{s+r,(r)},
 \qquad 0\le h\le T-r,
\]
simultaneously for all $h$ almost surely by continuity of the chosen
versions.  Here $\Xi^{s+r,(r)}$ is the stochastic convolution with weight
$f(s+r+\cdot)$ driven by $W_{r+\cdot}-W_r$.  Its entire path is
independent of $\mathcal F_r$.  The future boundary, in elapsed time
$h$, is therefore
\[
 \widehat\zeta_h:=\Theta_r^{s,\theta}(h)+\Xi_h^{s+r,(r)}(0)
 =\zeta_{r+h}^{s,\theta}.
\]
The first summand is the restriction of the new state to its horizon
variable, and has bounded variation by
\eqref{eq:boundary-variation-control}.  The second is a fresh Gaussian
path with the regularity used above.  Additivity and time translation
of the pathwise integrals show that
$\widehat X_h:=X_{r+h}^{s,x,\theta}$ satisfies
\[
 \widehat X_h=X_r^{s,x,\theta}
 +\int_0^h\drift(\widehat X_u,s+r+u)\dd u
 +\int_0^h\sigma(\widehat X_u,s+r+u)\dd\widehat\zeta_u.
\]
By uniqueness, this is the restarted solution from
$\mathcal X_r^{s,x,\theta}$ at time $s+r$.
The measurable solution maps constructed above apply to every state
on the event where the fresh path has the stated regularity.
The future is thus a measurable function of the
$\mathcal F_r$-measurable state and of fresh noise independent of
$\mathcal F_r$, which gives
\[
 \E\!\left[F(\mathcal X_q^{s,x,\theta})\mid\mathcal F_r\right]
 =P_{s+r,s+q}F\!\left(\mathcal X_r^{s,x,\theta}\right)
 \quad\text{almost surely},\qquad 0\le r\le q\le T.
\]
Applying this at an intermediate time also gives $P_{u,t}=P_{u,v}P_{v,t}$
for $u\le v\le t$.

\paragraph{Continuity in the lifted state.}
Let $(x_n,\theta_n)\to(x,\theta)$ in
$\R\times\cE_\gamma$, and couple all
initial states through the same fresh noise.  Then
$\Theta_r^{s,\theta_n}-\Theta_r^{s,\theta}=S_r(\theta_n-\theta)$, so by
\eqref{eq:boundary-variation-control} the curves converge in $\cE_\gamma$
uniformly on $[0,T]$, and $B_T(\theta_n-\theta)\to0$.

If $1<b<2$, put $K:=\sup_n\TV_{[0,T]}(\theta_n)<\infty$ and
$D_n(r):=|X_r^{s,x_n,\theta_n}-X_r^{s,x,\theta}|$.
Subtract the equations, integrating the difference of the diffusion
coefficients against $\zeta^{s,0}+\theta_n$ and the remaining term
against $\theta_n-\theta$.  With
$S:=\sup_{r\leq T}|\sigma(X_r^{s,x,\theta},s+r)|<\infty$,
this gives
\[
 D_n(r)\leq |x_n-x|+S\TV_{[0,T]}(\theta_n-\theta)
 +L_T\int_0^rD_n(u)
        \bigl(\dd u+|\dd\zeta_u^{s,0}|+|\dd\theta_n(u)|\bigr).
\]
The atomless-measure Gronwall inequality therefore yields
\[
 \begin{aligned}
 &\sup_{r\le T}\abs{X_r^{s,x_n,\theta_n}-X_r^{s,x,\theta}}\\
 &\quad\le \exp\!\left(L_T[T+\TV_{[0,T]}(\zeta^{s,0})+K]\right)\\
 &\qquad\times
 \left(\abs{x_n-x}+S\,\TV_{[0,T]}(\theta_n-\theta)\right)\longrightarrow0.
 \end{aligned}
\]
If $0<b\le1$, the continuity of the Young solution map, with $z=\zeta^{s,0}$ fixed and
$B_T(\theta_n-\theta)\to0$, gives the same uniform convergence almost
surely.  Hence, for $F\in C_b(\R\times\cE_\gamma)$, bounded convergence
shows that $P_{u,t}F$ is continuous.  Since
$\|P_{u,t}F\|_\infty\leq\|F\|_\infty$, it belongs to
$C_b(\R\times\cE_\gamma)$.

Finally, for $s=0$ and $\theta=0$, \eqref{eq:shifted-forward-explicit}
reduces to \eqref{eq:lift-explicit}, so $\zeta^{0,0}$ is the realization
\eqref{eq:lift-boundary}.
\end{proof}

\begin{remark}
For $\kappa>0$ and $\sigma_0\in\R$, the mean-reverting coefficients
$\drift(x,t)=-\kappa x$ and $\sigma(x,t)=\sigma_0$ satisfy
Assumption~\ref{ass:coefficients-mf}. The corresponding restarted solution is
\[
 X_r^{s,x,\theta}
 =e^{-\kappa r}x+\sigma_0\int_0^r e^{-\kappa(r-u)}\dd\zeta_u^{s,\theta},
\]
where the integral has the Young or Lebesgue--Stieltjes interpretation
specified in Theorem~\ref{thm:nonlinear-lift}.
For linear multiplicative diffusion,
$\drift(x,t)=\nu x$ and $\sigma(x,t)=\beta x$, with
$\nu,\beta\in\R$, the solution is
\[
 X_r^{s,x,\theta}
 =x\exp\!\left(\nu r+\beta[\zeta_r^{s,\theta}-\zeta_0^{s,\theta}]\right).
\]
Both formulas follow from the pathwise integration-by-parts and chain
rules.  The driving paths are continuous and have finite $p$-variation
for some $p<2$ in both regimes, so the chain rule has no quadratic-variation
term.
\end{remark}

\section{Spectral Ornstein--Uhlenbeck representation and multifactor approximation}\label{sec:spectral-ou}

The Laplace factorization expresses the boundary process as a superposition
of integrated Ornstein--Uhlenbeck factors. Retaining finitely many
relaxation rates gives a finite-dimensional Markovian approximation.
We couple it with the original process, estimate the spectral error in
the path norms needed by the nonlinear equation, and then control the
additional error from time discretization.

Throughout this section, $f$ is deterministic, nonnegative, and locally
integrable. Estimates are stated on a fixed interval $[0,T]$; the
boundedness of $f$ required for the convergence rates is imposed in
the corresponding results.

\paragraph{Spaces and notation.}
We use $\cU_b=L^2((0,\infty),\mu_b)$ for functions of the spectral
variable $\lambda$, and $\cE_\gamma$ for curves of the forward horizon
$r$, as defined in \eqref{eq:mu-intro} and \eqref{eq:E-intro}.
The former indexes the cylindrical noise; the latter is the state space
of the lift, whose boundary is the scalar path $t\mapsto\zeta_{t,f,b}$.
For comparisons with the lift, fix $0<\gamma<1$ with $b+\gamma>1$.
The finite-factor state space is $\R\times\R^M$ with its Euclidean
topology.

For a separable Banach space $E$ and $1\leq p<\infty$,
$L^p(\Omega;E)$ denotes the space of strongly measurable $E$-valued
random variables, identified almost surely, with finite norm
\[
 \|U\|_{L^p(\Omega;E)}:=\bigl(\E\|U\|_E^p\bigr)^{1/p}.
\]
The Bochner space $L^q(0,T;E)$, $1\leq q<\infty$, consists of strongly
measurable functions, identified almost everywhere, with finite norm
$(\int_0^T\|u(t)\|_E^q\dd t)^{1/q}$; we omit $E$ when $E=\R$.
Thus the mixed norm is
\[
 \|D\|_{L^p(\Omega;L^q(0,T))}
 :=\left(\E\left[\left(\int_0^T|D_t|^q\dd t\right)^{p/q}\right]
   \right)^{1/p}.
\]
For $1<b<2$, the $L^2(0,T)$ norm of a velocity error controls the total
variation of its primitive. For $0<b\leq1$, we estimate driver errors
in the H\"older spaces $\Holder{\eta}([0,T])$ of the notation paragraph;
their separable subspaces $c^\eta([0,T])$, defined below, provide
measurable realizations with the same norm. Under the weaker coefficient
assumptions, continuity of the Young solution map is used in
$C^{p\text{-var}}([0,T])$ with norm $N_p$, as defined in
Section~\ref{subsec:weaker-coefficients}.

\subsection{From Laplace features to Ornstein--Uhlenbeck factors}\label{subsec:ou-factors}

For $\lambda>0$, define
\[
 k_\lambda(r):=\frac{1-e^{-\lambda r}}{\lambda}
              =\int_0^r e^{-\lambda v}\dd v.
\]
For every $0<b<2$, including $b=1$, the Laplace feature satisfies
$\phi_r(\lambda)=\sqrt b\,k_\lambda(r)$. Let $\mathcal W$ be the
orthogonal Gaussian random measure associated with $W$, with control
measure $\mu_b(\dd\lambda)\dd u$. Thus
\[
 W_t(h)=\int_0^t\int_0^\infty h(\lambda)
 \mathcal W(\dd\lambda,\dd u),
 \qquad h\in\cU_b.
\]
The boundary process has the representation
\begin{equation}\label{eq:integrated-ou-rep}
 \zeta_{t,f,b}
 =\sqrt b\int_0^t\int_0^\infty
 \frac{1-e^{-\lambda(t-u)}}{\lambda}\sqrt{f(u)}
 \mathcal W(\dd\lambda,\dd u).
\end{equation}

The integrated kernel is square integrable over the full spectrum.
Indeed, $0\leq k_\lambda(r)\leq\min\{r,\lambda^{-1}\}$ and
$\mu_b(\dd\lambda)=\lambda^{1-b}\dd\lambda/\Gamma(2-b)$ imply,
by splitting the integral at $\lambda=r^{-1}$, that
\[
 \int_0^\infty k_\lambda(r)^2\mu_b(\dd\lambda)
 \leq\frac{1}{\Gamma(2-b)}
 \left(r^2\int_0^{1/r}\lambda^{1-b}\dd\lambda
       +\int_{1/r}^\infty\lambda^{-1-b}\dd\lambda\right)
 \leq C_b r^b,\qquad r>0.
\]
Consequently the variance in \eqref{eq:integrated-ou-rep} is bounded by
$C_b\int_0^t f(u)(t-u)^b\dd u$, which is finite on every bounded time
interval. This calculation explains why integration in time permits
the same representation throughout $0<b<2$, including $b=1$.

For $h\in\cU_b$, define
\begin{equation*}
 \mathsf M_t(h)
 :=\int_0^t\int_0^\infty h(\lambda)\sqrt{f(u)}
 \mathcal W(\dd\lambda,\dd u),
\end{equation*}
and
\begin{equation*}
 \mathsf Y_t(h)
 :=\int_0^t\int_0^\infty h(\lambda)e^{-\lambda(t-u)}
 \sqrt{f(u)}\mathcal W(\dd\lambda,\dd u).
\end{equation*}
For each $t$, these are bounded linear maps from $\cU_b$ into
$L^2(\Omega)$ with centered Gaussian values.

More precisely, the It\^o isometry gives
\[
 \max\bigl\{\|\mathsf M_t(h)\|_{L^2(\Omega)}^2,
             \|\mathsf Y_t(h)\|_{L^2(\Omega)}^2\bigr\}
 \leq \|h\|_{\cU_b}^2\int_0^t f(u)\dd u.
\]
These scalar evaluations specify the cylindrical processes. The
spectral space $\cU_b$ supplies their test functions, whereas
$\cE_\gamma$ contains the forward curves that form the state of the
transport lift. A $\cU_b$-valued realization of the cylindrical
processes is not needed for this construction.
Define the spectral multiplication operator on $\cU_b$ by
\[
 (\mathsf D_bh)(\lambda):=\lambda h(\lambda),
 \qquad
 \cD(\mathsf D_b)
 :=\{h\in\cU_b:\lambda h(\lambda)\in\cU_b\}.
\]
The operator $\mathsf D_b$ acts on the frequency variable $\lambda$,
while the translation generator $A$ acts on curves in $\cE_\gamma$.
For $h\in\cD(\mathsf D_b)$, the contraction bound for the exponential
multipliers gives
\[
 \int_0^t\left(\int_0^s f(u)
 \|e^{-\cdot(s-u)}\mathsf D_bh\|_{\cU_b}^2\dd u\right)^{1/2}\dd s
 \leq t\|\mathsf D_bh\|_{\cU_b}
       \left(\int_0^t f(u)\dd u\right)^{1/2}<\infty.
\]
Stochastic Fubini \cite[Section~4.5]{DaPratoZabczyk2014}
therefore yields the weak Ornstein--Uhlenbeck identity
\begin{equation*}
 \mathsf Y_t(h)
 =\mathsf M_t(h)-\int_0^t\mathsf Y_s(\mathsf D_bh)\dd s.
\end{equation*}
This identity defines the cylindrical Ornstein--Uhlenbeck dynamics
through their scalar evaluations on $\cU_b$.

Since $\mu_b((0,\infty))=\infty$, the constant function $1$ does not
belong to $\cU_b$. The boundary velocity must therefore be defined
and checked through its space--time stochastic integral below; it
cannot be justified simply by substituting $h=1$ in the cylindrical
identity. The variance calculation below separates the two regularity
regimes.

If $1<b<2$ and $f\in L^\infty_{\mathrm{loc}}(\Rplus)$, the boundary
velocity of Proposition~\ref{prop:velocity} is
\begin{equation*}
 V_t=\sqrt b\int_0^t\int_0^\infty
 e^{-\lambda(t-u)}\sqrt{f(u)}
 \mathcal W(\dd\lambda,\dd u),
\end{equation*}
with variance
\begin{equation}\label{eq:velocity-double-variance}
 \E\abs{V_t}^2
 =b2^{b-2}\int_0^t f(u)(t-u)^{b-2}\dd u.
\end{equation}
Let $0<b\leq1$ and $t>0$.  If there exist $\varepsilon\in(0,t]$ and $c>0$ such
that $f(u)\geq c$ for almost every $u\in(t-\varepsilon,t)$, then the integral
in \eqref{eq:velocity-double-variance} is infinite.  The integrated
representation \eqref{eq:integrated-ou-rep} is well defined for every
$0<b<2$.

Fix $M\in\mathbb N$ and $\Lambda>0$, and choose a deterministic partition
$0=\lambda_0<\lambda_1<\cdots<\lambda_M=\Lambda$.  For
$k=1,\ldots,M$, let
\[
 I_k:=(\lambda_{k-1},\lambda_k],\qquad \xi_k\in I_k,
\]
where $I_k$ is the $k$th spectral cell and $\xi_k$ its deterministic
quadrature node. Write
\[
 \delta_\lambda:=\max_{1\leq k\leq M}
                    (\lambda_k-\lambda_{k-1}).
\]
The partitions need not be uniform. For the spectral rates in $M$, we assume
$\delta_\lambda\leq C_{\mathrm{mesh}}\Lambda/M$, where
$C_{\mathrm{mesh}}\geq1$ is fixed independently of $M$ and $\Lambda$.
Only the largest cell is constrained. The constants in these rates may depend on
$C_{\mathrm{mesh}}$, but are uniform over the admissible partitions
and nodes.

No lower bound on individual cell widths is imposed, so local refinement
is allowed. Define
\begin{equation}\label{eq:factor-weights}
 w_k:=\mu_b(I_k)
 =\frac{\lambda_k^{2-b}-\lambda_{k-1}^{2-b}}
 {(2-b)\Gamma(2-b)},
 \qquad
 B_t^k:=w_k^{-1/2}\mathcal W(I_k\times[0,t]).
\end{equation}
The normalized indicators $w_k^{-1/2}\mathbf 1_{I_k}$ are
orthonormal in $\cU_b$. Accordingly, $B^1,\ldots,B^M$ are independent
standard Brownian motions.

The weights are exact spectral masses; only the relaxation rates are
approximated. Using the same $\mathcal W$ couples all resolutions
with the original process, although different resolutions need not
be independent. The feature approximation is
\begin{equation}\label{eq:finite-step-feature}
 \phi_r^{M,\Lambda}(\lambda)
 :=\sqrt b\sum_{k=1}^M k_{\xi_k}(r)\mathbf 1_{I_k}(\lambda),
\end{equation}
with value zero above $\Lambda$.
The $k$th factor is the unique solution of
\begin{equation*}
 Y_t^{k,M,\Lambda}
 =-\xi_k\int_0^tY_s^{k,M,\Lambda}\dd s
  +\int_0^t\sqrt{f(s)}\dd B_s^k,
 \qquad Y_0^{k,M,\Lambda}=0,
\end{equation*}
namely
\begin{equation*}
 Y_t^{k,M,\Lambda}
 =\int_0^t e^{-\xi_k(t-s)}\sqrt{f(s)}\dd B_s^k.
\end{equation*}
Define the aggregate factor process
\begin{equation*}
 \mathcal Y_t^{M,\Lambda}
 :=\sqrt b\sum_{k=1}^M\sqrt{w_k}\,Y_t^{k,M,\Lambda}
\end{equation*}
and the finite-factor driver
\begin{equation*}
 \zeta_t^{M,\Lambda}:=\int_0^t\mathcal Y_s^{M,\Lambda}\dd s.
\end{equation*}
Since $f\in L^1(0,t)$ and the exponential kernels are bounded by one,
stochastic Fubini gives
\begin{equation}\label{eq:finite-integrated-feature}
 \zeta_t^{M,\Lambda}
 =\sqrt b\sum_{k=1}^M\sqrt{w_k}
   \int_0^t k_{\xi_k}(t-u)\sqrt{f(u)}\dd B_u^k.
\end{equation}
Thus \eqref{eq:finite-integrated-feature} is obtained from
\eqref{eq:integrated-ou-rep} by replacing $\phi$ with
$\phi^{M,\Lambda}$: the spectral tail is discarded and the relaxation
rate on each retained cell is replaced by its node.

The corresponding forward curve is
\[
 \Theta_t^{M,\Lambda}(r)
 :=\sqrt b\sum_{k=1}^M\sqrt{w_k}
   \int_0^t k_{\xi_k}(t+r-u)\sqrt{f(u)}\dd B_u^k.
\]
The identity $k_\lambda(a+r)=k_\lambda(a)+e^{-\lambda a}k_\lambda(r)$
gives
\begin{equation}\label{eq:finite-forward-curve}
 \Theta_t^{M,\Lambda}(r)
 =\zeta_t^{M,\Lambda}
   +\sqrt b\sum_{k=1}^M\sqrt{w_k}\,k_{\xi_k}(r)
          Y_t^{k,M,\Lambda}.
\end{equation}

Its derivative is a finite sum of decaying exponentials, so
$\Theta_t^{M,\Lambda}\in\cE_\gamma$. The relation to a finite
Markovian state is seen directly by splitting the future driver
increment at the current time $t$: for $r\geq0$,
\[
 \zeta_{t+r}^{M,\Lambda}-\zeta_t^{M,\Lambda}
 =\sqrt b\sum_{k=1}^M\sqrt{w_k}
 \left(k_{\xi_k}(r)Y_t^{k,M,\Lambda}
 +\int_t^{t+r}k_{\xi_k}(t+r-u)\sqrt{f(u)}\dd B_u^k\right).
\]
The first term is determined by the current factors and equals
$\Theta_t^{M,\Lambda}(r)-\Theta_t^{M,\Lambda}(0)$; the second uses
only fresh Brownian increments. Thus the $M$ factors retain the
information from the past needed to evolve the driver. Since the
scalar equation uses driver increments, its current value together
with these factors is sufficient to restart the coupled system.

For $t>0$ and $g\in C([0,t])$,
\begin{equation*}
 \int_0^t g(s)\dd\zeta_s^{M,\Lambda}
 =\int_0^t g(s)\mathcal Y_s^{M,\Lambda}\dd s.
\end{equation*}
If the coefficients are continuous and globally Lipschitz in $x$ uniformly
on finite time strips, the finite-factor equation
\begin{equation*}
 X_t^{M,\Lambda}
 =x_0+\int_0^t\drift(X_s^{M,\Lambda},s)\dd s
 +\int_0^t\sigma(X_s^{M,\Lambda},s)
          \mathcal Y_s^{M,\Lambda}\dd s
\end{equation*}
has a unique global solution. These conditions hold under
Assumption~\ref{ass:coefficients-mf}. Set
\begin{equation*}
 \mathbf Y_t^{M,\Lambda}
 :=\bigl(Y_t^{1,M,\Lambda},\ldots,Y_t^{M,M,\Lambda}\bigr)\in\R^M.
\end{equation*}
Then $(X^{M,\Lambda},\mathbf Y^{M,\Lambda})$ is a
time-inhomogeneous Markov process on $\R\times\R^M$ and is $C_b$-Feller.
The coupling argument following Proposition~\ref{prop:ou-weak-coefficients}
establishes this property for arbitrary initial states.

\subsection{The finite-variation rate}

For $1<b<2$, both the boundary process and its finite-factor
approximation are absolutely continuous. Their difference has
derivative $V-\mathcal Y^{M,\Lambda}$, and hence
\[
 \TV_{[0,T]}(\zeta_{f,b}-\zeta^{M,\Lambda})
 =\int_0^T|V_t-\mathcal Y_t^{M,\Lambda}|\dd t
 \leq\sqrt T\,\|V-\mathcal Y^{M,\Lambda}\|_{L^2(0,T)}.
\]
An $L^2(0,T)$ estimate for the velocity therefore gives the
total-variation control required by the nonlinear equation. The next
result first displays the separate quadrature and tail contributions,
then balances them to obtain a rate in the number of factors.

\begin{theorem}\label{thm:finite-variation-rate}
Let $1<b<2$, set $\alpha:=2-b\in(0,1)$, and let
$T>0$, $M\in\mathbb N$, $\Lambda>0$, and $f\in L^\infty(0,T)$.
For a spectral partition of $[0,\Lambda]$ with maximum mesh width
$\delta_\lambda$, one has
\begin{align}
 \E\int_0^T\abs{V_t-\mathcal Y_t^{M,\Lambda}}^2\dd t
 \leq{}&
 \frac{b\norm{f}_{\infty,T} T^4}{12\alpha\Gamma(\alpha)}
 \delta_\lambda^2\Lambda^\alpha
 \label{eq:finite-variation-general}\\*
 &+\frac{b\norm{f}_{\infty,T} T}
 {2(1-\alpha)\Gamma(\alpha)}\Lambda^{\alpha-1}.
 \nonumber
\end{align}
Set $\Lambda=M^{2/3}$ and assume
$\delta_\lambda\leq C_{\mathrm{mesh}}M^{-1/3}$.
Then, for every $p\geq2$,
\begin{equation}\label{eq:finite-variation-velocity-rate}
 \left(\E\norm{V-\mathcal Y^{M,\Lambda}}_{L^2(0,T)}^p\right)^{1/p}
 \leq C M^{-(b-1)/3}
\end{equation}
and
\begin{equation}\label{eq:finite-variation-driver-rate}
 \left(\E\sup_{t\leq T}
 \abs{\zeta_{t,f,b}-\zeta_t^{M,\Lambda}}^p\right)^{1/p}
 \leq C M^{-(b-1)/3}.
\end{equation}
If, in addition, Assumption~\ref{ass:coefficients-mf} holds, the coupled
solutions with the same deterministic initial value $x_0$ satisfy
\begin{equation}\label{eq:finite-variation-X-rate}
 \left(\E\sup_{t\leq T}
 \abs{X_t-X_t^{M,\Lambda}}^p\right)^{1/p}
 \leq C M^{-(b-1)/3}.
\end{equation}
\end{theorem}

\begin{proof}
Under the coupling induced by $\mathcal W$, the It\^o isometry gives
\begin{align*}
 \E\abs{V_t-\mathcal Y_t^{M,\Lambda}}^2
 ={}&b\int_0^tf(u)\sum_{k=1}^M\int_{I_k}
 \abs{e^{-\lambda(t-u)}-e^{-\xi_k(t-u)}}^2
 \mu_b(\dd\lambda)\dd u\\
 &+b\int_0^tf(u)\int_\Lambda^\infty
 e^{-2\lambda(t-u)}\mu_b(\dd\lambda)\dd u.
\end{align*}
For the quadrature error, use
$\abs{e^{-\lambda r}-e^{-\xi r}}\leq r\delta_\lambda$ and
\[
 \mu_b([0,\Lambda])
 =\frac{\Lambda^\alpha}{\alpha\Gamma(\alpha)}.
\]
Integrating over $0\leq u\leq t\leq T$ and using
$\int_0^T\int_0^t(t-u)^2\dd u\dd t=T^4/12$ gives the first term in
\eqref{eq:finite-variation-general}.

For the tail, integrate first in $t$ and use
$\int_0^{T-u}e^{-2\lambda r}\dd r\leq(2\lambda)^{-1}$. Its
contribution is at most
\[
 \frac{b\|f\|_{\infty,T}T}{2\Gamma(\alpha)}
 \int_\Lambda^\infty\lambda^{\alpha-2}\dd\lambda
 =\frac{b\|f\|_{\infty,T}T}
        {2(1-\alpha)\Gamma(\alpha)}\Lambda^{\alpha-1},
\]
where $\alpha<1$ is essential for integrability. Under
$\delta_\lambda\leq C_{\mathrm{mesh}}\Lambda/M$, the quadrature
term is bounded by $CM^{-2}\Lambda^{\alpha+2}$. Equating its order
with that of the tail gives $\Lambda^3=M^2$. For
$\Lambda=M^{2/3}$, both terms in the squared error have order
$M^{-2(b-1)/3}$, which gives the stated rate after taking a square
root. More generally, the bound tends to zero along any sequence
such that
\[
 \Lambda_M\longrightarrow\infty,
 \qquad \delta_{\lambda,M}^{\,2}\Lambda_M^{2-b}\longrightarrow0.
\]
The powers of $M$ above balance the two terms in this bound; they do not
assert an optimal rate over other spectral approximations.
For higher moments, Minkowski's
inequality gives, for any centered Gaussian process $G$ with
$\E\int_0^T|G_t|^2\dd t<\infty$ and $p\geq2$,
\[
 \left(\E\|G\|_{L^2(0,T)}^p\right)^{2/p}
 \leq\int_0^T\left(\E|G_t|^p\right)^{2/p}\dd t
 =m_p^2\E\int_0^T|G_t|^2\dd t,
 \qquad m_p:=\left(\E|N|^p\right)^{1/p},
\]
where $N$ is a standard normal random variable. Apply this with
$G=V-\mathcal Y^{M,\Lambda}$ to obtain
\eqref{eq:finite-variation-velocity-rate}.

The total-variation estimate at the start of this subsection also
controls the uniform norm, since the driver error vanishes at zero.
It therefore gives \eqref{eq:finite-variation-driver-rate}.

For the nonlinear equation, put $D_t=\abs{X_t-X_t^{M,\Lambda}}$.
Let $L_{\drift,T}$ and $L_{\sigma,T}$ denote the uniform spatial
Lipschitz constants supplied by Assumption~\ref{ass:coefficients-mf}.
Write
\[
 S_M:=\sup_{t\leq T}|\sigma(X_t^{M,\Lambda},t)|
 \leq K_T+L_T\|X^{M,\Lambda}\|_{\infty;[0,T]}.
\]
Subtracting the equations yields
\begin{align*}
 D_t\leq{}&L_{\drift,T}\int_0^tD_s\dd s
 +L_{\sigma,T}\int_0^t\abs{V_s}D_s\dd s\\
 &+S_M
 \int_0^t\abs{V_s-\mathcal Y_s^{M,\Lambda}}\dd s.
\end{align*}

This decomposition leaves the exact velocity $V$ in the coefficient
of $D_s$. The resulting stability factor is therefore independent of
the spectral discretization, and the approximation enters only
through the final velocity-error term.

Gronwall's inequality with respect to the random finite measure
$L_{\drift,T}\dd s+L_{\sigma,T}\abs{V_s}\dd s$ gives
\[
 \sup_{t\leq T}D_t
 \leq\sqrt{T}S_M
 e^{L_{\drift,T}T+L_{\sigma,T}\norm{V}_{L^1(0,T)}}
 \norm{V-\mathcal Y^{M,\Lambda}}_{L^2(0,T)}.
 \]
By Proposition~\ref{prop:velocity},
$\E\|V\|_{L^2(0,T)}^2<\infty$. For $g\in L^2(0,T)$, stochastic Fubini writes
$\int_0^Tg(t)V_t\dd t$ as a scalar Gaussian stochastic integral
with deterministic kernel
\[
 K_g(\lambda,u)
 :=\sqrt{b f(u)}\int_u^T g(t)e^{-\lambda(t-u)}\dd t.
\]
The Fubini condition follows from
\[
 \int_0^T |g(t)|\bigl(\E|V_t|^2\bigr)^{1/2}\dd t
 \leq\|g\|_{L^2(0,T)}
       \bigl(\E\|V\|_{L^2(0,T)}^2\bigr)^{1/2}<\infty.
\]
Hence $V$ is a centered Gaussian random element of $L^2(0,T)$.
Fernique's theorem
\cite[Section~2.2.1]{DaPratoZabczyk2014} gives
$\E\exp(c\|V\|_{L^2(0,T)}^2)<\infty$ for some $c>0$.
Together with $\|V\|_{L^1}\leq\sqrt T\|V\|_{L^2}$ and
$a y\leq c y^2+a^2/(4c)$ for $a,y\geq0$, this gives moments of every
order for the exponential factor in the Gronwall bound.

It remains to control $S_M$ uniformly. The velocity rate implies
\[
 \sup_{M\geq1}\E\|\mathcal Y^{M,\Lambda}\|_{L^2(0,T)}^2<\infty,
 \qquad \Lambda=M^{2/3}.
\]
These velocities are centered Gaussian random elements of $L^2(0,T)$.
For any such element $U$ whose covariance has trace at most $C_0>0$, its
covariance eigenbasis expansion gives
\[
 \E e^{\varepsilon\|U\|_{L^2}^2}
 \leq\exp\!\left(\frac{\varepsilon C_0}{1-2\varepsilon C_0}\right),
 \qquad 0<\varepsilon<(2C_0)^{-1}.
\]
Indeed, apply $-\log(1-v)\leq v/(1-v)$ to the Gaussian product
formula and then pass to the limit by monotone convergence.
Consequently, the norms
$\|\mathcal Y^{M,\Lambda}\|_{L^1(0,T)}$ have uniform exponential
moments of every linear order. Linear growth and Gronwall give
\[
 1+\|X^{M,\Lambda}\|_{\infty;[0,T]}
 \leq(1+|x_0|)\exp\!\left(C_T
          [T+\|\mathcal Y^{M,\Lambda}\|_{L^1(0,T)}]\right).
\]
Thus $S_M$ has uniformly bounded moments of every finite order.
H\"older's inequality with three factors, using the velocity error,
$S_M$, and the exact-velocity exponential in $L^{3p}$, proves
\eqref{eq:finite-variation-X-rate} with the same rate.
\end{proof}
\subsection{The Young-regime rate}

For $0<b\leq1$, we work directly with increments of the integrated
driver. To apply the pathwise stability estimate below, we need control
in a H\"older norm with exponent greater than $1/2$. We strengthen
the pointwise mean-square estimate first to a uniform bound and then,
by interpolation, to a H\"older bound. For the solution estimates we choose
$1/2<\eta<\eta'<H_b^\circ$: the exponent $\eta$ is used for stability
of the nonlinear equation, whereas $\eta'$ provides the additional
regularity needed for interpolation. The strict upper bound on
$\eta'$ also covers the logarithmic case $b=1$, for which
the available bounds give every H\"older exponent below
$H_b^\circ=1$.

The kernel $k_\lambda$ from Section~\ref{subsec:ou-factors} satisfies
\begin{equation}\label{eq:k-lambda-Lipschitz}
 \abs{\partial_\lambda k_\lambda(r)}
 =\left\lvert\int_0^r v e^{-\lambda v}\dd v\right\rvert
 \leq\frac{r^2}{2}
\end{equation}
for all $\lambda>0$ and $r\geq0$.

For $\eta\in(0,1)$, let $c^\eta([0,T])$ denote the closure of the
piecewise linear paths in $\Holder{\eta}([0,T])$. This is a separable
Banach space with the inherited H\"older norm: piecewise linear paths
with rational values on the dyadic time grids form a dense countable
set.

The use of $c^\eta([0,T])$ concerns measurability: the full H\"older
space is not separable in its usual norm. Every path with H\"older
exponent strictly greater than $\eta$ belongs to $c^\eta([0,T])$,
as the interpolation argument in the proof below shows. Thus the
Gaussian paths under consideration admit strongly measurable
realizations in a separable Banach space without changing the norm
in any of the estimates.

\begin{lemma}\label{lem:uniform-gaussian-holder}
Let $T>0$ and let $(G^m)_{m\in\mathcal I}$ be a nonempty family of
centered continuous Gaussian processes on $[0,T]$ with $G_0^m=0$.
Suppose that, for some $H\in(0,1)$ and $0\leq K<\infty$,
\begin{equation}\label{eq:uniform-gaussian-metric-assumption}
 \sup_{m\in\mathcal I}\E\abs{G_t^m-G_s^m}^2
 \leq K\abs{t-s}^{2H},\qquad s,t\in[0,T].
\end{equation}
Then, for every $\eta\in(0,H)$, there exist constants $c,C>0$, depending only
on $K,H,\eta,T$, such that
\begin{equation}\label{eq:uniform-gaussian-holder-exp}
 \sup_{m\in\mathcal I}
 \E\exp\!\left(c\norm{G^m}_{\Holder{\eta}([0,T])}^2\right)
 \leq C.
\end{equation}
Each $G^m$ has a realization as a Borel measurable random element of
$c^\eta([0,T])$.
\end{lemma}

\begin{proof}
The case $K=0$ is immediate. For $K>0$, let $t_{n,k}=kT2^{-n}$ and
$\Delta_{n,k}G=G(t_{n,k+1})-G(t_{n,k})$.  For a continuous function $x$ with
$x(0)=0$, dyadic chaining gives
\begin{equation}\label{eq:dyadic-holder-control}
 \norm{x}_{\Holder{\eta}([0,T])}
 \leq C_{\eta,T}\sup_{n\geq0}
 2^{n\eta}\max_{0\leq k<2^n}\abs{\Delta_{n,k}x}.
\end{equation}
This follows by approximating the endpoints by their nested dyadic
ancestors and summing the geometric series.

Each $\Delta_{n,k}G^m$ has variance at most $KT^{2H}2^{-2nH}$, so
uniformly in $m$,
\begin{align*}
 &\mathbb P\!\left(
 \sup_{n\geq0}2^{n\eta}\max_{k<2^n}
 \abs{\Delta_{n,k}G^m}>x\right)\\
 &\qquad\leq
 2\sum_{n=0}^\infty 2^n
 \exp\!\left(-\frac{x^2 2^{2n(H-\eta)}}{2KT^{2H}}\right).
\end{align*}
For $x\geq1$, splitting each exponential into two equal factors bounds
this sum by $C_1e^{-c_1x^2}$, uniformly in $m$, since
$\sum_{n\geq0}2^n e^{-c2^{2n(H-\eta)}}$ is finite.
Combining this estimate with \eqref{eq:dyadic-holder-control} and integrating
the tail proves \eqref{eq:uniform-gaussian-holder-exp} for sufficiently small
$c>0$. The norm is measurable: for a continuous path its supremum and
H\"older seminorm may both be evaluated over rational times.

Fix $\eta_1\in(\eta,H)$. Applying the preceding estimate at exponent
$\eta_1$ shows that $G^m$ is almost surely $\eta_1$-H\"older. Its dyadic
piecewise linear interpolants $G^{m,n}$ satisfy
\[
 \|G^{m,n}-G^m\|_{\Holder{\eta}([0,T])}
 \leq C_{\eta,\eta_1,T}\,2^{-n(\eta_1-\eta)}
 [G^m]_{\eta_1;[0,T]}.
\]
This follows by estimating increments separately below and above the
mesh size. Each interpolant is a measurable function of finitely many
Gaussian values and belongs to $c^\eta([0,T])$. Its almost sure limit
in that space gives the Borel realization.
\end{proof}

\begin{theorem}
\label{thm:young-driver-rate}
Let $T>0$, $M\in\mathbb N$, $\Lambda>0$, $0<b\leq1$, and
$f\in L^\infty(0,T)$. Use any spectral partition of $[0,\Lambda]$
with maximum mesh width $\delta_\lambda$. Then
\begin{align}
 \sup_{t\leq T}\E\abs{\zeta_{t,f,b}-\zeta_t^{M,\Lambda}}^2
 \leq{}&
 \frac{b\norm{f}_{\infty,T} T^5}
 {20(2-b)\Gamma(2-b)}\delta_\lambda^2\Lambda^{2-b}
 \label{eq:young-pointwise-general}\\*
 &+\frac{\norm{f}_{\infty,T} T}{\Gamma(2-b)}\Lambda^{-b}.
 \nonumber
\end{align}
Set $\Lambda=M^{1/2}$ and assume
$\delta_\lambda\leq C_{\mathrm{mesh}}M^{-1/2}$. Then
\begin{equation}\label{eq:young-pointwise-rate}
 \sup_{t\leq T}
 \norm{\zeta_{t,f,b}-\zeta_t^{M,\Lambda}}_{L^2(\Omega)}
 \leq C M^{-b/4}.
\end{equation}
Moreover, for every $p\geq2$,
\begin{equation}\label{eq:young-uniform-rate}
 \left(\E\sup_{t\leq T}
 \abs{\zeta_{t,f,b}-\zeta_t^{M,\Lambda}}^p\right)^{1/p}
 \leq C_p M^{-b/4}\sqrt{\log(2+M)}.
\end{equation}
For every $p\geq2$ and $1/2<\eta<\eta'<H_b^\circ$,
\begin{equation}\label{eq:young-holder-rate}
 \left(\E\norm{\zeta_{f,b}-\zeta^{M,\Lambda}}_{\Holder{\eta}([0,T])}^p
 \right)^{1/p}
 \leq C_{p,\eta,\eta'}
 \left(M^{-b/4}\sqrt{\log(2+M)}\right)^{1-\eta/\eta'}.
\end{equation}
\end{theorem}

\begin{proof}
Let $e_t^{M,\Lambda}:=\zeta_{t,f,b}-\zeta_t^{M,\Lambda}$.  Under the common
orthogonal Gaussian random measure $\mathcal W$,
\begin{align*}
 \E\abs{e_t^{M,\Lambda}}^2
 ={}&b\int_0^tf(u)\sum_{k=1}^M\int_{I_k}
 \abs{k_\lambda(t-u)-k_{\xi_k}(t-u)}^2
 \mu_b(\dd\lambda)\dd u\\
 &+b\int_0^tf(u)\int_\Lambda^\infty
 k_\lambda(t-u)^2\mu_b(\dd\lambda)\dd u.
\end{align*}
By \eqref{eq:k-lambda-Lipschitz}, the first term is bounded by
\[
 \frac{b\norm{f}_{\infty,T}\delta_\lambda^2}{4}\mu_b([0,\Lambda])
 \int_0^t(t-u)^4\dd u,
\]
which is the first term in \eqref{eq:young-pointwise-general}.  Since
$k_\lambda(r)\leq1/\lambda$, the tail is bounded by
\[
 b\norm{f}_{\infty,T} T\int_\Lambda^\infty\lambda^{-2}
 \mu_b(\dd\lambda)
 =\frac{\norm{f}_{\infty,T} T}{\Gamma(2-b)}\Lambda^{-b}.
\]
Under $\delta_\lambda\leq C_{\mathrm{mesh}}\Lambda/M$, balancing
the resulting powers $M^{-2}\Lambda^{4-b}$ and $\Lambda^{-b}$
gives $\Lambda=M^{1/2}$ and \eqref{eq:young-pointwise-rate}.

More generally, \eqref{eq:young-pointwise-general} gives uniform
pointwise mean-square convergence whenever $\Lambda_M\to\infty$ and
$\delta_{\lambda,M}^{\,2}\Lambda_M^{2-b}\to0$. The balanced choice
above gives the explicit rate in $M$.

For the uniform estimate, let
\[
 d_M(s,t)^2:=\E\abs{e_t^{M,\Lambda}-e_s^{M,\Lambda}}^2,
\]
and write $h=t-s>0$. The spectral tail contribution is bounded by the
full increment variance, hence by \eqref{eq:increment-bound} for $b<1$
and by $C_{H,T}h^{2H}$ for every fixed $H\in(0,1)$ when $b=1$, by
Proposition~\ref{prop:critical-regularity}. For the quadrature error,
put $r=s-u$. Then
\[
 \abs{(k_\lambda-k_\xi)(r+h)-(k_\lambda-k_\xi)(r)}
 \leq\delta_\lambda\int_r^{r+h}v\dd v
 \leq\delta_\lambda h(r+h),
\]
and, on the interval $(s,t]$,
$\abs{k_\lambda(v)-k_\xi(v)}\leq\delta_\lambda v^2/2$.
The contribution from $u\in[0,s]$ is therefore bounded by
\begin{align*}
 &b\norm{f}_{\infty,T}\delta_\lambda^2
   \mu_b([0,\Lambda])h^2\int_0^s(s-u+h)^2\dd u\\
 &\qquad\leq C_T\delta_\lambda^2
   \mu_b([0,\Lambda])h^2,
\end{align*}
and the contribution from $u\in(s,t]$ is at most
\[
 \frac{b\norm{f}_{\infty,T}\delta_\lambda^2}{4}
 \mu_b([0,\Lambda])\int_s^t(t-u)^4\dd u
 =\frac{b\norm{f}_{\infty,T}}{20}
 \delta_\lambda^2\mu_b([0,\Lambda])h^5.
\]
For $\Lambda=M^{1/2}$ and
$\delta_\lambda\leq C_{\mathrm{mesh}}M^{-1/2}$,
\[
 \delta_\lambda^2\mu_b([0,\Lambda])
 \leq\frac{C_{\mathrm{mesh}}^2 M^{-b/2}}
 {(2-b)\Gamma(2-b)}\leq C_{b,C_{\mathrm{mesh}}}.
\]
If $b<1$, then $h^2\leq T^{1-b}h^{1+b}$ and
$h^5\leq T^{4-b}h^{1+b}$, so adding the tail contribution proves
\begin{equation}\label{eq:error-canonical-metric}
 d_M(s,t)^2\leq C\abs{t-s}^{1+b}
\end{equation}
uniformly in $M$.  If $b=1$, then for every fixed $H<1$ the same calculation,
using $h^2+h^5\leq C_{H,T}h^{2H}$, gives
\begin{equation}\label{eq:critical-error-canonical-metric}
 d_M(s,t)^2\leq C_{H,T}\abs{t-s}^{2H}.
\end{equation}
By \eqref{eq:young-pointwise-rate}, the canonical diameter is
$O(M^{-b/4})$. Take $H=(1+b)/2$ for $b<1$ and any fixed $H\in(0,1)$
for $b=1$, with $H>\eta'$ for the H\"older estimate. For $\varepsilon>0$,
let
$\mathcal N([0,T],d_M,\varepsilon)$ denote the minimum number of open
$d_M$-balls of radius $\varepsilon$ required to cover $[0,T]$. The
increment bounds give
\[
 \mathcal N([0,T],d_M,\varepsilon)
 \leq 1+C\varepsilon^{-1/H}.
\]
Apply Dudley's entropy bound
\cite[Theorem~1.3.3]{AdlerTaylor2007} to the centered Gaussian process
$(t,\kappa)\mapsto\kappa e_t^{M,\Lambda}$ on
$[0,T]\times\{-1,1\}$. Its covering number is at most
$2\mathcal N([0,T],d_M,\varepsilon)$ and its canonical diameter is
at most $a_M:=CM^{-b/4}$, after enlarging $C$.  Therefore,
\[
 \E\sup_{t\leq T}\abs{e_t^{M,\Lambda}}
 \leq C\int_0^{a_M}
 \sqrt{\log(1+C\varepsilon^{-1/H})}\dd\varepsilon
 \leq C M^{-b/4}\sqrt{\log(2+M)}.
\]
For the last inequality, substitute $\varepsilon=a_Mu$.  The logarithm is
bounded by $C\{\log(2+M)+\log(1/u)\}$, and
$\int_0^1\sqrt{\log(1/u)}\dd u<\infty$.
Moreover,
$\sup_{t\leq T}\|e_t^{M,\Lambda}\|_{L^2(\Omega)}\leq CM^{-b/4}$.
The Borell--Tsirelson--Ibragimov--Sudakov inequality
\cite[Theorem~2.1.1]{AdlerTaylor2007} yields, for every $p\geq2$,
\[
 \left\|\sup_{t\leq T}\abs{e_t^{M,\Lambda}}\right\|_{L^p(\Omega)}
 \leq \E\sup_{t\leq T}\abs{e_t^{M,\Lambda}}+C_pM^{-b/4},
\]
after integrating its Gaussian tail bound.  This proves
\eqref{eq:young-uniform-rate}.

Finally, the increment bounds and
Lemma~\ref{lem:uniform-gaussian-holder} give exponential moments of
$\|e^{M,\Lambda}\|_{\Holder{\eta'}([0,T])}$, uniformly in $M$, for
every $\eta'<H_b^\circ$. To see how the uniform error transfers to
the weaker H\"older norm, observe that
\[
 |x_t-x_s|\leq
 \min\bigl(2\|x\|_{\infty;[0,T]},
           [x]_{\eta';[0,T]}|t-s|^{\eta'}\bigr).
\]
Interpolating between these bounds, and putting
$\vartheta:=\eta/\eta'\in(0,1)$, gives
\[
 \norm{x}_{\Holder{\eta}([0,T])}
 \leq C\norm{x}_{\infty;[0,T]}^{1-\vartheta}
       \norm{x}_{\Holder{\eta'}([0,T])}^{\vartheta}.
\]
Consequently, H\"older's inequality yields
\begin{align*}
 \bigl\|\|e^{M,\Lambda}\|_{\Holder{\eta}}\bigr\|_{L^p(\Omega)}
 &\leq C
 \bigl\|\|e^{M,\Lambda}\|_\infty\bigr\|_{L^p(\Omega)}^{1-\vartheta}
 \bigl\|\|e^{M,\Lambda}\|_{\Holder{\eta'}}\bigr\|_{L^p(\Omega)}^{\vartheta}
 \\
 &\leq C_{p,\eta,\eta'}
 \left(M^{-b/4}\sqrt{\log(2+M)}\right)^{1-\vartheta}.
\end{align*}
Here all path norms are taken on $[0,T]$. This proves
\eqref{eq:young-holder-rate}. The factor $1-\eta/\eta'$ in the
exponent comes from interpolation with the stronger H\"older norm.

\end{proof}

\begin{assumption}
\label{ass:young-quantitative}
The coefficients $\drift,\sigma$ are continuous, and $\sigma$ is
continuously differentiable in $x$. For every $T>0$, there are constants
$0\leq K_T,L_T<\infty$ such that
\[
 \begin{aligned}
 |\drift(0,t)|&\leq K_T,\\
 |\drift(x,t)-\drift(y,t)|&\leq L_T|x-y|,\\
 |\sigma(x,t)|+|\partial_x\sigma(x,t)|&\leq K_T,\\
 |\sigma(x,t)-\sigma(y,s)|&\leq L_T(|x-y|+|t-s|),\\
 |\partial_x\sigma(x,t)-\partial_x\sigma(y,s)|
    &\leq L_T(|x-y|+|t-s|),
 \end{aligned}
 \qquad x,y\in\R,\quad s,t\in[0,T].
\]
\end{assumption}

Assumption~\ref{ass:young-quantitative} implies
Assumption~\ref{ass:coefficients-mf} with $\alpha=1$.
We next give the dependence of the solution error on the driver norms.

\begin{lemma}\label{lem:young-solution-stability}
Let $T>0$, $\eta\in(1/2,1)$, and suppose that
Assumption~\ref{ass:young-quantitative} holds.  For
$z,\widetilde z\in\Holder{\eta}([0,T])$ with
$z(0)=\widetilde z(0)=0$, let $x$ and $\widetilde x$ solve
\[
 x_t=x_0+\int_0^t\drift(x_s,s)\dd s
       +\int_0^t\sigma(x_s,s)\dd z_s
\]
and the corresponding equation with initial value $\widetilde x_0$ and
driver $\widetilde z$.  Put $q:=1/\eta$.  Then
\begin{align}
 \norm{x-\widetilde x}_{\infty;[0,T]}
 \leq{}&C\exp\!\left(C\left[1+
 \norm{z}_{\Holder{\eta}([0,T])}^{q}
 +\norm{\widetilde z}_{\Holder{\eta}([0,T])}^{q}\right]\right)
 \notag\\*
 &\times\left(
 \abs{x_0-\widetilde x_0}
 +\norm{z-\widetilde z}_{\Holder{\eta}([0,T])}\right),
 \label{eq:young-stability}
\end{align}
where $C$ depends only on $T$, $\eta$, $A_T=\sup_{t\leq T}|\drift(0,t)|$,
the spatial Lipschitz constant of $\drift$, the uniform bounds and
Lipschitz constants of $\sigma$ and $\partial_x\sigma$ in
Assumption~\ref{ass:young-quantitative}, and $|x_0|+|\widetilde x_0|$.
If the drift is bounded on
$\R\times[0,T]$, $C$ can be chosen independently of the initial values,
with additional dependence on
$\sup_{(x,t)\in\R\times[0,T]}|\drift(x,t)|$.
\end{lemma}

The proof uses a partition on which the local stability terms can be
absorbed. Its number of intervals is of order at most
$1+[z]_\eta^{1/\eta}+[\widetilde z]_\eta^{1/\eta}$ for fixed initial
values. This partition is used only in the proof.

\begin{proof}
\textit{Step 1: A priori bounds.}
Fix $p_0\in(1/\eta,2)$, chosen as a function of $\eta$.
Proposition~\ref{prop:weak-coefficient-young},
with $\alpha=1$, gives solutions of finite $p_0$-variation. The
Young variation estimate \cite[Lemma~2.4, Eq.~(2.8)]{CongDucHong2018} and the
Lipschitz bounds on $\sigma$ imply, on any interval $I=[a,c]$ of
length $\Delta\leq1$, with $K:=[z]_{\eta;[0,T]}$,
\[
 \|x\|_{p_0\text{-var};I}
 \leq C\Delta(1+|x_a|)
       +C(\Delta+K\Delta^\eta)\|x\|_{p_0\text{-var};I}
       +CK\Delta^\eta.
\]
Here the linear-growth drift contributes at most
$C\Delta(1+|x_a|+\|x\|_{p_0\text{-var};I})$, and we used
$\|z\|_{p_0\text{-var};I}\leq[z]_{\eta;[0,T]}\Delta^\eta$ and
$\|\sigma(x_\cdot,\cdot)\|_{p_0\text{-var};I}
\leq C(\|x\|_{p_0\text{-var};I}+\Delta)$.
For $\Delta\leq c_*(1+K)^{-1/\eta}$, with sufficiently small
deterministic $c_*>0$, absorption yields
\[
 \|x\|_{p_0\text{-var};I}
 \leq C\Delta(1+|x_a|)+CK\Delta^\eta.
\]
On an equal partition of mesh $h$ comparable to
$\min\{T,1,c_*(1+K)^{-1/\eta}\}$, the endpoint values therefore satisfy
\[
 1+|x_{t_{j+1}}|
 \leq(1+Ch)(1+|x_{t_j}|)+CKh^\eta.
\]
Discrete Gronwall, using $\prod_j(1+Ch)\leq e^{CT}$ and
$\sum_jKh^\eta\leq C_T(1+K^{1/\eta})$, gives
\[
 \|x\|_{\infty;[0,T]}
 \leq C_T(1+|x_0|+K^{1/\eta}).
\]
The same bound holds for $\widetilde x$. Put
\[
 R:=\left(1+|x_0|+|\widetilde x_0|
       +(1+[z]_{\eta;[0,T]}+[\widetilde z]_{\eta;[0,T]})^{1/\eta}
       \right)^\eta.
\]
Both solution norms are bounded by $CR^{1/\eta}$, and
$K\le R$. If $I$ has length
$\Delta\leq\min\{1,c_*R^{-1/\eta}\}$, the absorbed estimate on
each $[u,v]\subset I$ gives
\[
 |x_v-x_u|
 \le C R^{1/\eta}|v-u|+CR|v-u|^\eta.
\]
The same estimate holds for $\widetilde x$. Dividing by
$|v-u|^\eta$ and taking the supremum proves the local H\"older
regularity and gives
\begin{equation}\label{eq:young-local-solution-bound}
 \Delta^\eta\bigl([x]_{\eta;I}+[\widetilde x]_{\eta;I}\bigr)
 \leq C(R^{1/\eta}\Delta+R\Delta^\eta)
 \leq C R\Delta^\eta.
\end{equation}
The last inequality follows from
$R^{1/\eta-1}\Delta^{1-\eta}
=(R^{1/\eta}\Delta)^{1-\eta}\leq c_*^{1-\eta}$.

\medskip
\noindent\textit{Step 2: Local difference estimate.}
Write $e=x-\widetilde x$,
$h=z-\widetilde z$, and introduce the scaled norm
\[
 \|v\|_{\eta,I}^{*}:=\|v\|_{\infty;I}
                         +\Delta^\eta[v]_{\eta;I}.
\]
The Young--Loeve inequality \cite[Eq.~(2.6)]{CongDucHong2018} gives
\[
 \left|\int_s^t u_r\dd v_r-u_s(v_t-v_s)\right|
 \leq C_\eta[u]_{\eta;I}[v]_{\eta;I}|t-s|^{2\eta},
 \qquad a\leq s<t\leq c.
\]
Consequently, for $J(t)=\int_a^t u_r\dd v_r$,
\[
 \|J\|_{\eta,I}^{*}
 \leq C_\eta[v]_{\eta;I}\Delta^\eta
       \bigl(\|u\|_{\infty;I}+\Delta^\eta[u]_{\eta;I}\bigr).
\]

Split the difference of the Young integrals into
\[
 \int_a^\cdot\bigl(\sigma(x_r,r)-\sigma(\widetilde x_r,r)\bigr)\dd z_r
 +\int_a^\cdot\sigma(\widetilde x_r,r)\dd h_r.
\]
To estimate the first integrand, write
\[
 \sigma(x_r,r)-\sigma(\widetilde x_r,r)=A_re_r,
 \qquad
 A_r:=\int_0^1\partial_x\sigma
       ((1-\theta)\widetilde x_r+\theta x_r,r)\dd\theta.
\]
The coefficient assumptions imply $\|A\|_{\infty;I}\le C$ and
\[
 [A]_{\eta;I}\le C\bigl([x]_{\eta;I}
          +[\widetilde x]_{\eta;I}+\Delta^{1-\eta}\bigr).
\]
The product increment identity for $Ae$ therefore gives
\begin{align*}
 &[\sigma(x_\cdot,\cdot)-\sigma(\widetilde x_\cdot,\cdot)]_{\eta;I}\\
 &\qquad\leq C[e]_{\eta;I}
   +C\|e\|_{\infty;I}
        \bigl([x]_{\eta;I}+[\widetilde x]_{\eta;I}
                          +\Delta^{1-\eta}\bigr).
\end{align*}
Together with \eqref{eq:young-local-solution-bound}, this bounds the
first Young integral in the scaled norm by
\[
 C\bigl(R\Delta^\eta+(R\Delta^\eta)^2\bigr)
 \|e\|_{\eta,I}^{*}.
\]
The second integral is bounded by
\[
 C\Delta^\eta(1+R\Delta^\eta)
 \|h\|_{\Holder{\eta}([0,T])}.
\]
Finally, spatial Lipschitz continuity of the drift bounds the scaled
norm of its difference integral by $C\Delta\|e\|_{\eta,I}^{*}$.
Combining these estimates gives
\begin{align}
 \|e\|_{\eta,I}^{*}
 \leq{}&|e_a|+C\Delta^\eta(1+R\Delta^\eta)
          \|h\|_{\Holder{\eta}([0,T])}\notag\\
 &+C\bigl(\Delta+R\Delta^\eta+(R\Delta^\eta)^2\bigr)
       \|e\|_{\eta,I}^{*}.
 \label{eq:young-local-stability}
\end{align}

Choose $0<c_*\le1$, depending only on the coefficient bounds and
$\eta$, so that
$C(c_*+c_*^\eta+c_*^{2\eta})\le1/2$.
If $\Delta\leq\min(1,c_*R^{-1/\eta})$, then
$R\Delta^\eta\le c_*^\eta$ and the last term in
\eqref{eq:young-local-stability} can be absorbed.
\medskip
\noindent\textit{Step 3: Propagation over the full interval.}
Put $d:=\min(1,T,c_*R^{-1/\eta})$, set $n:=\lceil T/d\rceil$,
and partition $[0,T]$ into $n$ equal intervals. Then
\[
 n\leq C_T(1+R^{1/\eta}),
\]
and the absorbed estimate on each interval gives
\[
 \|e\|_{\infty;[t_j,t_{j+1}]}
 \leq C\bigl(|e_{t_j}|+\|h\|_{\Holder{\eta}([0,T])}\bigr),
 \qquad j=0,\ldots,n-1,
\]
with $C$ independent of $R$ and $j$. Iteration yields
\[
 \|e\|_{\infty;[0,T]}
 \leq C^n(1+n)
        \bigl(|e_0|+\|h\|_{\Holder{\eta}([0,T])}\bigr).
\]
Since $q=1/\eta$ and
$R^q\leq C_q(1+|x_0|+|\widetilde x_0|+\|z\|_{\Holder{\eta}}^q
                   +\|\widetilde z\|_{\Holder{\eta}}^q)$,
this is \eqref{eq:young-stability}, with the deterministic initial-value
factor included in $C$.

If the drift is bounded, its variation contribution in the first estimate
is $C\Delta$, independently of $x_a$. Absorption then gives
$\|x\|_{p_0\text{-var};I}\leq C(\Delta+K\Delta^\eta)$.
In this case one can use
$R=1+[z]_{\eta;[0,T]}+[\widetilde z]_{\eta;[0,T]}$ throughout the
local H\"older and stability estimates. Both
\eqref{eq:young-local-solution-bound} and \eqref{eq:young-stability}
are then uniform in the initial values.
\end{proof}

\begin{theorem}\label{thm:young-X-rate}
Let $T>0$ and $M\in\mathbb N$, set $\Lambda=M^{1/2}$, and assume
$0<b\leq1$, $f\in L^\infty(0,T)$, and
Assumption~\ref{ass:young-quantitative}. Use a spectral partition with
$\delta_\lambda\leq C_{\mathrm{mesh}}M^{-1/2}$ and arbitrary
nodes $\xi_k\in I_k$. Let $X$ and $X^{M,\Lambda}$ have the same
deterministic initial value $x_0$ and be coupled through $\mathcal W$. For every
$p\geq2$ and every
$1/2<\eta<\eta'<H_b^\circ$,
\begin{equation}\label{eq:young-X-rate}
 \left(\E\sup_{t\leq T}\abs{X_t-X_t^{M,\Lambda}}^p\right)^{1/p}
 \leq C
 \left(M^{-b/4}\sqrt{\log(2+M)}\right)^{1-\eta/\eta'}.
\end{equation}
The constant may depend on $x_0$ but is independent of $M$.
\end{theorem}

\begin{proof}
Apply Lemma~\ref{lem:young-solution-stability} with $q=1/\eta\in(1,2)$,
$z=\zeta_{f,b}$, and $\widetilde z=\zeta^{M,\Lambda}$.
The increment bounds in Propositions~\ref{prop:young-regularity}
and~\ref{prop:critical-regularity}, together with
\eqref{eq:error-canonical-metric}--\eqref{eq:critical-error-canonical-metric}
and Lemma~\ref{lem:uniform-gaussian-holder}, give $c_0>0$ such that
\begin{equation}\label{eq:uniform-young-exponential-integrability}
 \E e^{c_0\norm{\zeta_{f,b}}_{\Holder{\eta}([0,T])}^2}<\infty,
 \qquad
 \sup_{M\geq1}\E e^{c_0
 \norm{\zeta_{f,b}-\zeta^{M,\Lambda}}_{\Holder{\eta}([0,T])}^2}<\infty.
\end{equation}
Since
$\zeta^{M,\Lambda}=\zeta_{f,b}-(\zeta_{f,b}-\zeta^{M,\Lambda})$, H\"older's
inequality and a reduction of $c_0$ yield
\[
 \sup_{M\geq1}\E e^{c_0
 \norm{\zeta^{M,\Lambda}}_{\Holder{\eta}([0,T])}^2}<\infty.
\]
For $q<2$ and $1\leq r<\infty$, the inequality
$x^q\leq\varepsilon x^2+C_{q,\varepsilon}$, with $\varepsilon$ chosen in
terms of $r$ and $c_0$, now gives moments of order $r$ for the exponential
factor in the stability estimate, uniformly in $M$.

Denote this factor by $K_M$. For the common initial value,
\eqref{eq:young-stability} and H\"older's inequality give
\[
 \bigl\|\|X-X^{M,\Lambda}\|_\infty\bigr\|_{L^p(\Omega)}
 \leq \|K_M\|_{L^{2p}(\Omega)}
       \bigl\|\|\zeta_{f,b}-\zeta^{M,\Lambda}\|_{\Holder{\eta}}
                   \bigr\|_{L^{2p}(\Omega)}.
\]
The first factor is bounded independently of $M$, and the second
satisfies \eqref{eq:young-holder-rate} with moment order $2p$.
This proves \eqref{eq:young-X-rate}.

\end{proof}

\begin{remark*}
Linear multiplicative diffusion also admits a direct quantitative bound.
Use the balanced spectral assumptions of Theorem~\ref{thm:young-driver-rate}:
$\Lambda=M^{1/2}$ and
$\delta_\lambda\leq C_{\mathrm{mesh}}M^{-1/2}$.
For constants $\nu,\beta\in\R$, take
$\drift(x,t)=\nu x$ and $\sigma(x,t)=\beta x$.
The explicit solutions
$X_t=x_0e^{\nu t+\beta\zeta_{t,f,b}}$ and
$X_t^{M,\Lambda}=x_0e^{\nu t+\beta\zeta_t^{M,\Lambda}}$ give
\[
 \|X-X^{M,\Lambda}\|_\infty
 \leq |x_0\beta|e^{|\nu|T+
           |\beta|\max\{\|\zeta_{f,b}\|_\infty,
                         \|\zeta^{M,\Lambda}\|_\infty\}}
             \|\zeta_{f,b}-\zeta^{M,\Lambda}\|_\infty.
\]
Uniform Gaussian exponential bounds and H\"older's inequality therefore
yield, for $p\geq2$,
\[
 \left(\E\|X-X^{M,\Lambda}\|_\infty^p\right)^{1/p}
 \leq C M^{-b/4}\sqrt{\log(2+M)}.
\]
This uses the explicit linear solution, rather than the bounded-diffusion
stability estimate of Lemma~\ref{lem:young-solution-stability}.
\end{remark*}

\begin{remark}
Our bounds use the exact cell masses and the maximum spectral cell
width, allowing nonuniform partitions. For other quadrature rules,
see the geometric partitions in \cite[p.~99]{CarmonaCoutin1998} and
\cite[Section~4]{Harms2021}, Gaussian quadrature for the fractional
kernel in \cite[Section~2.1]{BayerBreneis2023}, and the multifactor
approximations of completely monotone kernels in
\cite[Sections~3--4]{AlfonsiKebaier2024}.
\end{remark}

The preceding solution rate uses the joint Lipschitz regularity of
$\partial_x\sigma$ to obtain an explicit stability constant. Under
the weaker Markov--Feller assumptions, continuity in the
$p$-variation topology gives convergence in probability of the solutions.
Uniform moment bounds then give convergence in every finite $L^r$.

\begin{proposition}
\label{prop:ou-weak-coefficients}
Let $0<b<2$, $T>0$, and $f\in L^\infty(0,T)$, and assume
Assumption~\ref{ass:coefficients-mf}. Couple the original driver and
the finite OU drivers using the same orthogonal Gaussian random measure.
Set
\[
 \Lambda_M=\begin{cases}M^{1/2},&0<b\leq1,\\M^{2/3},&1<b<2.\end{cases}
\]
Use spectral partitions whose maximum cell width satisfies
$\delta_{\lambda,M}\leq C_{\mathrm{mesh}}\Lambda_M/M$ with a
fixed $C_{\mathrm{mesh}}\geq1$. Let the corresponding nonlinear
solutions have the same deterministic
initial value $x_0$. Then, for every $1\leq r<\infty$,
\begin{equation}\label{eq:ou-weak-coefficient-convergence}
 \lim_{M\to\infty}
 \E\!\left[\sup_{0\leq t\leq T}
 |X_t^{M,\Lambda_M}-X_t|^r\right]=0.
\end{equation}
In particular, convergence is uniform in time in probability. For
$1<b<2$, the quantitative bound \eqref{eq:finite-variation-X-rate}
holds.
\end{proposition}
\begin{proof}
First suppose $0<b\leq1$. Choose
\[
 (H_b^\circ)^{-1}<p<\min(2,1+\alpha),\qquad
 p^{-1}<\eta<\eta'<H_b^\circ,
\]
where $H_b^\circ$ and $\alpha$ are as in
Assumption~\ref{ass:coefficients-mf}.

The coefficient threshold guarantees that these intervals are nonempty.

Lemma~\ref{lem:uniform-gaussian-holder} gives measurable realizations of
the drivers in $c^\eta([0,T])$. The continuous inclusion into
$C^{p\text{-var}}([0,T])$ follows from
\[
 \|z\|_{p\text{-var};[0,T]}
 \leq T^\eta[z]_{\eta;[0,T]},\qquad p\eta>1.
\]
Since the drivers start at zero, \eqref{eq:young-holder-rate}, used at
moment order $2$, implies convergence in probability in the
$p$-variation topology. The continuous
solution map of
Proposition~\ref{prop:weak-coefficient-young} therefore implies
$\|X^{M,\Lambda_M}-X\|_\infty\to0$ in probability.

For the moment estimate, let $z$ be a deterministic $\eta$-H\"older
driver with $z_0=0$, put
$R=1+[z]_{\eta;[0,T]}$ and let $x$ denote its solution.
On an interval $I=[u,v]$ of length $h\leq1$, the Young variation
bound \cite[Lemma~2.4, Eq.~(2.8)]{CongDucHong2018}, applied with both variation
exponents equal to $p<2$, and the spatial and temporal Lipschitz bounds
on $\sigma$ give
\[
 \|x\|_{p\text{-var};I}
 \leq C(h+Rh^\eta)
       \bigl(1+|x_u|+\|x\|_{p\text{-var};I}\bigr).
\]
Here both coefficient values grow at most linearly, and
$\|\sigma(x_\cdot,\cdot)\|_{p\text{-var};I}
\leq C\|x\|_{p\text{-var};I}+Ch$.
For $h\leq\min\{1,cR^{-1/\eta}\}$, with a fixed sufficiently small
$c>0$, absorption gives
\[
 \|x\|_{p\text{-var};I}
 \leq C(h+Rh^\eta)(1+|x_u|).
\]
Take an equal partition with mesh comparable to
$\min\{T,1,cR^{-1/\eta}\}$. Its endpoint values satisfy
\[
 1+|x_{t_{j+1}}|
 \leq\bigl(1+C[h+Rh^\eta]\bigr)(1+|x_{t_j}|).
\]
Since $\prod_j(1+C[h+Rh^\eta])\leq
\exp(C\sum_j[h+Rh^\eta])$ and
$\sum_jRh^\eta\leq C_T(1+R^{1/\eta})$, discrete Gronwall and the
within-cell variation estimate give
\begin{equation}\label{eq:linear-growth-solution-moment}
 \|x\|_{\infty;[0,T]}
 \leq C_T(1+|x_0|)
 \exp\!\left(C_T[1+[z]_{\eta;[0,T]}^{1/\eta}]\right).
\end{equation}
Since $1/\eta<2$, the uniform Gaussian quadratic exponential bounds
from Lemma~\ref{lem:uniform-gaussian-holder} control every fixed power
of this estimate. Thus
\[
 \sup_{M\geq1}\E\|X^{M,\Lambda_M}\|_\infty^q
       +\E\|X\|_\infty^q<\infty
       \qquad\text{for every }1\leq q<\infty.
\]
Taking $q>r$ gives uniform integrability of
$\|X^{M,\Lambda_M}-X\|_\infty^r$; its convergence in probability to
zero now proves \eqref{eq:ou-weak-coefficient-convergence}.

For $1<b<2$, Theorem~\ref{thm:finite-variation-rate} uses linear growth,
spatial Lipschitz continuity and uniform moments of the solutions.
It therefore gives
\eqref{eq:finite-variation-X-rate} for $r\geq2$; the remaining
$1\leq r<2$ follow from the $r=2$ bound.
\end{proof}

\begin{remark}
For fixed $M$ and $\Lambda$, the state
$(X^{M,\Lambda},\mathbf Y^{M,\Lambda})$ is Markov and $C_b$-Feller
under Assumption~\ref{ass:coefficients-mf}. Start the system at a
deterministic time $s$ from
$(x,\mathbf y)\in\R\times\R^M$. For $t\geq s$, each factor is
\[
 Y_t^k=e^{-\xi_k(t-s)}y_k
        +\int_s^t e^{-\xi_k(t-u)}\sqrt{f(u)}\dd B_u^k,
\]
and the scalar component solves the ordinary integral equation
with continuous velocity
$v_t:=\sqrt b\sum_k\sqrt{w_k}Y_t^k$.
Its solution is a measurable function of the initial state and the
future Brownian increments. Uniqueness gives the restart identity,
and independence of those increments from $\mathcal F_s$ proves the
Markov property.

For the Feller property, couple two initial states using the same
increments. The factor differences are deterministic:
\[
 Y_t^k-\widetilde Y_t^k
   =e^{-\xi_k(t-s)}(y_k-\widetilde y_k).
\]
On $[s,T]$, Gronwall's inequality gives
\begin{align*}
 \|X-\widetilde X\|_{\infty;[s,T]}
 \leq{}&\exp\!\left(L_{\drift,T}(T-s)
                   +L_{\sigma,T}\|v\|_{L^1(s,T)}\right)\\
 &\times\left(|x-\widetilde x|
              +S_{\widetilde X,T}
                    \|v-\widetilde v\|_{L^1(s,T)}\right).
\end{align*}
Here $L_{\drift,T}$ and $L_{\sigma,T}$ are the spatial Lipschitz
constants of the coefficients on $\R\times[0,T]$.
The factor $S_{\widetilde X,T}:=
\sup_{t\in[s,T]}|\sigma(\widetilde X_t,t)|$ is finite by continuity.
For a convergent sequence of initial states, take the tilded state as
the limit. Its solution and $S_{\widetilde X,T}$ are then fixed, and the
velocities $v$ converge uniformly on $[s,T]$ under the common-noise coupling.
Thus the full state depends continuously, path by path, on
$(x,\mathbf y)$. Bounded convergence then shows that every transition
operator maps $C_b(\R\times\R^M)$ into itself.
\end{remark}
\subsection{Fully discrete approximation and the joint limit}
\label{subsec:factor-euler}

We now discretize the finite-factor dynamics in time. The error has
three contributions: the spectral approximation, Euler discretization
of the Ornstein--Uhlenbeck factors, and left-point evaluation of the
coefficients in the scalar equation. We keep the dependence on
$\Lambda$ explicit because larger cutoffs introduce faster mean
reversion and require finer time steps. Both partitions may be
nonuniform; the bounds use their maximum cell widths.

Fix $M\in\mathbb N$, $\Lambda>0$, and a deterministic time partition
$\pi=\{0=t_0<t_1<\cdots<t_N=T\}$, with $N\in\mathbb N$. Set
\[
 \Delta_n:=t_{n+1}-t_n,
 \qquad |\pi|:=\max_{0\le n<N}\Delta_n,
\]
and define the left-endpoint map
$\eta_\pi(t):=t_n$ for $t\in[t_n,t_{n+1})$ and
$\eta_\pi(T):=t_{N-1}$. The value at $T$ is immaterial for the
Lebesgue integrals below. For $k=1,\ldots,M$, set
\begin{equation*}
 W_f^k(t):=\int_0^t\sqrt{f(s)}\dd B_s^k.
\end{equation*}
Define
\begin{equation*}
 \Delta W_{f,n}^k
 :=W_f^k(t_{n+1})-W_f^k(t_n),
 \qquad 1\leq k\leq M,\quad 0\leq n<N.
\end{equation*}
The variables $\Delta W_{f,n}^k$ are independent and centered Gaussian,
with variances
\begin{equation*}
 v_{f,n}:=\int_{t_n}^{t_{n+1}}f(s)\dd s.
\end{equation*}
Their simulation requires only the cell integrals $v_{f,n}$, so no
temporal continuity of $f$ is assumed. The exact and discrete factors
are coupled through the same Brownian motions.

The continuous-time Euler approximation of the factors is
\begin{equation*}
 Y_t^{k,M,\Lambda,\pi}
 =-\xi_k\int_0^tY_{\eta_\pi(s)}^{k,M,\Lambda,\pi}\dd s+W_f^k(t),
 \qquad Y_0^{k,M,\Lambda,\pi}=0.
\end{equation*}
At the grid points,
\begin{equation}\label{eq:factor-euler-recursion}
 Y_{t_{n+1}}^{k,M,\Lambda,\pi}
 =(1-\xi_k\Delta_n)Y_{t_n}^{k,M,\Lambda,\pi}
  +\Delta W_{f,n}^k.
\end{equation}
The aggregate factor process and its piecewise linear driver are
\begin{equation*}
 \mathcal Y_t^{M,\Lambda,\pi}
 :=\sqrt b\sum_{k=1}^M\sqrt{w_k}\,Y_t^{k,M,\Lambda,\pi}
\end{equation*}
and
\begin{equation*}
 \zeta_t^{M,\Lambda,\pi}
 :=\int_0^t\mathcal Y_{\eta_\pi(s)}^{M,\Lambda,\pi}\dd s.
\end{equation*}
The Euler approximation of the scalar equation is
\begin{equation}\label{eq:X-factor-euler-integral}
 X_t^{M,\Lambda,\pi}
 =x_0+\int_0^t
   \drift\bigl(X_{\eta_\pi(s)}^{M,\Lambda,\pi},\eta_\pi(s)\bigr)\dd s
 +\int_0^t
   \sigma\bigl(X_{\eta_\pi(s)}^{M,\Lambda,\pi},\eta_\pi(s)\bigr)
   \mathcal Y_{\eta_\pi(s)}^{M,\Lambda,\pi}\dd s.
\end{equation}
Equivalently,
\begin{equation*}
\begin{split}
 X_{t_{n+1}}^{M,\Lambda,\pi}
 ={}&X_{t_n}^{M,\Lambda,\pi}
 +\drift(X_{t_n}^{M,\Lambda,\pi},t_n)\Delta_n\\
 &+\sigma(X_{t_n}^{M,\Lambda,\pi},t_n)
   \mathcal Y_{t_n}^{M,\Lambda,\pi}\Delta_n.
\end{split}
\end{equation*}

The restriction $\Lambda|\pi|\leq1$ ensures
$0\leq1-\xi_k\Delta_n\leq1$ for every factor and every time cell.
Thus all Euler relaxation factors are nonnegative contractions,
uniformly over the spectral nodes. This uniformity is needed when
$\Lambda$ grows with $M$.

\begin{proposition}\label{prop:factor-temporal-error}
Let $T>0$, $0<b<2$, $f\in L^\infty(0,T)$, $M\in\mathbb N$,
$\Lambda>0$, and $p\geq2$.  Suppose that
\begin{equation*}
 \Lambda|\pi|\leq1.
\end{equation*}
There is a constant $C=C(T,b,p,\|f\|_{\infty,T})$, independent of
$M$, $\Lambda$, $\pi$, and the spectral partition and nodes, such that
\begin{equation}\label{eq:factor-euler-pointwise-uniform}
 \max_{1\leq k\leq M}\sup_{t\leq T}
 \left\|Y_t^{k,M,\Lambda}
 -Y_{\eta_\pi(t)}^{k,M,\Lambda,\pi}\right\|_{L^p(\Omega)}
 \leq C|\pi|^{1/2}
\end{equation}
and
\begin{equation}\label{eq:euler-aggregate-explicit}
 \sup_{t\leq T}
 \left\|\mathcal Y_t^{M,\Lambda}
 -\mathcal Y_{\eta_\pi(t)}^{M,\Lambda,\pi}\right\|_{L^p(\Omega)}
 \leq C\Lambda^{(2-b)/2}|\pi|^{1/2}.
\end{equation}
Moreover, for every $\rho\in(0,1)$, there is a constant $C_\rho$ with
the same independence properties such that
\begin{equation}\label{eq:euler-driver-holder-explicit}
 \left(
 \E\norm{\zeta^{M,\Lambda}
 -\zeta^{M,\Lambda,\pi}}_{\Holder{\rho}([0,T])}^p
 \right)^{1/p}
 \leq C_\rho\Lambda^{(2-b)/2}|\pi|^{1/2}.
\end{equation}
\end{proposition}

\begin{proof}

The first two bounds place the supremum outside the probability norm.
We obtain the path-norm bound by integrating the aggregate error in time.

Write $h:=|\pi|$ and fix one factor with node $\xi=\xi_k$.  From
\eqref{eq:factor-euler-recursion},
\[
 Y_{t_n}^{k,M,\Lambda,\pi}
 =\sum_{j=0}^{n-1}
 \left(\prod_{\ell=j+1}^{n-1}(1-\xi\Delta_\ell)\right)
 \Delta W_{f,j}^k.
\]
Since $0\leq1-\xi\Delta_\ell\leq1$, the product is bounded by
$e^{-\xi(t_n-t_{j+1})}$.  Consequently,
\begin{align*}
 \E\abs{Y_{t_n}^{k,M,\Lambda,\pi}}^2
 &\leq\norm{f}_{\infty,T}
 \sum_{j=0}^{n-1}e^{-2\xi(t_n-t_{j+1})}\Delta_j\\
 &\leq\norm{f}_{\infty,T}
 \min\left\{T,\frac{e^2}{2\xi}\right\}.
\end{align*}
Indeed, comparison on each mesh cell and $\xi h\leq1$ give
\[
 \sum_{j=0}^{n-1}e^{-2\xi(t_n-t_{j+1})}\Delta_j
 \leq e^{2\xi h}\int_0^{t_n}e^{-2\xi(t_n-s)}\dd s
 \leq\frac{e^2}{2\xi};
\]
the alternative bound by $T$ is immediate.

For $0\leq t\leq T$, put
\[
 R_t^k:=Y_t^{k,M,\Lambda,\pi}
 -Y_{\eta_\pi(t)}^{k,M,\Lambda,\pi}.
\]
On a cell with left endpoint $t_n=\eta_\pi(t)$, the continuous
interpolation gives
\[
 R_t^k=-\xi(t-t_n)Y_{t_n}^{k,M,\Lambda,\pi}
       +W_f^k(t)-W_f^k(t_n).
\]
The drift contribution has $L^p$ norm at most
$C\xi h\min\{T^{1/2},\xi^{-1/2}\}\leq Ch\sqrt{\xi}
\leq Ch^{1/2}$, and the increment has variance at most
$\|f\|_{\infty,T}h$. Gaussian moment equivalence therefore gives
\begin{equation}\label{eq:euler-local-uniform-new}
 \sup_{t\leq T}\norm{R_t^k}_{L^p(\Omega)}
 \leq C h^{1/2},
\end{equation}
uniformly in $\xi$.

Let $E_t^k:=Y_t^{k,M,\Lambda}-Y_t^{k,M,\Lambda,\pi}$.
Subtracting the two factor equations and solving the resulting linear
equation gives
\[
 E_t^k=-\xi\int_0^t e^{-\xi(t-s)}R_s^k\dd s.
\]
Since $\xi\int_0^t e^{-\xi(t-s)}\dd s\leq1$,
\eqref{eq:euler-local-uniform-new} gives
$\sup_{t\leq T}\norm{E_t^k}_{L^p(\Omega)}\leq Ch^{1/2}$.
Together with
$Y_t^{k,M,\Lambda}-Y_{\eta_\pi(t)}^{k,M,\Lambda,\pi}
=E_t^k+R_t^k$, this proves
\eqref{eq:factor-euler-pointwise-uniform}.

The errors for distinct $k$ are independent centered Gaussian variables.
Thus Gaussian moment equivalence, \eqref{eq:factor-weights}, and
\eqref{eq:factor-euler-pointwise-uniform} give
\[
 \left\|\mathcal Y_t^{M,\Lambda}
 -\mathcal Y_{\eta_\pi(t)}^{M,\Lambda,\pi}\right\|_{L^p(\Omega)}^2
 \leq C_ph\sum_{k=1}^Mw_k
 =C_ph\,\frac{\Lambda^{2-b}}{(2-b)\Gamma(2-b)},
\]
which is \eqref{eq:euler-aggregate-explicit}.

Independence makes the variances add, so the aggregate error carries
the square root of the total spectral mass
$\sum_k w_k=\mu_b((0,\Lambda])$. This accounts for the factor
$\Lambda^{(2-b)/2}$ without introducing a separate loss in $M$.

The driver error vanishes at zero and has almost-everywhere derivative
\[
 D_t:=\mathcal Y_t^{M,\Lambda}
 -\mathcal Y_{\eta_\pi(t)}^{M,\Lambda,\pi}.
\]
Choose $q_0>1/(1-\rho)$ and put $r:=\max\{p,q_0\}$.  H\"older's
inequality in time gives, pathwise,
\[
 \left|\int_s^tD_u\dd u\right|
 \leq |t-s|^{1-1/q_0}\|D\|_{L^{q_0}(0,T)}.
\]
Since $1-1/q_0>\rho$, the full $C^\rho$ norm of the primitive is
bounded by $C_{T,\rho,q_0}\|D\|_{L^{q_0}(0,T)}$.
Monotonicity of probability norms and Minkowski's inequality, with
$r\geq p,q_0$, give
\[
 \norm{D}_{L^p(\Omega;L^{q_0}(0,T))}
 \leq T^{1/q_0}\sup_{t\leq T}\norm{D_t}_{L^r(\Omega)}.
\]
Applying \eqref{eq:euler-aggregate-explicit} at moment order $r$ proves
\eqref{eq:euler-driver-holder-explicit}.
\end{proof}

\begin{assumption}\label{ass:coefficients}
Assume Assumption~\ref{ass:young-quantitative}. For every $T>0$, there
are finite constants $K_T,L_T$ such that
\[
 \begin{aligned}
 |\drift(x,t)|&\leq K_T,\\
 |\drift(x,t)-\drift(y,s)|&\leq L_T(|x-y|+|t-s|),
 \end{aligned}
 \qquad x,y\in\R,\quad s,t\in[0,T].
\]
\end{assumption}

The bounded drift makes the flow estimates uniform in the initial
value; its additional Lipschitz regularity in time controls the
drift-freezing error. The next lemma bounds each local Euler defect
by $\Delta_n^{2\eta}$ times a driver-dependent constant. Summing over
the partition gives the power $|\pi|^{2\eta-1}$, with an exponential
factor of order $1/\eta$ as in \eqref{eq:young-stability}.

\begin{lemma}\label{lem:young-euler-scheme}
Let $T>0$, $\eta\in(1/2,1)$, and suppose that
Assumption~\ref{ass:coefficients} holds.  Let $z\in\Holder{\eta}([0,T])$
with $z(0)=0$, and let $x$ solve
\[
 x_t=x_0+\int_0^t\drift(x_s,s)\dd s
 +\int_0^t\sigma(x_s,s)\dd z_s.
\]
For $t\in[t_n,t_{n+1}]$, define the left-point Euler interpolation by
\begin{equation}\label{eq:abstract-young-euler}
\begin{split}
 x_t^\pi={}&x_{t_n}^\pi
 +\drift(x_{t_n}^\pi,t_n)(t-t_n)\\
 &+\sigma(x_{t_n}^\pi,t_n)(z_t-z_{t_n}),
 \qquad x_0^\pi=x_0.
\end{split}
\end{equation}
Then
\begin{equation}\label{eq:young-euler-pathwise-rate}
 \norm{x-x^\pi}_{\infty;[0,T]}
 \leq C\exp\!\left(
 C\left[1+\norm{z}_{\Holder{\eta}([0,T])}^{1/\eta}\right]\right)
 |\pi|^{2\eta-1}.
\end{equation}
The constant $C$ depends only on $T$, $\eta$, and the coefficient bounds
on $\R\times[0,T]$; it is independent of the partition, the driver,
and the initial value.
\end{lemma}

\begin{proof}
Put $K:=1+[z]_{\eta;[0,T]}$. All constants in this proof depend only
on $T$, $\eta$, and the stated coefficient bounds. Denote by
$\varphi_{s,t}(\xi)$ the solution at time $t$ started from $\xi$ at time
$s$, using the increments of the same path $z$.
The bounded-drift case of Lemma~\ref{lem:young-solution-stability},
applied on subintervals with this common driver, gives
\begin{equation}\label{eq:young-flow-lipschitz}
 \sup_{0\le s\le t\le T}\operatorname{Lip}(\varphi_{s,t})
 \le C\exp(CK^{1/\eta}).
\end{equation}
The time shifts in this application preserve the coefficient bounds,
and the shifted driver has H\"older norm at most $C_TK$.

We first estimate the error of a single Euler step uniformly in its
initial value. The bounded-drift case of the local solution estimate
\eqref{eq:young-local-solution-bound} gives
\[
 [\varphi_{s,\cdot}(\xi)]_{\eta;[s,s+\ell]}\le CK
 \quad\text{if}\quad \ell\le d:=\min\{T,cK^{-1/\eta}\},
\]
where $0<c\le1$ is chosen sufficiently small in terms of the deterministic
bounds. In particular, the displacement on such an interval is at most
$CK\ell^\eta\le C$. Concatenating at most $C_T(1+K^{1/\eta})$
intervals therefore gives
\begin{equation}\label{eq:young-flow-displacement}
 \sup_{s\le t,\,\xi}|\varphi_{s,t}(\xi)-\xi|
 \le C(1+K^{1/\eta}).
\end{equation}

Write
\[
 \Psi_{s,t}(\xi):=\xi+\drift(\xi,s)(t-s)
                          +\sigma(\xi,s)(z_t-z_s),
 \qquad \Delta:=t-s.
\]
For $\Delta\le d$, let $x_r=\varphi_{s,r}(\xi)$. Subtracting the frozen
coefficients from the equation gives the defect
\[
 \varphi_{s,t}(\xi)-\Psi_{s,t}(\xi)
 =\int_s^t[\drift(x_r,r)-\drift(\xi,s)]\dd r
  +\int_s^t[\sigma(x_r,r)-\sigma(\xi,s)]\dd z_r.
\]
Joint Lipschitz continuity of the drift bounds the first term by
$C(K\Delta^{1+\eta}+\Delta^2)$. For the second,
the integrand vanishes at $s$ and has $\eta$-H\"older seminorm at most
$CK$. Young's inequality bounds it by $CK^2\Delta^{2\eta}$.
Since $\eta<1$ and $\Delta\le T$, the two drift contributions are
bounded by $C_TK^2\Delta^{2\eta}$. Thus the local defect has this
same bound for $\Delta\le d$.

If $\Delta>d$, then necessarily $d=cK^{-1/\eta}$. Boundedness of
the coefficients and \eqref{eq:young-flow-displacement} give a defect
at most $C_T(1+K^{1/\eta})$. Also
$\Delta^{2\eta}>c^{2\eta}K^{-2}$, and hence the defect is at most
$C_TK^{2+1/\eta}\Delta^{2\eta}$. Since $K\ge1$, this and the
small-interval bound are both bounded by
$C\exp(CK^{1/\eta})\Delta^{2\eta}$. We have proved
\begin{equation}\label{eq:young-flow-defect}
 \sup_\xi|\varphi_{s,t}(\xi)-\Psi_{s,t}(\xi)|
 \le C\exp(CK^{1/\eta})\Delta^{2\eta},
 \qquad 0\le s<t\le T.
\end{equation}
Thus no restriction on the temporal mesh is needed for this pathwise
bound.

For a fixed mesh point $t_m$, set
$F_k:=\varphi_{t_k,t_m}(x_{t_k}^\pi)$ for $k=0,\ldots,m$.
Then $F_0=x_{t_m}$ and $F_m=x_{t_m}^\pi$.
The flow composition rule and the Euler recursion give
\begin{align*}
 |F_k-F_{k+1}|
 &\le \operatorname{Lip}(\varphi_{t_{k+1},t_m})
 \left|\varphi_{t_k,t_{k+1}}(x_{t_k}^\pi)
              -\Psi_{t_k,t_{k+1}}(x_{t_k}^\pi)\right|.
\end{align*}
Summing these differences and using
\eqref{eq:young-flow-lipschitz} and
\eqref{eq:young-flow-defect} gives
\begin{equation}\label{eq:lejay-euler-applied}
 \max_m|x_{t_m}-x_{t_m}^\pi|
 \le C\exp(CK^{1/\eta})\sum_k\Delta_k^{2\eta}
 \le CT\exp(CK^{1/\eta})|\pi|^{2\eta-1}.
\end{equation}
This telescoping uses exact flows after each local error, so it does
not require a H\"older estimate for the numerical path.

To include times between the nodes, the local estimate and
\eqref{eq:young-flow-displacement} give
$[x]_{\eta;[0,T]}\le C(1+K^{1+1/\eta})$: use the local bound for
separations at most $d$, and divide the global displacement bound by
$d^\eta$ for larger separations. Bounded coefficients in
\eqref{eq:abstract-young-euler} also give
\[
 |x_t^\pi-x_{t_n}^\pi|\le C_TK|\pi|^\eta,
 \qquad t_n\le t\le t_{n+1}.
\]
Hence the uniform error is bounded by the nodal error plus
$C(1+K^{1+1/\eta})|\pi|^\eta$. Since
$\eta>2\eta-1$ and $|\pi|\le T$, this last term is absorbed in
\eqref{eq:young-euler-pathwise-rate}; all polynomial factors in $K$
are absorbed into its exponential.
\end{proof}

We combine the spectral and temporal estimates next. In the
finite-variation regime we use Gronwall's inequality; in the Young
regime we use Lemmas~\ref{lem:young-solution-stability}
and~\ref{lem:young-euler-scheme}. The three terms in
\eqref{eq:joint-young} bound the spectral, factor discretization,
and coefficient-freezing errors, respectively.

\begin{theorem}\label{thm:joint-fully-discrete}
Let $T>0$, $p\geq2$, and $f\in L^\infty(0,T)$.
Assume that the spectral partition satisfies
$\delta_\lambda\leq C_{\mathrm{mesh}}\Lambda/M$, with fixed
$C_{\mathrm{mesh}}\geq1$ and arbitrary nodes $\xi_k\in I_k$.
Couple the original process and its approximations through $\mathcal W$
and use the same deterministic initial value $x_0$.
For $0<b\le1$, impose Assumption~\ref{ass:coefficients};
for $1<b<2$, assume that both coefficients are bounded and jointly
Lipschitz on $\R\times[0,T]$.

For $1<b<2$, set $\Lambda=M^{2/3}$, and suppose that
$M^{2/3}|\pi|\leq1$.  Then
\begin{align}
 &\left(\E\sup_{t\leq T}
 \abs{X_t-X_t^{M,\Lambda,\pi}}^p\right)^{1/p}
 \notag\\*
 &\qquad\leq
 C M^{-(b-1)/3}
 +C M^{(2-b)/3}|\pi|^{1/2}
 \leq C M^{-(b-1)/3}.
 \label{eq:joint-finite-variation}
\end{align}

For $0<b\leq1$, set $\Lambda=M^{1/2}$, and choose
\[
 \frac12<\eta<\eta'<H_b^\circ.
\]
If $M^{(2-b)/2}|\pi|\leq1$, then
\begin{align}
 &\left(\E\sup_{t\leq T}
 \abs{X_t-X_t^{M,\Lambda,\pi}}^p\right)^{1/p}
 \notag\\*
 &\quad\leq
 C\left(M^{-b/4}\sqrt{\log(2+M)}\right)^{1-\eta/\eta'}
 +C M^{(2-b)/4}|\pi|^{1/2}
 +C|\pi|^{2\eta-1}.
 \label{eq:joint-young}
\end{align}
All constants in \eqref{eq:joint-finite-variation}--\eqref{eq:joint-young}
may depend on $C_{\mathrm{mesh}}$, but are independent of $M$, $\pi$,
and the admissible spectral partitions and nodes.
\end{theorem}

\begin{proof}
In the finite-variation case, abbreviate
\[
 Z:=\mathcal Y^{M,\Lambda},
 \qquad
 Z^\pi:=\mathcal Y_{\eta_\pi(\cdot)}^{M,\Lambda,\pi}.
\]
Proposition~\ref{prop:factor-temporal-error} gives
\begin{equation}\label{eq:finite-joint-aggregate}
 \norm{Z-Z^\pi}_{L^p(\Omega;L^1(0,T))}
 \leq C M^{(2-b)/3}|\pi|^{1/2}.
\end{equation}
Theorem~\ref{thm:finite-variation-rate} and the triangle inequality imply
that the centered Gaussian random elements $Z$ are uniformly bounded in
$L^2(\Omega;L^2(0,T))$.  Integrating
\eqref{eq:euler-aggregate-explicit} with $p=2$ over time and using
$M^{2/3}|\pi|\leq1$ gives the same bound for $Z-Z^\pi$, and hence for
$Z^\pi$. For $U\in\{Z,Z^\pi\}$, let $Q_{M,\pi}$ be its covariance
operator on $L^2(0,T)$. Then
\[
 \|Q_{M,\pi}\|_{\mathrm{op}}
 \leq\Tr Q_{M,\pi}=\E\|U\|_{L^2(0,T)}^2\leq C_T
\]
uniformly in $M$ and $\pi$.  If
$0<\varepsilon_F<(4C_T)^{-1}$ and $(\lambda_j)$ are the eigenvalues of
$Q_{M,\pi}$, then
\[
 \E e^{\varepsilon_F\|U\|_{L^2(0,T)}^2}
 =\prod_j(1-2\varepsilon_F\lambda_j)^{-1/2}
 \leq\exp\!\left(
   \frac{\varepsilon_F\Tr Q_{M,\pi}}
        {1-2\varepsilon_F\|Q_{M,\pi}\|_{\mathrm{op}}}
 \right)\leq e^{2\varepsilon_F C_T}.
\]
The product formula follows from the orthonormal eigenbasis expansion
and monotone convergence. Thus $Z$ and $Z^\pi$ have uniformly bounded
moments of every order in $L^2(0,T)$. Since
$\|U\|_{L^1}\leq\sqrt{T}\|U\|_{L^2}$, their $L^1(0,T)$ norms
also have uniform exponential moments of every linear order.

Let
\[
 D_t:=\sup_{r\leq t}
 \abs{X_r^{M,\Lambda}-X_r^{M,\Lambda,\pi}}.
\]
On every mesh cell,
\[
 \abs{X_s^{M,\Lambda,\pi}
 -X_{\eta_\pi(s)}^{M,\Lambda,\pi}}
 \leq C|\pi|(1+\abs{Z_s^\pi}).
\]
Subtracting the two scalar equations and using the joint Lipschitz
bounds in space and time gives
\begin{align*}
 D_t\leq{}&
 C\int_0^t(1+\abs{Z_s})D_s\dd s\\
 &+C|\pi|\int_0^T(1+\abs{Z_s})(1+\abs{Z_s^\pi})\dd s
 +C\int_0^T\abs{Z_s-Z_s^\pi}\dd s.
\end{align*}
The local error term satisfies
\[
 \int_0^T(1+|Z_s|)(1+|Z_s^\pi|)\dd s
 \leq T+\sqrt{T}(\|Z\|_{L^2}+\|Z^\pi\|_{L^2})
       +\|Z\|_{L^2}\|Z^\pi\|_{L^2},
\]
whose right-hand side has uniformly bounded moments of every order.
Gronwall's inequality, the exponential moment bounds, H\"older's
inequality, and \eqref{eq:finite-joint-aggregate} at a higher moment
order therefore give
\[
 \left(\E\sup_{t\leq T}
 \abs{X_t^{M,\Lambda}-X_t^{M,\Lambda,\pi}}^p\right)^{1/p}
 \leq C M^{(2-b)/3}|\pi|^{1/2}.
\]
Here the term of order $|\pi|$ is absorbed using
$|\pi|\leq|\pi|^{1/2}\leq M^{(2-b)/3}|\pi|^{1/2}$.
Combining this bound with \eqref{eq:finite-variation-X-rate} proves
the first inequality in \eqref{eq:joint-finite-variation}; the mesh
condition gives the second.

In the Young case, write
\[
 g^{M,\pi}:=\zeta^{M,\Lambda}-\zeta^{M,\Lambda,\pi},
 \qquad
 a_{M,\pi}:=M^{(2-b)/4}|\pi|^{1/2}.
\]
The mesh condition $M^{(2-b)/2}|\pi|\leq1$ implies
$\Lambda|\pi|\leq1$ and $a_{M,\pi}\leq1$, since $0<b\leq1$.
Proposition~\ref{prop:factor-temporal-error} therefore gives
\begin{equation}\label{eq:young-joint-driver}
 \left(\E\norm{g^{M,\pi}}_{\Holder{\eta}([0,T])}^p\right)^{1/p}
 \leq Ca_{M,\pi}.
\end{equation}
Moreover, Minkowski's inequality applied to the time integral of
$\mathcal Y^{M,\Lambda}
-\mathcal Y_{\eta_\pi(\cdot)}^{M,\Lambda,\pi}$ and
\eqref{eq:euler-aggregate-explicit} at $p=2$ yield
\[
 \norm{g_t^{M,\pi}-g_s^{M,\pi}}_{L^2(\Omega)}
 \leq Ca_{M,\pi}|t-s|.
\]
Fix $H_0\in(\eta,1)$. Since $a_{M,\pi}\leq1$ and
$|t-s|\leq T^{1-H_0}|t-s|^{H_0}$,
Lemma~\ref{lem:uniform-gaussian-holder} gives an exponential moment
bound for $\|g^{M,\pi}\|_{\Holder{\eta}([0,T])}^2$, uniform in
$M$ and $\pi$. Combining it with
\eqref{eq:uniform-young-exponential-integrability} and H\"older's
inequality gives, for some $c_0>0$,
\[
 \sup_{M,\pi}\E\exp\!\left(
 c_0\|\zeta^{M,\Lambda,\pi}\|_{\Holder{\eta}([0,T])}^2
 \right)<\infty.
\]
Since $1/\eta<2$, it follows that
\[
 \sup_{M,\pi}
 \E\exp\!\left(
 c\|\zeta^{M,\Lambda,\pi}\|_{\Holder{\eta}([0,T])}^{1/\eta}
 \right)<\infty
 \qquad\text{for every fixed }c>0,
\]
where both suprema range over the admissible meshes. The same bounds
hold for the exact finite-factor drivers.

Let $\widetilde X^{M,\Lambda,\pi}$ solve the Young equation driven by
$\zeta^{M,\Lambda,\pi}$ with initial value $x_0$. This intermediate
solution isolates the error in the driver from the Euler error in
the scalar equation:
\begin{align*}
 \|X-X^{M,\Lambda,\pi}\|_\infty
 \leq{}&\|X-X^{M,\Lambda}\|_\infty\\
 &+\|X^{M,\Lambda}-\widetilde X^{M,\Lambda,\pi}\|_\infty
 +\|\widetilde X^{M,\Lambda,\pi}-X^{M,\Lambda,\pi}\|_\infty.
\end{align*}
The first term is the spectral error already estimated in
\eqref{eq:young-X-rate}. We estimate the other two terms separately.

Lemma~\ref{lem:young-solution-stability}, the preceding exponential
bounds, \eqref{eq:young-joint-driver} at moment order $2p$, and
H\"older's inequality give
\[
 \left(\E\sup_{t\leq T}
 \abs{X_t^{M,\Lambda}-\widetilde X_t^{M,\Lambda,\pi}}^p
 \right)^{1/p}
 \leq Ca_{M,\pi}.
\]
Since $\zeta^{M,\Lambda,\pi}$ is piecewise linear and has derivative
$Z^\pi:=\mathcal Y_{\eta_\pi(\cdot)}^{M,\Lambda,\pi}$, the scheme
\eqref{eq:X-factor-euler-integral} is exactly the Euler
interpolation \eqref{eq:abstract-young-euler} for this driver.  Hence
Lemma~\ref{lem:young-euler-scheme} and the uniform exponential bound imply
\[
 \left(\E\sup_{t\leq T}
 \abs{\widetilde X_t^{M,\Lambda,\pi}
 -X_t^{M,\Lambda,\pi}}^p\right)^{1/p}
 \leq C|\pi|^{2\eta-1}.
\]
The triangle inequality and \eqref{eq:young-X-rate} now prove
\eqref{eq:joint-young}.
\end{proof}

For $1<b<2$, the condition $|\pi_M|\leq M^{-2/3}$ already makes
the temporal error no larger than the spectral error. For
$0<b\leq1$, the mesh restriction in the theorem ensures uniform
stability. To make the temporal terms on the right-hand side
of \eqref{eq:joint-young} tend to zero, one chooses
$M^{(2-b)/2}|\pi_M|\to0$. The following corollary gives time-step
bounds under which neither temporal error exceeds the order of the
spectral error. The two entries in the maximum
defining $\kappa_{b,\eta,\eta'}$ correspond to the factor error
and the coefficient-freezing error, respectively.

\begin{corollary}\label{cor:joint-meshes}
Under the assumptions of Theorem~\ref{thm:joint-fully-discrete}, let
$\pi=\pi_M$.

If $1<b<2$ and $|\pi_M|\leq M^{-2/3}$, then
\[
 \left(\E\sup_{t\leq T}
 \abs{X_t-X_t^{M,M^{2/3},\pi_M}}^p\right)^{1/p}
 \leq C M^{-(b-1)/3}.
\]

If $0<b\leq1$ and $1/2<\eta<\eta'<H_b^\circ$, define
\[
 r_{b,\eta,\eta'}
 :=\frac b4\left(1-\frac{\eta}{\eta'}\right)
\]
and
\begin{equation}\label{eq:young-joint-kappa}
 \kappa_{b,\eta,\eta'}
 :=\max\left\{
 \frac{2-b}{2}+2r_{b,\eta,\eta'},
 \frac{r_{b,\eta,\eta'}}{2\eta-1}
 \right\}.
\end{equation}
If $|\pi_M|\leq M^{-\kappa_{b,\eta,\eta'}}$, then
\[
 \left(\E\sup_{t\leq T}
 \abs{X_t-X_t^{M,M^{1/2},\pi_M}}^p\right)^{1/p}
 \leq
 C M^{-r_{b,\eta,\eta'}}
 \bigl(\log(2+M)\bigr)^{(1-\eta/\eta')/2}.
\]
These mesh choices preserve the spectral convergence rate.
The choice
$b=1$, $\eta=0.7$, and $\eta'=0.9$ gives
\[
 r_{1,0.7,0.9}=\frac1{18},
 \qquad
 \kappa_{1,0.7,0.9}=\frac{11}{18}.
\]
\end{corollary}

\begin{proof}
The finite-variation statement is \eqref{eq:joint-finite-variation}.  In the
Young regime, \eqref{eq:young-joint-kappa} gives
\[
 \frac{\kappa_{b,\eta,\eta'}}2-\frac{2-b}{4}
 \geq r_{b,\eta,\eta'},
 \qquad
 \kappa_{b,\eta,\eta'}(2\eta-1)
 \geq r_{b,\eta,\eta'}.
\]
The first inequality also implies the mesh restriction in
Theorem~\ref{thm:joint-fully-discrete}. Substitution in
\eqref{eq:joint-young} proves the claim.
\end{proof}

\begin{remark}\label{rem:joint-complexity}

The convergence estimates depend only on the maximal temporal cell
width. To relate the computational cost to that width, we also
bound the number of time steps.

Call the temporal meshes quasi-uniform if there is a fixed $c>0$ such
that $c|\pi_M|\leq\Delta_n\leq|\pi_M|$ for every cell and every $M$.
Let $N_M$ denote the number of cells of $\pi_M$. For positive sequences,
$a_M\asymp b_M$ means that $cb_M\leq a_M\leq Cb_M$ for all sufficiently
large $M$, with constants $0<c\leq C<\infty$ independent of $M$.
For such meshes, the finite-variation choice
$|\pi_M|\asymp M^{-2/3}$ uses $N_M\asymp M^{2/3}$ time steps.  Updating all
$M$ factors at every step therefore has arithmetic cost
$O(MN_M)=O(M^{5/3})$.  In the Young regime, the choice
$|\pi_M|\asymp M^{-\kappa_{b,\eta,\eta'}}$ uses
$N_M\asymp M^{\kappa_{b,\eta,\eta'}}$ steps and has cost
$O(M^{1+\kappa_{b,\eta,\eta'}})$. These counts apply to spectral
partitions satisfying the maximum-mesh condition, once the cell
variances $v_{f,n}$ are available. Adapted quadratures or exact
Ornstein--Uhlenbeck transitions may change the stability restriction and the
cost balance.

Taking $\eta\downarrow1/2$ with $\eta'$ fixed improves the spectral
exponent $r_{b,\eta,\eta'}$, but makes
$r_{b,\eta,\eta'}/(2\eta-1)$ diverge. Retaining the spectral rate
therefore requires increasingly fine temporal meshes; that rate
alone does not determine the accuracy per arithmetic operation.
The estimates above are stated for $\eta,\eta'$ fixed independently
of $M$. Their constants may depend on these exponents, so choosing
$\eta=\eta(M)$ requires additional control of that dependence.
 \end{remark}
\section*{Statements and Declarations}

\bmhead{Funding}
J.~H.~Ram\'irez-Gonz\'alez acknowledges support from the S\~ao Paulo
Research Foundation (FAPESP), grant 2025/03804-0.

\bmhead{Competing interests}
The authors have no relevant financial or non-financial interests to disclose.

\bmhead{Data availability}
Data sharing is not applicable to this article as no datasets were
generated or analysed during the current study.

\bmhead{Use of AI assistance}
AI tools were used to improve the wording and flow of the manuscript.

\end{document}